%% file: ex_article.tex
\documentclass[onefignum,onetabnum]{siamart220329}

\usepackage{latexsym,amssymb,amsmath,mathtools,epsfig,xcolor} 
\usepackage{url}
\usepackage{subfig}
\usepackage{ragged2e}
\usepackage[]{graphicx} 
\usepackage{bm}
\usepackage{mathrsfs}
\usepackage{float}
\usepackage{hyperref}
\usepackage{cite}
\usepackage{appendix}
\usepackage[font=small,labelfont=bf]{caption}
\usepackage{pifont,enumitem}
\usepackage{setspace}

\usepackage{algorithm,algorithmic}

\usepackage{booktabs,array,multirow}    
\usepackage{tabularx}
\usepackage{longtable}

\usepackage[usestackEOL]{stackengine}   

\DeclareMathOperator*{\minimize}{minimize}
\DeclareMathOperator{\atantwo}{atan2}

\newcommand{\vb}{{\boldsymbol b}}

\newcommand{\mA}{{\boldsymbol A}}
\newcommand{\mR}{{\boldsymbol R}}
\newcommand{\mI}{{\boldsymbol I}}
\newcommand{\mP}{{\boldsymbol P}}
\newcommand{\mB}{{\boldsymbol B}}
\newcommand{\mD}{{\boldsymbol D}}
\newcommand{\mO}{{\boldsymbol O}}
\newcommand{\mF}{{\boldsymbol F}}
\newcommand{\mE}{{\boldsymbol E}}

\newcommand{\mQ}{{\boldsymbol Q}}
\newcommand{\mU}{{\boldsymbol U}}

\newcommand{\mV}{{\boldsymbol V}}

\newcommand{\mC}{{\boldsymbol C}}
\newcommand{\vw}{{\boldsymbol w}}
\newcommand{\va}{{\boldsymbol a}}

\newcommand{\vf}{{\boldsymbol f}}
\newcommand{\vg}{{\boldsymbol g}}

\newcommand{\vk}{{\boldsymbol k}}
\newcommand{\vr}{{\boldsymbol r}}

\newcommand{\vp}{{\boldsymbol p}}

\input{ex_shared}

\ifpdf
\hypersetup{
  pdftitle={Orthogonal Matrix Retrieval with Spatial Consensus for 3D Unknown-View Tomography},
  pdfauthor={Shuai Huang, Mona Zehni, Ivan Dokmani\'c, and Zhizhen Zhao}
}
\fi

\begin{document}

\maketitle

\begin{abstract}
Unknown-view tomography (UVT) reconstructs a 3D density map from its 2D projections at unknown, random orientations. A line of work starting with Kam (1980) employs the method of moments (MoM) with rotation-invariant Fourier features to solve UVT in the frequency domain, assuming that the orientations are uniformly distributed. This line of work includes the recent orthogonal matrix retrieval (OMR) approaches based on matrix factorization, which, while elegant, either require side information about the density that is not available, or fail to be sufficiently robust. 
For OMR to break free from those restrictions, we propose to jointly recover the density map and the orthogonal matrices by requiring that they be mutually consistent. 
We regularize the resulting non-convex optimization problem by a denoised reference projection and a nonnegativity constraint. This is enabled by the new closed-form expressions for spatial autocorrelation features. Further, we design an easy-to-compute initial density map which effectively mitigates the non-convexity of the reconstruction problem. Experimental results show that the proposed OMR with spatial consensus is more robust and performs significantly better than the previous state-of-the-art OMR approach in the typical low-SNR scenario of 3D UVT.
\end{abstract}



\section{Introduction}
\label{sec:intro}
Unknown view tomography (UVT) arises in applications such as single-particle cryo-electron microscopy (cryo-EM), where noisy projections of biological macromolecules are taken at \emph{random, unknown} view angles and then used to reconstruct the 3D molecular density map \cite{CryoEM:Primer:2015}. The unknown particle orientations and the low signal-to-noise-ratio (SNR) of projection images make this reconstruction a challenging task. Depending on whether the particle orientations need to be estimated, there are generally two approaches to reconstruction. The first approach proceeds by alternating between reconstructing the density map and estimating the particle orientations. Since the projection images typically contain more noise than signal ($\text{SNR} < 1$), it is hard to estimate the particle orientations accurately. Soft assignments are thus adopted in the form of posterior distributions over view angles. The density map is computed as the maximum-a-posteriori (MAP) estimate via the expectation maximization algorithm \cite{Dempster:EM:1977,Sigworth:MLE_CryoEM:1998,Scheres:BayesianCryoEM:2012,Scheres:Relion:2012}. Since all projection images need to be matched to reference templates at each iteration, this approach is computationally expensive.

The second approach bypasses estimation of particle orientations and recovers the density map through \textit{autocorrelation analysis}, an instance of the method of moments (MoM) \cite{KAM1977, KAM1980, Kam1985}. For uniformly distributed orientations, it is known that the MoM achieves generic list recovery---determining the density map up to a finite list of candidate densities---from first-, second-, and third-order moments \cite{bandeira2018estimation}. Autocorrelation analysis, however, primarily uses the first- and second-order moments. Levin et al. \cite{Eithan2017} have shown that using second-order moments suffices for unique recovery if they are complemented by two projections with known view angles. For a known, non-uniform rotation distribution, Sharon et al. \cite{Sharon_2020} showed that the first- and second-moments suffice to determine a finite-list of possible structures. Higher-order moments come with a price: while it is possible to use the third or even fourth-order moments, the computational complexity scales exponentially with the moment-cutoff order,
and the noise amplification in higher-order moments drastically increases the number of required projection images \cite{CryoEM:Singer:2018}. For uniformly distributed orientations, the sample complexity of the MoM with moments of order up to $m$ scales at a rate of 1/SNR$^m$ \cite{bandeira2018estimation}. In general, the moments are calculated with respect to single-particle projection images cropped from micrographs. Bendory et al. \cite{Bendory18} also showed that it is possible to recover the density map through autocorrelation analysis on the micrographs directly.

Kam proposed an MoM approach with Fourier autocorrelation functions which are rotation-invariant second-order moments \cite{KAM1980}. The same moments were used for reconstruction in a number of follow-up works \cite{Bhamre:OMR:2015,Eithan2017,Sharon_2020}. Fourier autocorrelations can be estimated from the projection images assuming that particle orientations are uniformly distributed. They determine the spherical harmonic expansion coefficients of the Fourier transform of the density map up to a set of \emph{unknown orthogonal matrices}. They are useful in the low SNR regime where it is difficult to get accurate orientation estimations. 

Our proposed approach builds on the line of work spearheaded by Kam and leverages the method of moments with rotation-invariant features \cite{KAM1980, Kam1985, Bhamre:OMR:2015, Eithan2017, Sharon_2020, Singer:WilsonStat:2021}. 
Kam's spherical harmonic expansion coefficients can be either solved for directly \cite{Kam1985,Sharon_2020}, or recovered by solving an orthogonal matrix retrieval (OMR) problem \cite{Bhamre:OMR:2015,Eithan2017}. Existing OMR methods, however, require additional information that is usually not available. For example, Bhamre et al. \cite{Bhamre:OMR:2015} describe two algorithms: the orthogonal extension method which requires the structure of a similar molecule and the orthogonal replacement method which requires projection images from two unknown structures and assumes the differences between the two structures are known. Levin et al. \cite{Eithan2017} describe an improved OMR by projection matching (OMR-PM) which requires (at least) two denoised projection images to perform the reconstruction. However, it is not robust enough to handle complicated density maps (cf. Section \ref{sec:exp}).

On the other hand, earlier work on unassigned distance geometry and unknown view tomography shows that the density map can be directly optimized in the spatial domain \cite{RPSFDD:Huang:2021,2DUVT:Zehni:2019,3DUVT:Zehni:2020}. That earlier work, however, employs a parametric point-source density map which does not scale easily to realistic molecular density maps.

\subsection{Our Contributions and Paper Outline}
We formulate the reconstruction in terms of both the orthogonal matrices and the density map for a consistent recovery.
\begin{itemize}
\item We propose novel radial and autocorrelation features in the spatial domain. Compared to Fourier autocorrelations, the proposed spatial autocorrelations have simpler closed-form expressions in terms of the density map.
\item Previous OMR approaches recover only the orthogonal matrices, and it has been difficult to incorporate spatial information through the orthogonal matrices.  
We relate the orthogonal matrices to the density map and recover them jointly via alternating optimization, so that the spatial consensus on the density can be enforced among the orthogonal matrices.
\item To make this work, we construct an initial density by solving a convex optimization program that involves the spatial radial features and a denoised reference projection image. This initialization provides the basis for the orthogonal matrices to ``reach a consensus'' on a density that additionally satisfies the nonnegativity constraints and matches the reference projection.

\end{itemize}
Experimental results show that the proposed orthogonal matrix retrieval with spatial consensus (OMR-SC) is more robust across a range of different density maps than the previous state-of-the-art OMR-PM, and that it excels in the low-SNR regime that is common in unknown-view tomography. Reproducible code and data are available at `` \urlstyle{tt}\url{https://github.com/shuai-huang/OMR-SC} ''.

\paragraph{Paper Outline} In Section \ref{sec:uvt_mom}, we first set up the mathematical model for unknown view tomography and review Kam's autocorrelation analysis in the Frequency domain. The proposed spatial radial and autocorrelation features and their relations to Fourier features are then introduced in Section \ref{sec:spatial_mom_features}, and the density map is parameterized in Section \ref{sec:parametric_density_map}. Building on the parametric radial and autocorrelation features, we propose the OMR-SC approach in Section \ref{sec:omr}. We then compare the OMR-SC and OMR-PM approaches on the recovery of random and protein density maps in Section \ref{sec:exp}, and conclude the paper with a discussion in Section \ref{sec:conclusion}. Additional experimental results are given in the Supplementary Material.

\paragraph{Notation} We use non-boldface lower- and upper-case letters to represent scalars, boldface lower-case letters to represent vectors, and boldface upper-case letters to represent matrices. For readers' convenience, we provide a list of notations of the variables in Table SM1 in the Supplementary Material.

\section{Unknown View Tomography via the Method of Moments}
\label{sec:uvt_mom}

\subsection{Problem Formulation}
\label{subsec:problem_formulation}

\begin{figure}[tbp]
\centering
\includegraphics[width=\textwidth]{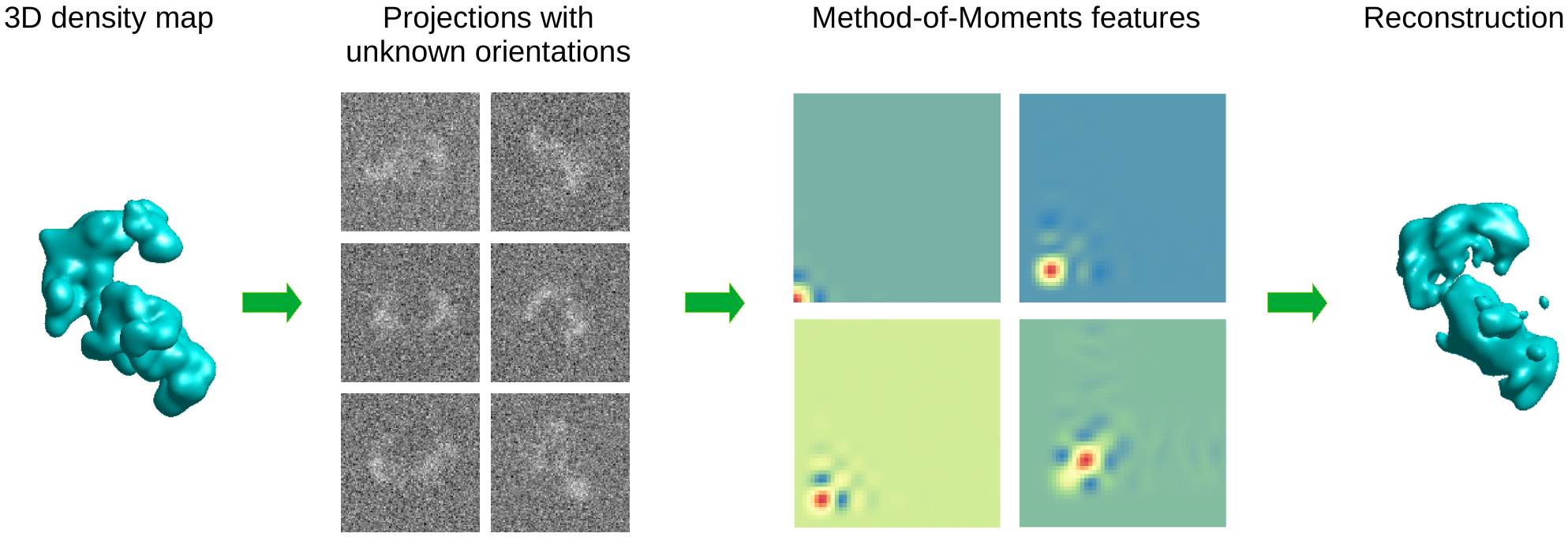}
\caption{The unknown view tomography via the Method-of-Moments.}
\label{fig:uvt_mom}
\end{figure}

Let $\vr \mapsto \rho(\vr)$ denote the 3D density map to be estimated, where $\vr = [x\ y\ z]^T \in \mathbb{R}^3$ contains the Cartesian coordinates. We assume $\rho(\vr)$ is compactly supported within a ball, 
approximately bandlimited with an effective bandwidth $\pi$, and square-integrable.\footnote{Since $\text{supp}(\rho)$ is compact, square integrability implies integrability.}

As illustrated in Fig. \ref{fig:uvt_mom}, our input data consists of $N$ projection images, each containing a projection of $\rho$ at some unknown orientation. Let $\mR_n$ denote the $3\times 3$ rotation matrix representing a 3D rotation $\chi_n \in \mathrm{SO}(3)$. The rotated density model is then $\rho(\mR_n^T\vr)$. The noiseless projection $P_n(x,y)$ along the $z$-direction is
\begin{align}
    \textnormal{Noiseless projection:}\quad P_n(x,y)=\int_{-\infty}^{\infty}\rho\left(\mR_n^T[x\ y\ z]^T\right)\ dz\,.
\end{align}
It is further corrupted by additive noise $\varepsilon_n(x,y)$ to produce the observed noisy projection image $S_n(x,y)$,
\begin{align}
    \label{eq:proj_model}
    \textnormal{Noisy projection:}\quad S_n(x,y) = P_n(x, y) + \varepsilon_n(x,y)\,.
\end{align}
The noise $\varepsilon_n$ is modeled as white Gaussian noise with variance estimated from the measured data \cite{CryoEM:Bendory:2020}. In practice, images of individual particles are cropped from the micrographs. Additional in-plane translation misalignment may occur when the particle is not centered during particle picking; a number of methods are available to mitigate it \cite{WANG2016325,Zhu2017ADC,Wagner:PP:2019,Bepler:PP:2019,Heimowitz:2021}. In this paper, we extract rotation-invariant features from the projection images and assume that the particles have been properly centered in the projection images. Although the distribution of particle orientations is not strictly uniform in practice, the uniform assumption is widely adopted by MoM approaches and has enjoyed empirical successes~\cite{KAM1980,Bhamre:OMR:2015,Eithan2017}. Here we also assume that the unknown rotations $\chi_n\in \mathrm{SO}(3)$ are uniformly distributed in $\mathrm{SO}(3)$. Our goal is then to reconstruct the 3D density map $\rho(\vr)$ from a collection of 2D noisy projection images $\{S_n\ |\ n=1,\ldots,N\}$ with unknown, uniformly distributed view directions.

\subsection{Kam's Autocorrelation Analysis in the Frequency Domain}
\label{subsec:kam_autocorrelaton}
We next review the recent developments in Kam's autocorrelation analysis and discuss the limitations of previous OMR approaches in this subsection. The MoM approach proposed by Kam performs the reconstruction in the frequency domain using autocorrelation features extracted from 2D projections \cite{KAM1977,KAM1980}. The Fourier transform $\widehat{\rho}(\vk)$ of the 3D density map $\rho(\vr)$ is given by
\begin{align}
\label{eq:ft}
\widehat{\rho}(\vk)
=\iiint e^{-\mathfrak{i} \vk \cdot \vr}\rho(\vr)\ d\vr =\iiint e^{-\mathfrak{i}kr\cos \alpha}\rho(\vr)\ d\vr\,,
\end{align}
where $\mathfrak{i}$ is the imaginary unit, $\vk=[k_x\ k_y\ k_z]^T\in\mathbb{R}^3$ is the frequency, $k=\|\vk\|$ is the norm of $\vk$, $r=\|\vr\|$ the norm of $\vr$, and $\alpha$ the angle between $\vk$ and $\vr$. We shall write $\va = a \ddot{\va}$ with $a = \|\va\|$ and $\|\ddot{\va}\| = 1$ for the polar representation of a generic vector $\va\in\mathbb{R}^3$.

The Fourier autocorrelation features used by MoM are expressed in terms of the spherical harmonic expansion coefficients $\{A_{lm}(k)\}_{lm}$ of $\widehat{\rho}(\vk)$. We have
\begin{align}
\label{eq:ft_sh_expansion}
    \widehat{\rho}(\vk)=\sum_{l=0}^\infty\sum_{m=-l}^lA_{lm}(k)\cdot Y_{lm}(\theta_\vk,\varphi_\vk)\,,
\end{align}
where $(\theta_\vk,\varphi_\vk)$ is the angular direction of $\vk$, and $Y_{lm}(\cdot)$ is the \emph{real} spherical harmonic function of degree $l$ and order $m$. Note that $A_{lm}(k)$ is purely real for even $l$ and purely imaginary for odd $l$. We have
\begin{align}
\label{eq:A_lm_fourier_eval}
    A_{lm}(k) = \iiint \widehat{\rho}(\tilde{\vk})\cdot\delta(k-\tilde{k})\cdot Y_{lm}(\theta_{\tilde{\vk}},\varphi_{\tilde{\vk}})\cdot \sin\theta_{\tilde{\vk}}\ d\tilde{k} d\theta_{\tilde{\vk}} d\varphi_{\tilde{\vk}}\,,
\end{align}
where $\delta(\cdot)$ is the Dirac impulse.

Features are computed from the 2D projection images. The 2D Fourier transform $\widehat{S}_n$ of a projection $S_n$ is
\begin{align}
    \widehat{S}_n(k_x,k_y)=\iint \exp\left(-\mathfrak{i}\cdot[k_x\ k_y]\left[\begin{array}{c} x\\ y\end{array}\right]\right)\cdot S_n(x,y)\ dx\,dy.
\end{align}
From the central slice theorem, the above $\widehat{S}_n$ corresponds to a central slice of $\widehat{\rho}(\vk)$ \cite{Messiah:1961},
\begin{align}
\label{eq:CST}
\begin{split}
\widehat{S}_n(k_x, k_y) &=  \widehat{\rho}\left(\mR_n^T [k_x\ k_y\ 0]^T\right)  + \widehat{\epsilon}(k_x, k_y) \\
&= \sum_{l=0}^\infty \sum_{m=-l}^l \sum_{m' = -l}^l A_{lm}(k)\cdot Y_{lm'} \left(\frac{\pi}{2}, \varphi_\vk \right) \cdot D^l_{mm'}(\chi_n) + \widehat{\epsilon}(k_x, k_y) \,,
\end{split}
\end{align}
where $\varphi_\vk = \atantwo(k_y,k_x)$, and $D_{mm'}^l(\chi_n)$ is an element of the Wigner D-matrix.
However, since the rotation $\chi_n$ (that is to say, $\mR_n$) is unknown, 
we do not know a priori which central slice $\widehat{S}_n$ corresponds to. What we do know is that the relative position of two frequencies in the Fourier slice $\widehat{S}_n$ is the same as the relative position of their ``true'' 3D counterparts since $[k_x, k_y, 0] [k'_x, k'_y, 0]^T = [k_x, k_y, 0] \mR \mR^T [k'_x, k'_y, 0]^T$ for any rotation $\mR$. (In particular, this allows us to compute the correct length of a frequency vector in the slice.) 

Let $\vk_1$ and $\vk_2$ denote two frequency vectors in the same 2D Fourier slice $\widehat{S}_n$, $\varphi_{\vk_1}=\measuredangle \vk_1$ be the azimuth angle of $\vk_1$, 
and $\psi$ be the angle between $\vk_1$ and $\vk_2$ such that $\varphi_{\vk_1}+\psi=\measuredangle \vk_2$. Assuming that the rotations $\{\chi_n\}_{n=1}^N$ are drawn uniformly from the rotation group $\mathrm{SO}(3)$, Kam proposed to estimate the autocorrelation function $C_N(k_1,k_2,\psi)$ of the Fourier transform $\widehat{\rho}(\vk)$ in the frequency domain by averaging $\widehat{S}_n(\vk_1)\widehat{S}_n^*(\vk_2)$ over all the projections, or, equivalently, over all the Fourier slices \cite{KAM1977},
\begin{align}
\label{eq:C_k1_k2_psi}
\begin{split}
    C_N\left(k_1,k_2,\psi\right) &= \frac{1}{N}\sum_{n=1}^N\frac{1}{2\pi}\int_0^{2\pi}\widehat{S}_n\left(k_1,\varphi_{\vk_1}\right)\cdot\widehat{S}_n^*\left(k_2,\varphi_{\vk_1}+\psi\right)\ d\varphi_{\vk_1}\\
    &\overset{N\rightarrow\infty}{\longrightarrow}\underbrace{\frac{1}{4\pi}\sum_{l=0}^\infty P_l(\cos\psi)\sum_{m=-l}^lA_{lm}(k_1)\cdot A_{lm}^*(k_2)}_{\textstyle =: C(k_1, k_2, \psi)} + \ \zeta(k_1, k_2, \psi)\,,
\end{split}
\end{align}
where $P_l(\cdot)$ is the Legendre polynomial of degree $l$, and $C(k_1,k_2,\psi)$ is the asymptotic noiseless autocorrelation function of the Fourier transform $\widehat{\rho}(\vk)$ in the frequency domain over the rotation group $\textrm{SO}(3)$.
$C(k_1, k_2, \psi)$ is a nonlinear quadratic function of $\widehat{\rho}(\vk)$.  
When the projections are corrupted with noise, $C_N(k_1,k_2,\psi)$ is a \emph{biased} estimator of $C(k_1,k_2,\psi)$. 
The bias is captured in the additional term $\zeta(k_1,k_2,\psi)$ in \eqref{eq:C_k1_k2_psi}, its asymptotic expression as $N\to\infty$ is derived in Appendix \ref{app:subsec:debias_autocorrelation_features}.

The direct computation of the covariance $C_N(k_1,k_2,\psi)$ via \eqref{eq:C_k1_k2_psi} has high complexity. Alternatively, $C(k_1,k_2,\psi)$ can be computed from the rotation-invariant covariance matrix and mean of the clean projection images. The covariance matrix and mean image can be efficiently estimated from noisy images using the fast steerable PCA \cite{zhao2016fast,zhao2013fourier}. In addition, for images modified by the contrast transfer functions, one can use covariance Wiener filtering (CWF) for the estimation \cite{bhamre2016denoising}. We shall use $\widetilde{C}(k_1,k_2,\psi)$ to denote the estimated $C(k_1,k_2,\psi)$.

Kam further computed the contribution of the subspace of all spherical harmonics with degree $l$, using orthogonality of Legendre polynomials \cite{KAM1977,KAM1980},
\begin{align}
\label{eq:l_subspace_autocorrelation}
\begin{split}
    C_l\left(k_1,k_2\right)&=2\pi(2l+1)\int_0^{\pi}C\left(k_1,k_2,\psi\right)\cdot P_l(\cos\psi)\cdot\sin\psi\ d\psi\\
    &=\sum_{m=-l}^lA_{lm}\left(k_1\right)\cdot A_{lm}^*\left(k_2\right)\,.
\end{split}
\end{align}

Discretizing $k$ into $U$ levels, $k\in\{u_1,u_2,\cdots,u_U\}$, the (discretized) $l$-subspace features $C_l(k_1,k_2)$ can be organized in a $U\times U$ matrix $\mC_l$ of rank $2l+1$ (assuming $U\geq 2l+1$):
\begin{align}
    \mC_l=\left[\begin{array}{cccc}
    C_l(u_1,u_1) &C_l(u_1,u_2) &\cdots &C_l(u_1,u_U)\\
    C_l(u_2,u_1) &C_l(u_2,u_2) &\cdots &C_l(u_2,u_U)\\
    \vdots &\vdots &\ddots &\vdots\\
    C_l(u_U,u_1) &C_l(u_U,u_2) &\cdots &C_l(u_U,u_U)
    \end{array}\right]\,.
\end{align}
Let $\mA_l$ denote the $U\times (2l+1)$ matrix of discretized coefficients $\{A_{lm}(k)\}_{lm}$, 
\begin{align}
\label{eq:A_mat}
    \mA_l= \left[\begin{array}{cccc}
    A_{ll}(u_1) &A_{l(l-1)}(u_1) &\cdots &A_{l(-l)}(u_1)\\
    A_{ll}(u_2) &A_{l{l-1}}(u_2) &\cdots &A_{l(-l)}(u_2)\\
    \vdots &\vdots &\ddots &\vdots \\
    A_{ll}(u_U) &A_{l(l-1)}(u_U) &\cdots &A_{l(-l)}(u_U)
    \end{array}\right]\,.
\end{align}
We then have $\mC_l=\mA_l\mA_l^*$, where $\mA_l$ is real for even $l$ and purely imaginary for odd $l$. OMR aims to recover $\mA_l$ from the Cholesky decomposition of $\mC_l$: $\mC_l=\mF_l\mF_l^*$ \cite{Bhamre:OMR:2015}. The matrix $\mF_l$ returned by the Cholesky decomposition is in general different from $\mA_l$ since for any orthogonal 
matrix $\mO_l$ it holds that $\mC_l = (\mA_l\mO_l^T) (\mA_l\mO_l^T)^*$. OMR then attempts to compute the orthogonal matrices $\{\mO_l\}_{l=1}^L$ that result in the \emph{true} but unknown $\mA_l$, i.e., such that
for all $l$, $\mF_l\mO_l \approx \mA_l$ (with possible discrepancies due to sampling and noise). Since the Cholesky decomposition of $\mC_l$ is independent for different $l$, the challenge faced by OMR is to coordinate the orthogonal matrices to reach a consensus so that $\{\mF_l\mO_l\}_l$ collectively produce the correct Fourier transform $\widehat{\rho}(\vk)$ via \eqref{eq:ft_sh_expansion}.

The earlier OMR methods rely on additional information about the density map that is usually unavailable and thus have limited applicability \cite{Bhamre:OMR:2015}. Levin et al. later introduced an improved OMR by projection matching (OMR-PM) that requires (at least) two denoised projection images to retrieve $\{\mO_l\}_{l=1}^L$. Each projection image can be used to determine every other column
of $\mO_l$, and the relative orientation associated with one of the images is needed to merge the results for a completed $\mO_l$. However, the retrieval of $\{\mO_l\}_{l=1}^L$ from a single image is a nonconvex problem without a known closed-form solution and the estimated relative orientation based on the retrieved $\{\mO_l\}_{l=1}^L$ generally contains error. For these reasons, the OMR method is unstable. Additionally, since the orthogonal matrices $\{\mO_l\}_{l=1}^L$ are coupled in determining the density map $\rho(\vr)$, there is no easy way to impose constraints on $\{\mO_l\}_{l=1}^L$ to ensure the recovered density $\widetilde{\rho}(\vr)$ corresponds to a nonnegative physical density. As a result, the recovered $\widetilde{\rho}(\vr)$ generally has negative values, leading to additional reconstruction errors.

\section{Spatial Features for the Method of Moments}
\label{sec:spatial_mom_features}

We propose to formulate the recovery problem in terms of both the density map $\rho(\vr)$ and the orthogonal matrices $\{\mO_l\}_{l=1}^L$. This allows us to alternate between recovering $\rho(\vr)$ subject to nonnegative summation constraints and updating the orthogonal matrices $\{\mO_l\}_{l=1}^L$ with respect to the recovered density $\widetilde{\rho}(\vr)$ using closed-form solutions. Starting from an initialization density,
we then seek a consensus among the orthogonal matrices on the spatial density map that satisfies the nonnegative summation constraints and matches a reference projection image.
To this end, we propose spatial radial and spatial autocorrelation features that are linear and quadratic functionals of $\rho(\vr)$. Using the derived connection between spatial and Fourier autocorrelations, we can finally link the Fourier autocorrelations to $\rho(\vr)$ as well.

\subsection{Spatial Radial Features}
\label{subsec:radial_integration_feature}
From the Fourier slices $\{\widehat{S}_n\}_n$, we compute the first-order moment of the Fourier transform $\widehat{\rho}(\vk)$ by averaging $\widehat{S}_n(\vk)$ over all directions of $\vk$ with the same norm $k$:
\begin{align}
\label{eq:radial_integration_feature_step_0}
\begin{split}
M_N(k)&=\frac{1}{N}\sum_{n=1}^N\frac{1}{2\pi}\int_0^{2\pi}\widehat{S}_n(\vk)\ d\varphi\\
&\overset{N\rightarrow\infty}{\longrightarrow} \frac{1}{4\pi k^2}\int_0^{2\pi}\int_0^{\pi}\widehat{\rho}(\vk)\cdot k^2\sin\theta_{\vk}d\varphi_{\vk}d\theta_{\vk} =\!  \underbrace{\iiint\frac{\sin(kr)}{kr}\rho(\vr)\ d\vr}_{\textstyle =: M(k)}\,,
\end{split}
\end{align}
where $\varphi$ is the azimuth angle of $\vk$ in the Fourier slice $\widehat{S}_n$. The detailed derivation of \eqref{eq:radial_integration_feature_step_0} is given in Appendix \ref{app:sec:radial_integration_feature}. As derived in Appendix \ref{app:subsec:debias_radial_features}, debiasing is not needed for the (linear) Fourier radial feature $M_N(k)$ in \eqref{eq:radial_integration_feature_step_0}. 
Orthogonality of the sine functions yields the sought spatial radial features.

\paragraph{Spatial radial features} The integration of the 3D density map $\rho(\vr)$ on the sphere with radius $r$ is given by
\begin{align}
\label{eq:radial_integration_linear_1}
\begin{split}
    W(r)&=\frac{2r}{\pi}\int_0^\infty k\cdot M(k)\cdot\sin(kr)\ dk =\iiint\rho(\tilde{\vr})\cdot\delta(r-\tilde{r})\ d\tilde{\vr}\,.
\end{split}
\end{align}

Radial feature extraction is summarized in Algorithm \ref{alg:radial_integration}. The radial feature $W(r)$ is a linear functional of the density map $\rho(\vr)$. We additionally compute the total mass $W_\rho$ of $\rho(\vr)$ in the real space,
\begin{align}
\label{eq:total_mass}
    W_\rho =\iiint\rho(\vr)\ d\vr = \int_0^\infty W(r)\ dr\,,
\end{align}
which is a linear constraint on $\rho(\vr)$, and $W(r)$ can be evaluated from data by plugging $M_N(k)$ defined in~\eqref{eq:radial_integration_feature_step_0} into~\eqref{eq:radial_integration_linear_1}.

\begin{algorithm}[tbp]
\caption{Spatial Radial Feature Extraction}
\label{alg:radial_integration}
\begin{algorithmic}[1]
\REQUIRE The collection of 2D projection images $\{S_n|n=1,\cdots,N\}$.
\STATE Compute $\widehat{S}_n(k, \varphi)$ from $S_n$ using non-uniform FFT.
\STATE Estimate the first order moment $M_N(k)$ in \eqref{eq:radial_integration_feature_step_0}, and use it to approximate $M(k)$.
\STATE Compute the radial feature $W(r)$ in \eqref{eq:radial_integration_linear}.
\STATE {\bfseries Return} $W(r)$.
\end{algorithmic}
\end{algorithm}

\subsection{Spatial Autocorrelation Features}
\label{subsec:autocorrelation_features}
We first expand the density map $\rho(\vr)$ using the real spherical harmonics $Y_{lm}(\cdot)$ in the spatial domain, 
\begin{align}
\label{eq:rho_sh_expansion_complex}
    \rho(\vr)&=\sum_{l=0}^\infty\sum_{m=-l}^lB_{lm}(r)\cdot Y_{lm}(\theta_{\vr},\varphi_{\vr})\,,
\end{align}
where $(\theta_\vr,\varphi_\vr)$ is the angular direction of $\vr$, $B_{lm}(r)$ is the spatial spherical harmonic expansion coefficient given by
\begin{align}
\label{eq:B_lm_r}
    B_{lm}(r) = \iiint\rho(\tilde{\vr})\cdot\delta(r-\tilde{r})\cdot Y_{lm}(\theta_{\tilde{\vr}},\varphi_{\tilde{\vr}})\cdot\sin\theta_{\tilde{\vr}}\ d\tilde{r} d\theta_{\tilde{\vr}} d\varphi_{\tilde{\vr}}\,,
\end{align}
where $\delta(\cdot)$ is the Dirac impulse. We then have the following definition:
\paragraph{Spatial autocorrelation feature} The inner product of the spherical harmonic coefficient vectors $\left\{B_{lm}(r_1)\right\}_m$ and $\left\{B_{lm}(r_2)\right\}_m$ is given by
\begin{align}
\label{eq:autocorrelation_1}
\begin{split}
E_l(r_1,r_2)
&=\sum_{m=-l}^l B_{lm}(r_1)\cdot B_{lm}(r_2)\,.
\end{split}
\end{align}
We call $E_l(r_1,r_2)$ the spatial autocorrelation feature in that we can compute the autocorrelation function $E(r_1,r_2,\psi)$ from it as follows
\begin{align}
\label{eq:spatial_autocorrelation_function}
\begin{split}
    &E(r_1,r_2,\psi)=\sum_{l=0}^\infty E_l(r_1,r_2)\cdot P_l(\cos\psi)\\
    &=\frac{1}{2\pi r_1^2 r_2^2}\iiint\iiint\rho(\tilde{\vr}_1)\cdot\rho(\tilde{\vr}_2)\cdot\delta(r_1-\tilde{r}_1)\cdot\delta(r_2-\tilde{r}_2)\cdot\delta(\psi_{\tilde{\vr}_1,\tilde{\vr}_2}-\psi)\ d\tilde{\vr}_1 d\tilde{\vr}_2\,,
\end{split}
\end{align}
where $\psi_{\tilde{\vr}_1,\tilde{\vr}_2}$ is the angle between $\tilde{\vr}_1$ and $\tilde{\vr}_2$, $\psi\in[0,\pi]$, and $E(r_1,r_2,\psi)$ is the autocorrelation function of the density map $\rho(\vr)$ in the spatial domain over the rotation group $\mathrm{SO}(3)$.
The detailed derivation is given in Appendix \ref{app:sec:autocorrelation_fun_spatial}. We note that Kam proposed a more complicated spatial correlation function calculated from the projection images in \cite{KAM1980} which is different from the autocorrelation function in \eqref{eq:spatial_autocorrelation_function}.

Suppose $r$ is sampled from $\{v_1,v_2,\cdots,v_V\}$. We group the (discretized) spatial autocorrelations into a $V\times V$ matrix $\mE_l$ of rank at most $2l+1$ (assuming $V\geq 2l+1$),
\begin{align}
\label{eq:spatial_autocorrelation_matrix}
    \mE_l=\left[\begin{array}{cccc}
    E_l(v_1,v_1) &E_l(v_1,v_2) &\cdots &E_l(v_1,v_V)\\
    E_l(v_2,v_1) &E_l(v_2,v_2) &\cdots &E_l(v_2,v_V)\\
    \vdots &\vdots &\ddots &\vdots\\
    E_l(v_V,v_1) &E_l(v_V,v_2) &\cdots &E_l(v_V,v_V)
    \end{array}\right]\,.
\end{align}
According to the definition in \eqref{eq:autocorrelation_1}, we can write $\mE_l$ as 
\begin{align}
    \mE_l=\mB_l\mB_l^T\,,
\end{align}
where $\mB_l$ is a $V\times (2l+1)$ matrix
\begin{align}
    \label{eq:B_mat}
    \mB_l= \left[\begin{array}{cccc}
    B_{ll}(v_1) &B_{l(l-1)}(v_1) &\cdots &B_{l(-l)}(v_1)\\
    B_{ll}(v_2) &B_{l(l-1)}(v_2) &\cdots &B_{l(-l)}(v_2)\\
    \vdots &\vdots &\ddots &\vdots \\
    B_{ll}(v_V) &B_{l(l-1)}(v_V) &\cdots &B_{l(-l)}(v_V)
    \end{array}\right]\,.
\end{align}
The matrix $\mB_l$ contains the spherical harmonic coefficients $\{B_{lm}(r)\}_{lm}$ of $\rho(\vr)$ in the spatial domain.

\subsection{Connection between Spatial and Fourier Autocorrelations}
\label{subsec:connection_spatial_fourier}
The proposed spatial autocorrelations are related
to Fourier autocorrelations via a spherical Bessel transform. This transform can be derived from the connection between spherical harmonic expansion coefficients of $\rho(\vr)$ and $\widehat{\rho}(\vk)$. Just like the above spatial expansion coefficients $B_{lm}(r)$ in \eqref{eq:B_lm_r}, the Fourier expansion coefficients $A_{lm}(k)$ in \eqref{eq:ft_sh_expansion} can also be computed using the density map $\rho(\vr)$. We begin by expanding the plane wave $e^{\mathfrak{i}\langle\vk,\vr\rangle}$ in spherical harmonics via the Rayleigh equation \cite{Mehrem2011ThePW},
\begin{align}
\label{eq:pw_expansion}
    e^{\mathfrak{i}\langle\vk,\vr\rangle}=4\pi\sum_{l=0}^\infty\sum_{m=-l}^l\mathfrak{i}^l\cdot j_l(kr)\cdot (Y_l^m)^*(\theta_\vk,\varphi_\vk)\cdot Y_l^m(\theta_\vr,\varphi_\vr),
\end{align}
where $j_l(kr)$ is the spherical Bessel function of order $l$, and $Y_l^m(\cdot)$ is the \emph{complex} spherical harmonic. We then expand $\widehat{\rho}(\vk)$ using the complex $Y_l^m(\cdot)$,
\begin{align}
\label{eq:ft_sh_expansion_complex}
    \widehat{\rho}(\vk)=\sum_{l=0}^\infty\sum_{m=-l}^lA_l^m(k)\cdot Y_l^m(\theta_\vk,\varphi_\vk)\,,
\end{align}
where $A_l^m(k)$ are different from real spherical harmonic coefficients $A_{lm}(k)$ in \eqref{eq:ft_sh_expansion}. Combining \eqref{eq:ft},\eqref{eq:pw_expansion},\eqref{eq:ft_sh_expansion_complex}, we get 
\begin{align}
    A_l^m(k) = 4\pi\left(\mathfrak{i}^l\right)^*\cdot\iiint\rho(\vr)\cdot j_l(kr)\cdot (Y_l^m)^*(\theta_\vr,\varphi_\vr)\ d\vr\,.
\end{align}
Expressing $A_{lm}(k)$  in terms of $A_l^m(k)$, we compute
\begin{align}
\label{eq:A_lm_k}
    A_{lm}(k) = 4\pi\left(\mathfrak{i}^l\right)^*\cdot\iiint\rho(\vr)\cdot j_l(kr)\cdot Y_{lm}(\theta_\vr,\varphi_\vr)\ d\vr\,.
\end{align}
Putting together \eqref{eq:B_lm_r} and \eqref{eq:A_lm_k}, we obtain the following proposition:
\begin{proposition}
\label{prop:fourier_spatial}
The spherical harmonic coefficients $(B_{lm}(r))_{l, m}$ of $\rho(\vr)$ and $(A_{lm}(k))_{l, m}$ of $\widehat{\rho}(\vk)$ are related by the following spherical Bessel transforms
\begin{align}
    \label{eq:spatial_from_fourier}
    B_{lm}(r)&=\frac{1}{2\pi^2(\mathfrak{i}^l)^*}\int_0^\infty A_{lm}(k)\cdot j_l(kr)\cdot k^2\ dk,\\
    \label{eq:fourier_from_spatial}
    A_{lm}(k)&=\int_0^\infty 4\pi(\mathfrak{i}^l)^*\cdot B_{lm}(r)\cdot j_l(kr)\cdot r^2\ dr\,,
\end{align}
where $B_{lm}(r)$ is the spherical harmonic coefficient of the density map $\rho(\vr)$ and $A_{lm}(k)$ the spherical harmonic coefficient of the Fourier transform $\widehat{\rho}(\vk)$.
\end{proposition}

Using the spherical Bessel transforms in Proposition \ref{prop:fourier_spatial}, we can connect the spatial and Fourier autocorrelation features as follows
\begin{align}
    \label{eq:spatial_autocorrelation_from_Fourier_autocorrelation}
    E_l(r_1,r_2) &= \frac{1}{4\pi^4}\int_0^\infty \left( \int_0^\infty C_l(k_1,k_2)\cdot j_l(k_1 r_1)k_1^2\ dk_1 \right) \cdot j_l(k_2 r_2) k_2^2\ dk_2\\
    \label{eq:Fourier_autocorrelation_from_spatial_autocorrelation}
    C_l(k_1,k_2)&=8\pi^2\int_0^\infty\left( \int_0^\infty E_l(r_1,r_2)\cdot j_l(k_1 r_1)r_1^2\ d r_1 \right) \cdot  j_l(k_2 r_2)r_2^2\ d r_2\,.
\end{align}
As summarized in Algorithm \ref{alg:spatial_autocorrelation}, we extract the spatial autocorrelation features from the Fourier autocorrelation features according to \eqref{eq:spatial_autocorrelation_from_Fourier_autocorrelation}. 

\begin{algorithm}[tbp]
\caption{Spatial Autocorrelation Feature Extraction }
\label{alg:spatial_autocorrelation}
\begin{algorithmic}[1]
\REQUIRE The collection of 2D projection images $\{S_n|n=1,\cdots,N\}$.
\STATE Compute $\widehat{S}_n(k, \varphi)$ from $S_n$ using non-uniform FFT.
\STATE Estimate the Fourier  autocorrelation function $C_N(k_1,k_2,\psi)$ in \eqref{eq:C_k1_k2_psi}.
\STATE Calculate the debiased and denoised $\widetilde{C}(k_1,k_2,\psi)$, and use it to approximate $C(k_1,k_2,\psi)$.
\STATE Extract the $l$-subspace features $C_l(k_1,k_2)$ in \eqref{eq:l_subspace_autocorrelation}, and save them in a matrix $\mC_l$.
\STATE Compute the autocorrelation feature matrix $\mE_l$:

    \textbullet$\ $ Transform each column of $\mC_l$ according to \eqref{eq:spatial_from_fourier}, and save the transformed matrix as $\mC_l^\prime$.
    
    \textbullet$\ $ Transform each row of $\mC_l^\prime$ according to \eqref{eq:spatial_from_fourier}, and save the transformed matrix as $\mE_l$. 
\STATE {\bfseries Return} $\mE_l$.
\end{algorithmic}
\end{algorithm}

To formulate the reconstruction problem in terms of the density map $\rho(\vr)$, the autocorrelation features must be expressed as functions of $\rho(\vr)$. The proposed spatial autocorrelation feature $E_l(r_1,r_2)$ in \eqref{eq:autocorrelation_1} is a quadratic functional of $\rho(\vr)$. Using \eqref{eq:Fourier_autocorrelation_from_spatial_autocorrelation}, we can also write the Fourier autocorrelation feature $C_l(k_1,k_2)$ as a quadratic functional of $\rho(\vr)$.

It is convenient to formulate \eqref{eq:Fourier_autocorrelation_from_spatial_autocorrelation} in matrix form. To this end, we approximate the spherical Bessel transform in \eqref{eq:fourier_from_spatial} by Gauss--Legendre quadrature (GLQ). Let $\{q_1,\cdots,q_V\}$ denote the GLQ weights associated with the GLQ sampling locations $r\in\{v_1,\cdots,v_V\}$ in the spatial domain, and the frequency sampling radii $k\in\{u_1,\cdots,u_U\}$. The $U\times (2l+1)$ matrix $\mA_l$ in \eqref{eq:A_mat} that contains the coefficients $\{A_{lm}(k)\}_{lm}$ can be computed as
\begin{align}
    \label{eq:evaluate_A_matrix_form}
    \mA_l=\mQ_l^*\mB_l\,,
\end{align}
where the $(i,j)$-th entry of the $V\times U$ matrix $\mQ_l$ is determined according to \eqref{eq:fourier_from_spatial},
\begin{align}
    Q_l(v_i,u_j) = \left(4\pi\mathfrak{i}^l\cdot j_l(v_iu_j)\cdot v_i^2\right)\cdot q_i\,.
\end{align}
We can then write the Fourier autocorrelation features $\mC_l$ as 
\begin{align}
\label{eq:Fourier_autocorrelation_from_spatial_autocorrelation_matrix_form}
    \mC_l=\mA\mA_l^*=\mQ_l^*\mB_l\mB_l^T\mQ_l=\mQ_l^*\mE_l\mQ_l\,.
\end{align}
We will use this link between Fourier autocorrelations $\mC_l$ and $\rho(\vr)$ through the spatial autocorrelations $\mE_l$ to encourage consistency between the orthogonal matrices and the estimate of $\rho$.

\section{Parametric Density Map}
\label{sec:parametric_density_map}
Since the density map has finite spatial support in practice, the 3D spatial domain is discretized into a $G\times G\times G$ Cartesian grid. For simplicity we assume $G$ is odd and let the discrete coordinates range from $-(G-1)/2$ to $(G-1)/2$. We fix the center of mass at the central cell $(0, 0, 0)$. The 3D density map $\rho(\vr)$ is then sampled at grid points within a radius of $\frac{G-1}{2}$ from the origin, leading to the following discrete representation:
\begin{align}
\label{eq:discrete_representation}
    \rho(\vr)=\sum_{d=1}^D w_d\cdot h(\vr-\boldsymbol\mu_d)\,,
\end{align}
where $D$ is the number of sampling locations, $\boldsymbol\mu_d=\left[\mu_d(x)\ \mu_d(y)\ \mu_d(z)\right]^T\in\mathbb{R}^3$ is the coordinate of the $d$-th sampling location on the Cartesian grid, $h(\cdot)$ is a nonnegative bump function associated with the sampling grid, and $ w_d\geq 0$ is the weight corresponding to the $d$-th sampling location. By abuse of notation, for simplicity we also use $\rho(\vr)$ to denote the parametric density map. 

\paragraph{Reference projection} To reduce the computational complexity, we can choose an arbitrary projection image as the reference projection with an identity rotation matrix $\mR=\mI$, and prune away the grid points that are inconsistent with this reference projection. When the SNR level is relatively high, the images can be denoised via a low-pass filter. Here we shall denoise the reference projection using the multi-frequency vector diffusion maps (MFVDM) \cite{MFVDM:Fan:2021}. However, MFVDM (like many other denoising methods) introduces unknown bias to the denoised image which makes it no longer suitable to estimate features. Notwithstanding, it can still be used for support estimation and pruning. Let $\overline{S}(x,y)$ denote the denoised reference projection, and $\mathcal{M}$ denote the set of grid points. We then have 
\begin{align}
\label{eq:location_pruning}
    \mathcal{M}=\left\{\boldsymbol\mu_d\ \left|\ \overline{S}\big(\mu_d(x),\mu_d(y)\big)>\overline{s}, \textnormal{ and } \|\boldsymbol\mu_d\|_2\leq\frac{G-1}{2}\right.\right\}\,,
\end{align}
where $\overline{s}\geq 0$ is a threshold chosen to filter out those small perturbations in the denoised image $\overline{S}(x,y)$. As discussed later in Section \ref{sec:omr}, we also use the denoised reference projection as additional linear features that play an important role in forming the consensus among the orthogonal matrices on the spatial density map.

The bump function $h(\cdot)$ in \eqref{eq:discrete_representation} should be isotropic and have a controlled effective spherical-harmonic bandwidth. We use the isotropic Gaussian function 
\begin{align}
\label{eq:gauss_basis}
     h(\vr)=\frac{1}{(2\pi)^{\frac{3}{2}}\tau^3}\exp\left(-\frac{1}{2}\frac{\|\vr\|_2^2}{\tau^2}\right)\,,
\end{align}
where $\tau$ is the usual width parameter. Assuming that the grid cell has size $1\times 1\times 1$, setting $\tau=\frac{\sqrt{3}}{2}$ empirically yielded optimal performance. The isotropic Gaussian $h(\cdot)$ is conveniently rotation invariant. Its spherical harmonic coefficients decay exponentially with increasing $l$ and are available in closed-form \cite{Atkinson:SH:2012}, enabling efficient computation of autocorrelation features. 
The problem of reconstructing the 3D density map $\rho(\vr)$ is then equivalent to recovering the weights $\{ w_d\ |\ d=1,\cdots,D\}$ in the sampled discrete representation \eqref{eq:discrete_representation}. 

Using the parametric density map defined by \eqref{eq:discrete_representation} and \eqref{eq:gauss_basis}, we can write $\rho(\vr)$ in terms of \emph{real} spherical harmonics as
\begin{align}
\label{eq:3d_density_sh_expansion}
\begin{split}
    \rho(\vr)&=\sum_{d=1}^Dw_d\cdot\sum_{l=0}^\infty\sum_{m=-l}^lg_{lm}(r,\boldsymbol\mu_d)\cdot Y_{lm}(\theta_\vr,\varphi_\vr)\,,
\end{split}
\end{align}
where $g_{lm}(r,\boldsymbol\mu_d)$ is the spherical harmonic expansion coefficient of $h(\vr-\boldsymbol\mu_d)$ for which a closed-form expression exists, and $Y_{lm}(\cdot)$ is the real spherical harmonic function. The  derivation of \eqref{eq:3d_density_sh_expansion} is detailed in Appendix \ref{app:sec:she_gaussian}. 
Plugging \eqref{eq:3d_density_sh_expansion} into \eqref{eq:radial_integration_linear_1} and \eqref{eq:autocorrelation_1}, we finally express the proposed spatial radial and autocorrelation features in terms of the weight vector $\vw$.

\paragraph{Parametric spatial radial feature} Assuming that the density follows the introduced parametric model, there exist real vectors $\vg(r)$ that let us express $W(r)$ as a linear functional of the weight vector $\vw=[w_1,\cdots,w_D]^T$,
\begin{align}
\label{eq:radial_integration_linear}
\begin{split}
    W(r)=\vg(r)^T\vw\,.
\end{split}
\end{align}
The expression for $\vg(r)$ is derived in Appendix \ref{app:sec:radial_integration_feature}. 

\paragraph{Parametric spatial autocorrelation feature} Similarly, assuming the parametric model holds, there exist real vectors $\vg_{lm}(r)$ such that $E_l$ are quadratics in $\vw$,
\begin{align}
\label{eq:autocorrelation_quad}
    E_l(r_1,r_2)&=\vw^T\cdot\left(\sum_{m=-l}^l\vg_{lm}(r_1)\cdot\vg_{lm}^T(r_2)\right)\cdot\vw\,.
\end{align}
We note from \eqref{eq:autocorrelation_quad} that $E_l(r_1,r_2)=E_l(r_2,r_1)$. Hence we only need to compute spatial autocorrelations for those triplets $\{(l,r_1,r_2)\}$ that satisfy $r_1\leq r_2$.
The derivations for $\vg_{lm}(r)$ and \eqref{eq:autocorrelation_quad} are given in Appendix \ref{app:sec:quad_autocorrelation}. 

As mentioned in Section \ref{subsec:problem_formulation}, we assume that $\rho(\vr)$ has compact support and an effective bandwidth of $\pi$. To reduce the computational complexity, we further set a cutoff threshold $L$ on the spherical harmonic degree $l$ when $\rho(\vr)$ in \eqref{eq:3d_density_sh_expansion} can be well approximated by a function of bandwidth $L$.
Summarizing, the domains of $r,k,l$ used to compute the features are
$0\, \leq r \leq \frac{G-1}{2}$, $0\, \leq k \leq \pi$, and $0\, \leq l \leq L$. 
We approximate the integrals with respect to $\varphi$, $\psi$, $k$, and $r$ (cf. Sections \ref{sec:uvt_mom} and \ref{sec:spatial_mom_features}) by the Gauss--Legendre quadrature (GLQ) \cite{Golub:GLQ:1969,Winckel:LGWT:2004}.

\section{Orthogonal Matrix Retrieval with Spatial Consensus}
\label{sec:omr}

In this section we introduce the proposed orthogonal matrix retrieval with spatial consensus (OMR-SC). Starting from the initial density, we update the orthogonal matrices simultaneously with the estimated density so that they are mutually consistent, while respecting the nonnegativity and total-mass constraints on the density, and agreeing with the denoised reference projection.

Since we model the density $\rho(\vr)$ by a weighted sum  \eqref{eq:discrete_representation} of Gaussians \eqref{eq:gauss_basis} on a grid, the reconstruction of $\rho(\vr)$ is cast as a constrained recovery of the mixture-weight vector $\vw$. We estimate the weight vector $\vw$ and update the orthogonal matrices $\{\mO_l\}_{l=1}^L$ in an alternating fashion. We begin by presenting OMR-SC with Fourier autocorrelations and then show how to simply adapt it to spatial autocorrelations.

\subsection{OMR-SC Using Fourier Autocorrelations}
\label{subsec:OMR-SC-F}
As discussed in Section \ref{subsec:connection_spatial_fourier}, Fourier autocorrelations are efficiently computed from spatial autocorrelations using the spherical Bessel transform (cf. Proposition \ref{prop:fourier_spatial}). Let us express the spatial autocorrelations in terms of the weight vector $\vw$. Under the parametric representation in \eqref{eq:autocorrelation_quad}, the $(i,j)$-th entry of the spatial autocorrelation matrix $\mE_l$ in \eqref{eq:spatial_autocorrelation_matrix} can be written as 
\begin{align}
    E_l(v_i,v_j) = \vb_l^T(v_i,\vw)  \vb_l(v_j,\vw)\,,
\end{align}
where $\vb_l(v,\vw)$ is a vector with $(2l+1)$ elements,
\begin{align*}
    \vb_l(v,\vw):=\left[\begin{array}{c}
        \vg_{ll}^T(v)\vw   \\
        \vg_{l(l-1)}^T(v)\vw  \\
        \vdots\\
        \vg_{l(-l)}^T(v)\vw
    \end{array}\right]=\left[\begin{array}{c}
        \vg_{ll}^T(v)  \\
        \vg_{l(l-1)}^T(v)  \\
        \vdots\\
        \vg_{l(-l)}^T(v)
    \end{array}\right]\vw\,.
\end{align*}
The spherical harmonic coefficient matrix $\mB_l$ of size $V\times(2l+1)$ is then
\begin{align}
    \label{eq:B_mat_parametric}
    \mB_l=\left[\begin{array}{c}
    \vb_l^T(v_1,\vw)\\
    \vb_l^T(v_2,\vw)\\
    \vdots\\
    \vb_l^T(v_V,\vw)\end{array}\right]\,.
\end{align}
Each entry in $\mB_l$ is a linear function of the weight vector $\vw$. 

Using $\mA_l=\mQ_l^*\mB_l$ in \eqref{eq:evaluate_A_matrix_form}, we compute the Fourier autocorrelation matrix $\mC_l$ from $\vw$ as
\begin{align}
    \mC_l=\mA\mA_l^*=\mQ_l^*\mB_l\mB_l^T\mQ_l\,,
    \tag{\ref{eq:Fourier_autocorrelation_from_spatial_autocorrelation_matrix_form} revisited}
\end{align}
where $\mA_l$ is purely real for even $l$ and purely imaginary for odd $l$. Since the recovery of $\vw$ from its quadratic functionals $\mC_l$ is nonconvex and in general challenging, the prospect of recovering $\vw$ from the linear functionals in $\mA_l$ is appealing. However, although $\mA_l$ is unique, the decomposition of $\mC_l$ is not. We have for any orthogonal $\mO_l$
of size $(2l+1)\times(2l+1)$ that
\begin{align}
    \mC_l=\mA_l\mO_l^T(\mA_l\mO_l^T)^*=\mF_l\mF_l^*\,,
\end{align}
where $\mF_l=\mA_l\mO_l^T$. The Cholesky decomposition of $\mC_l$ yields $\mF_l$ for some unknown orthogonal matrix $\mO_l$ that needs to be recovered. The orthogonal matrices $\{\mO_l\}_{l=1}^L$ must be consistent across the different degrees $l\in\{0,\cdots,L\}$ so that for all $l$,  $\mF_l\mO_l=\mA_l$ and they generate the correct density. When the sampling radii $k\in\{u_1,\cdots,u_U\}$ are fixed, the matrix $\mA_l$ only depends on $\vw$. 
We thus jointly estimate the weight vector $\vw$ and update the orthogonal matrices $\{\mO_l\}_{l=1}^L$ by alternating between the two tasks. 

As mentioned in Section \ref{sec:parametric_density_map}, we supplement the radial and autocorrelation features by one denoised reference projection image. Since we cannot recover the absolute orientation from uniform unknown-view projections, we can associate an arbitrary projection with the viewing direction $\mR = \mI$ without loss of generality. This reference projection acts as a regularizer that is empirically crucial for successful reconstruction.

Let $\mB_l(\cdot)$ denote the linear operator on $\vw$ such that $\mB_l(\vw)=\mB_l$ in \eqref{eq:B_mat_parametric}, $\vg(v)$ the measurement vector that produces the $v$-th radial feature $W(v)$ in \eqref{eq:radial_integration_linear}, $\vw$ the weight vector to be recovered, $\overline{S}(x,y)$ the denoised reference projection, and $P_{\vw}(x,y)$ the projection of the density map represented by $\vw$ along the $z$-direction corresponding to $\overline{S}$. We then formulate the following nonconvex OMR-SC problem to recover the density in the spatial domain:
\begin{equation}
\tag{OMR-SC-F}
    \begin{aligned}
    \label{eq:fourier_objective_function}
    \minimize_{\vw,\{\mO_l\}_{l=0}^{L}}&\quad f(\vw,\mO_l):=\sum_{l=0}^L\left\|\mF_l\mO_l-\mQ_l^*\mB_l(\vw)\right\|_2^2\\
    &\quad\quad\quad\quad\quad\quad\quad+\lambda\cdot\sum_{v=1}^V\left(\vg(v)^T\vw-W(v)\right)^2\\
    &\quad\quad\quad\quad\quad\quad\quad+\xi\cdot\sum_{x,y}\big(P_{\vw}(x,y)-\overline{S}(x,y)\big)^2,\\
    \textnormal{subject to}&\quad 0\leq w_d\leq W_\rho,\\    &\quad\sum_{d=1}^Dw_d=W_\rho\,,\\
    &\quad \mO_l^T \mO_l =\mO_l\mO_l^T =\mI, \quad l \in \{0, \ldots, L\}\,.
    \end{aligned}
\end{equation}
where $L$ is the spherical harmonic cutoff degree, $\lambda$ and $\xi$ are the regularization parameters corresponding to the mean-squared-error (MSE) of radial and projection features which can be tuned on some training data acquired under the same setting.
The constraints on $\vw$ in \eqref{eq:fourier_objective_function} come from the requirements that the entries of $\vw$ are nonnegative and the integration of $\rho(\vr)$ in $\mathbb{R}^3$ equals the total mass $W_\rho$ in \eqref{eq:total_mass}, and they can be easily enforced in the spatial domain. Together, they define a convex set $\mathcal{S}$ that is a simplex,
\begin{align}
\label{eq:simplex}
    \mathcal{S}=\left\{\vw\ \left|\ 0\leq w_d\leq W_\rho, \textnormal{ and } \sum_{d=1}^D w_d=W_\rho\right.\right\}\,.
\end{align}

\paragraph{Initialization} Since the problem in \eqref{eq:fourier_objective_function} is nonconvex, the initialization $\vw_0$ directly affects the final solutions $\vw$ and $\{\mO\}_{l=0}^L$. We first reconstruct a low-resolution \emph{ab initio} model from the downsampled projection images, and then use it as the initialization to \eqref{eq:fourier_objective_function}. Noting that both the radial and the projection features are linear measurements, the \emph{ab initio} model is initialized by solving the following convex problem: 
\begin{align}
\label{eq:spatial_initialization}
\begin{split}
    \minimize_{\vw^\prime}&\quad \sum_{v=1}^V\left(\vg(v)^T\vw^\prime-W(v)\right)^2+\sum_{x,y}\big(P_{\vw^\prime}(x,y)-\overline{S}(x,y)\big)^2,\\
    \textnormal{subject to}&\quad 0\leq w^\prime_d\leq W_\rho,\\
    &\quad \sum_{d=1}^D w^\prime_d=W_\rho\,.
\end{split}
\end{align}
The above \eqref{eq:spatial_initialization} is an underdetermined linear problem with multiple globally optimal solutions. Such a problem can be solved by the projected gradient descent initialized at $\boldsymbol 0$ that favors the minimum $l_2$-norm solution \cite{ImplicitReg:Gunasekar:2017,Duchi:2008}.

In practice, a single initialization often fails to produce a good reconstruction due to the nonconvexity of \eqref{eq:fourier_objective_function}. 
To obtain multiple initializations for multiple trials,
we note that we have access to multiple denoised projection images which generally yield different solutions to \eqref{eq:spatial_initialization}. We can thus attempt to solve \eqref{eq:fourier_objective_function} starting from multiple initial points, and choose the reconstructions $\vw,\{\mO\}_{l=0}^L$ that minimize the MSE of autocorrelation features: $\sum_{l=0}^L\left\|\mF_l\mO_l-\mQ_l^*\mB_l(\vw)\right\|_2^2$. For the experiments later in Section \ref{sec:exp}, we find that 10 random choices for a reference projection generally suffice.

Given an initialization $\vw_0$, we propose to minimize \eqref{eq:fourier_objective_function} by alternating between solving for $\mO_l$ and $\vw$:
\begin{itemize}
    \item \textbf{($\mO$-update)}. Fix $\vw$, and update $\{\mO_l\}_{l=1}^L$ with respect to $\vw$:
    \begin{align}
    \label{eq:fourier_objective_function_O}
    \begin{split}
        \minimize_{\{ \mO_l\}_{l=0}^{L}}\quad &f_1(\mO_l):=\sum_{l=0}^L\left\|\mF_l\mO_l-\mQ_l^*\mB_l(\vw)\right\|_2^2\\
        \textnormal{subject to}&\quad \mO_l^T \mO_l =\mO_l\mO_l^T =\mI, \quad l \in \{0, \ldots, L\}\,.
    \end{split}
    \end{align}
    This is an orthogonal Procrustes problem, and the closed-form solution for $\{\mO_l\}_{l=1}^L$ is \cite{Schonemann:1966},
    \begin{align}
    \label{eq:orth_mat_fourier}
        \mO_l = \mV_l\mU_l^T\,,
    \end{align}
    where $\mV_l$ and $\mU_l$ are obtained from the singular value decomposition (SVD) of $\mB_l^T(\vw)\mQ_l\mF_l$
    \begin{align}
        \mB_l^T(\vw)\mQ_l\mF_l=\mU_l\boldsymbol\Sigma_l\mV_l^T\,.
    \end{align}
    \item \textbf{($\vw$-update)}. Fix $\{\mO_l\}_{l=1}^L$, and estimate $\vw$ with respect to $\{\mO_l\}_{l=1}^L$:
    \begin{align}
    \label{eq:fourier_objective_function_z}
    \begin{split}
        \minimize_{\vw}&\quad f_2(\vw):=\sum_{l=0}^L\left\|\mF_l\mO_l-\mQ_l^*\mB_l(\vw)\right\|_2^2\\
        &\quad\quad\quad\quad\quad+\lambda\cdot\sum_{v=1}^V\left(\vg(v)^T\vw-W(v)\right)^2\\
        &\quad\quad\quad\quad\quad+\xi\cdot\sum_{x,y}\big(P_{\vw}(x,y)-\overline{S}(x,y)\big)^2,\\
        \textnormal{subject to}&\quad 0\leq w_d\leq W_\rho,\\
        &\quad\sum_{d=1}^Dw_d=W_\rho\,.
    \end{split}
    \end{align}
    The above \eqref{eq:fourier_objective_function_z} is a convex problem which can be solved using projected gradient descent (PGD). In the $(t+1)$-th iteration, we have
    \begin{align}
    \label{eq:pgd_update}
        \vw_{t+1}=\mathcal{P}_\mathcal{S}\big(\vw_t-\eta\cdot\nabla f_2(\vw_t)\big)\,,
    \end{align}
    where $\eta>0$ is the step size, $\nabla f_2(\vw_t)$ is the gradient at $\vw_t$ from the previous $t$-th iteration, the operator $\mathcal{P}_\mathcal{S}(\cdot)$ projects the gradient descent update onto the convex set $\mathcal{S}$. We compute the projection $\mathcal{P}_\mathcal{S}(\cdot)$ efficiently using the method proposed in \cite{Duchi:2008}.
\end{itemize}
The consensus among the orthogonal matrices on the density map begins with the initialization $\vw_0$. Under the requirement that $\vw$ matches the denoised reference projection subject to nonnegative summation constraints, the consensus in the spatial domain is finally reached through the alternating updates of $\{\mO_l\}_{l=1}^L$ and $\vw$. The proposed OMR-SC with Fourier autocorrelations (OMR-SC-F) is summarized in Algorithm \ref{alg:OMR-SC-F}.

\begin{algorithm}[tbp]
\caption{OMR with Spatial Consensus Using Fourier Autocorrelations}
\label{alg:OMR-SC-F}
\begin{algorithmic}[1]
\REQUIRE Denoised reference projections $\{\overline{S}_i\ |\ i=1,\cdots,I\}$, step size $\eta$, convergence threshold $\varsigma$.
\STATE Extract the spatial radial features $\{W(v)\}_v$ and Fourier autocorrelation features $\{\mC_l\}_l$.
\STATE Perform Cholesky decompositions of Fourier autocorrelation matrices $\{\mC_l\}_l$.
\FOR{$i=\{1,\cdots,I\}$}
\STATE Compute the initialization $\vw_0(i)$ from the spatial radial features $\{W(v)\}_v$ and the $i$-th reference projection $\overline{S}_i$.
\FOR{$t=\{0,1,\cdots,T\}$}
    \STATE Fix $\vw_t(i)$, and update $\{\mO_l(i)\}_l$ with respect to $\vw_t(i)$ via singular value decomposition.
    \STATE Fix $\{\mO_l(i)\}_l$, and estimate $\vw_{t+1}(i)$ with respect to $\{\mO_l(i)\}_l$ via projected gradient descent.
    \IF {$\frac{\|\vw_{t+1}(i)-\vw_t(i)\|_2}{\|\vw_t(i)\|_2}<\varsigma$}
        \STATE Convergence is reached, set $\vw(i)=\vw_{t+1}(i)$ and \textbf{break}.
    \ENDIF
\ENDFOR
\STATE Save the $i$-th set of solutions $\{\vw(i), \{\mO_l(i)\}_l\}$.
\ENDFOR
\STATE Find the set of solutions that minimizes the MSE of autocorrelation features:
\begin{align}
    \tilde{i} = \arg\min_i\quad\sum_{l=0}^L\left\|\mF_l\mO_l(i)-\mQ_l^*\mB_l(\vw(i))\right\|_2^2
\end{align}
\STATE {\bfseries Return} $\tilde{\vw}=\vw(\tilde{i})$.
\end{algorithmic}
\end{algorithm}

\subsubsection{Computational Complexity}
Let $N$ denote the number of projection images and $G\times G$ the size of the image. The complexities of the proposed OMR-SC steps are then as follows:
\begin{enumerate}[label={\arabic*)}]
\item \emph{Feature extraction.} 
The complexity of feature extraction is $\mathcal{O}\left(NG^3+LG^3\right)$.

\item \emph{Optimization with respect to $\{\mO_l\}_{l = 1}^L$.} The overall complexity of computing $\{\mO_l\}_{l = 1}^L$ per iteration is $\mathcal{O}(L^3(LG+G^2))$.

\item \emph{Optimization with respect to $\vw$.} The overall complexity of computing $\vw$ per iteration is $\mathcal{O}(L^2G^3+G^3\log G)$.
\end{enumerate}
Detailed derivations are given in Appendix \ref{app:sec:complexity}. We only need to go through the projection images once to extract the features, which is more efficient compared to approaches that estimate the particle orientations of projection images at every iteration. Additionally, the non-uniform FFT of projection images can be parallelized and computed on the fly during data acquisition. For the experiments performed later in Section \ref{sec:exp}, we observed that convergence was reached after $\sim500$ iterations of optimizing $\{\mO_l\}_{l=1}^L$ and $\vw$ in an alternating fashion.

\subsection{OMR-SC Using Spatial Autocorrelations}
The proposed OMR-SC using spatial autocorrelations (OMR-SC-S) is derived in Appendix \ref{app:sec:omr_sc_s}. We empirically observe that the Fourier autocorrelations perform better than spatial autocorrelations (cf. Section \ref{sec:exp}). We conjecture that this is due to the implicit spherical-harmonic frequency marching effect that is significantly stronger when using Fourier autocorrelations. Indeed, Fig. \ref{fig:spatial_fourier_gradient_norm} in Appendix \ref{app:sec:compare_spatial_fourier_autocorrelations} shows that the gradient norms of Fourier features for lower spherical-harmonic frequencies (degrees) $\ell$ are much larger than those for higher frequencies at initialization and remain so through the iterations. This suggests that the low-frequency Fourier features are matched first and given priority during optimization, which gives rise to the aforementioned implicit frequency marching effect. By contrast, the differences among the gradient norms of spatial features are much smaller, and the corresponding frequency marching effect is thus weaker.

\subsection{Comparison with Earlier Orthogonal Matrix Retrieval Methods}
\label{subsec:compare_with_previous_omr}
As we briefly discussed in Section \ref{subsec:kam_autocorrelaton}, the earlier orthogonal matrix retrieval (OMR) approach is based on Kam's (autocorrelation) method, where the $l$-subspace features $C_l(k_1,k_2)$ in \eqref{eq:l_subspace_autocorrelation} are used to recover the spherical harmonic coefficients $\{A_{lm}(k)\}_{lm}$ in the frequency domain. The feature matrix $\mC_l$ and the coefficient matrix $\mA_l$ are related as $\mC_l=\mA_l\mA_l^*$. OMR performs Cholesky decomposition on $\mC_l$ to produce $\mC_l=\mF_l\mF_l^*$, and attempts to recover a set of orthogonal matrices $\{\mO_l\}_{l=1}^L$ such that $\mF_l\mO_l=\mA_l$. The main difficulty OMR faces is how to relate the orthogonal matrices $\{\mO_l\}_{l=1}^L$ of the different degrees $l$. A na\"ive attempt leads to the following ill-posed optimization problem
\begin{align}
\label{eq:omr_naive}
    \minimize_{\mO_l,\mA_l}\quad\|\mF_l\mO_l-\mA_l\|_2^2,\quad \text{for all}~l\in\{0,\cdots,L\}\,.
\end{align}
Since the above \eqref{eq:omr_naive} is an independent problem for different degrees $l$, there is no reason for the recovered $(\mO_l,\mA_l)$ to be consistent with the true underlying density. 
To resolve this, different authors proposed to use different kinds of side information  \cite{Bhamre:OMR:2015,Eithan2017}. The orthogonal extension (OE) method \cite{Bhamre:OMR:2015}, for example, assumes a similar 3D structure is known so that we can compute its spherical harmonic coefficients $\mA_l^\prime$, and recover $\mO_l$ by replacing $\mA_l$ with $\mA_l^\prime$,
\begin{align}
    \label{eq:omr_oe}
    \min_{\mO_l}\quad\|\mF_l\mO_l-\mA_l^\prime\|_2^2,\quad \text{for all}~ l\in\{0,\cdots,L\}\,.
\end{align}
Such additional information is usually unavailable, which limits the method's practical value.

A different approach known as OMR by projection matching (OMR-PM) was introduced in \cite{Eithan2017} which instead relies on (at least) two denoised projections with estimated relative rotations to compute the orthogonal matrices. Let $\{\mO_{l;1}\}_l$ denote the orthogonal matrices recovered from the first projection, $\{\mO_{l;2}\}_l$ the orthogonal matrices recovered from the second projection, $\chi$ the unknown relative rotation associated with the second projection, and $\mD^{(\chi)}_l$ the corresponding Wigner D-matrix at degree $l$. It was shown in \cite{Eithan2017} that every other column of $\mO_{l;1}\mD^{(\chi)}_l$ should equal the corresponding column of $\mO_{l;2}$. The relative rotation $\chi$ is estimated by a dense grid search over the rotation group $\text{SO}(3)$ which tries to match every other column of $\mO_{l;1}\mD^{(\chi)}_l$ and $\mO_{l;2}$. The estimated rotation $\chi$ is then used to \emph{merge} $\{\mO_{l;1}\}_l$ and $\{\mO_{l;2}\}_l$ to produce the final solution $\{\mO_l\}_{l=1}^L$. Since the retrieval of $\{\mO_{l;1}\}_l$ (or $\{\mO_{l;2}\}_l$) from a single projection is a nonconvex problem that does not have a known closed-form solution,
this introduces errors in the estimated rotation $\chi$ which makes the final merging step unstable. Additionally, the orthogonal matrices $\{\mO_l\}_{l=1}^L$ are coupled to one another in determining the density map $\rho(\vr)$ since, for example, they should result in a nonnegative density, but in this formulation there is no easy way to implement the nonnegativity constraint.

Our proposed OMR-SC approach alternates between recovering the spatial density map $\rho(\vr)$ and updating the orthogonal matrices $\{\mO_l\}_{l=1}^L$ with respect to the recovered density $\widetilde{\rho}(\vr)$. 
It is thus straightforward to enforce nonnegativity via a projection onto the simplex in \eqref{eq:simplex}. Compared to a single projection used by OMR-PM, the 3D density map in the OMR-SC formulation, together with the appropriate constraints, contains all the information needed to determine the orthogonal matrices, and the updates of $\{\mO_l\}_{l=1}^L$ have closed-form solutions in \eqref{eq:orth_mat_fourier}. With a suitable initialization as computed in \eqref{eq:spatial_initialization}, OMR-SC is more robust in recovering the density map and the orthogonal matrices. This is enabled by the new spatial autocorrelation features and their relation to Fourier autocorrelations in Proposition \ref{prop:fourier_spatial}.

\begin{figure*}[tb]
\centering
\subfloat[]{
\label{fig:proj_clean}
\includegraphics[height=1.2in]{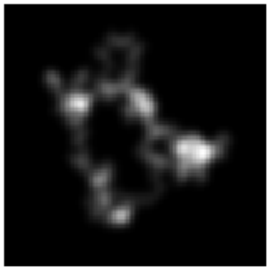}}
\subfloat[]{
\label{fig:proj_noisy}
\includegraphics[height=1.2in]{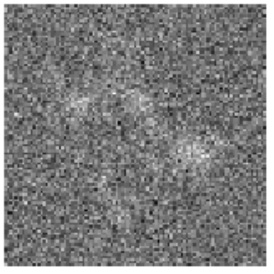}}
\caption{(a) A noiseless 2D projection image; (b) A noisy 2D projection image with SNR=$0.1$. 
}
\label{fig:proj_clean_noisy}
\end{figure*}

\section{Experimental Results}
\label{sec:exp}
In this section we compare the proposed OMR-SC approach with the OMR-PM approach on the recovery of 3D density maps.  
\begin{itemize}
\item We first reconstruct 10 random density maps\footnote{Due to the high computational complexity of UVT applications such as cryo-EM, the number of reconstructions is generally limited to several density maps.}: each groundtruth density map is a mixture of Gaussian components whose means are generated using a 3D random walk with 500 steps and variances are set to 1. The model is further scaled to fit within the ball embedded in a $101\times 101 \times 101$ Cartesian grid.
\item We then reconstruct 3 protein density maps: 1) the human calcium-sensing receptor (CaS) density model from Electron Microscopy Data Bank (EMDB) \cite{HCSR:2021}; 2) the Holliday junction complex (HJC) density model;  
3) the human patched 1 protein (PTCH1) density model.
HJC and PTCH1 maps are synthesized in Chimera~\cite{chimera} from their atomic models in the Protein Data Bank (PDB). The 3 models are downsampled and scaled to fit within the ball embedded in a $101\times 101\times 101$ Cartesian grid, with voxels corresponding to cubes with physical side lengths 2.1672 \r{A}, 2.5 \r{A} and 1.5 \r{A}, respectively. 
\end{itemize}
Without loss of generality, the density maps are normalized so that the total mass of each density $W_\rho=50$. The 2D projection images are generated from $N$ uniformly-distributed 3D rotations. As shown in Fig. \ref{fig:proj_clean_noisy}, we generate $N=10,000$ noiseless projection images of size $101\times 101$, and corrupt them with additive white Gaussian noise (AWGN) so that the resulting $\mathrm{SNR}=0.1$ as in the typical low-SNR scenario of UVT, where
\[\mathrm{SNR}=\frac{\textnormal{Power of Signal}}{\textnormal{Power of Noise}}=\frac{\mathbb{E}\left[\sum_{x,y}S(x,y)^2\right]}{\mathbb{E}\left[\sum_{x,y}\epsilon(x,y)^2\right]}\,.\]

\subsection{Feature Extraction}
\label{subsec:feature_extraction} 
The integrations involved in feature extraction are computed according to the Gauss--Legendre quadrature (GLQ) rule. Using a holdout random density as the training data, we empirically (through trial and error) determined that the following choices provide a good balance between complexity and accuracy:
$k_{\max}=\frac{\pi}{2}$ as the cutoff frequency (approximate bandlimit) of the 3D density and its projection images, $\Phi=401$ GLQ points for $\varphi\in[0,2\pi]$, $V=101$ GLQ points for $r\in[0,50]$, $U=51$ uniformly sampled points for the frequency bandwidth $k\in[0,\frac{\pi}{2}]$, and the spherical harmonic bandwidth $L=10$. 
The non-uniform FFT of the projection image is calculated using the FINUFFT package \cite{Barnett:NUFFT_1:2019,Barnett:NUFFT_2:2021}.

\begin{figure*}[tb]
\centering
\subfloat[Noiseless $\mC_{1}$]{
\label{fig:autocorrelation_l1_noiseless}
\includegraphics[width=0.31\textwidth]{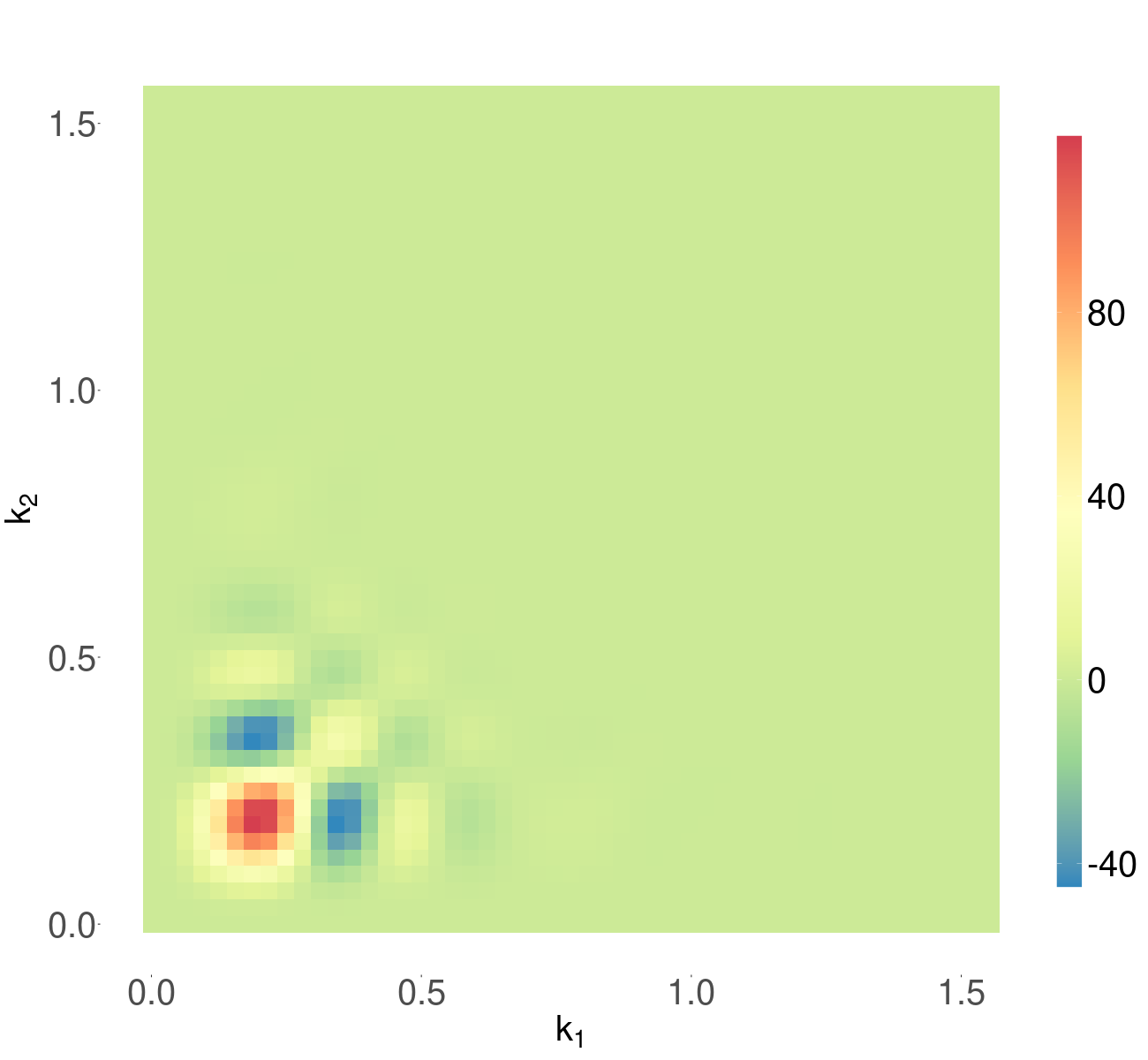}}
\subfloat[Noiseless $\mC_{5}$]{
\label{fig:autocorrelation_l5_noiseless}
\includegraphics[width=0.31\textwidth]{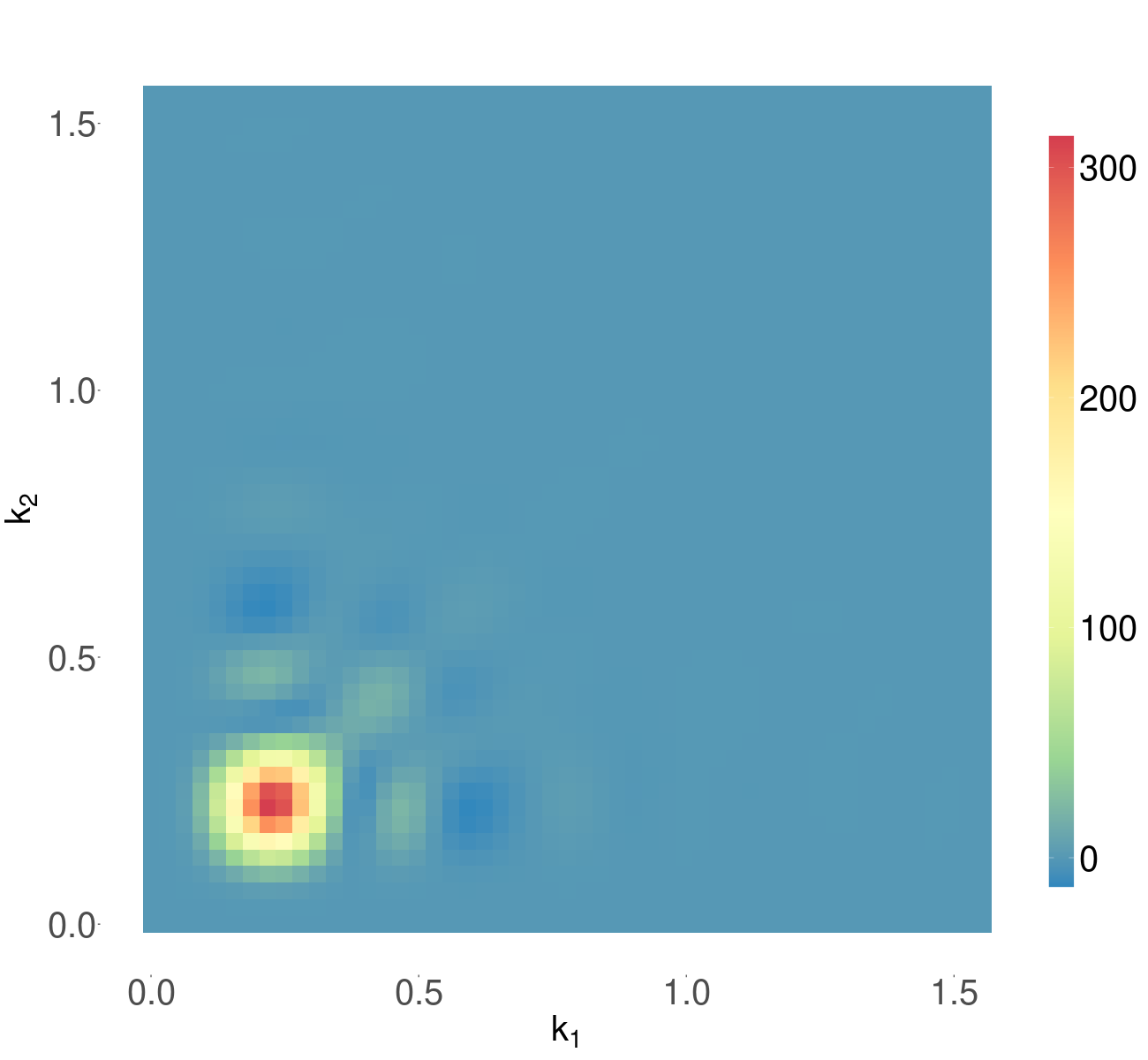}}
\subfloat[Noiseless $\mC_{10}$]{
\label{fig:autocorrelation_l10_noiseless}
\includegraphics[width=0.31\textwidth]{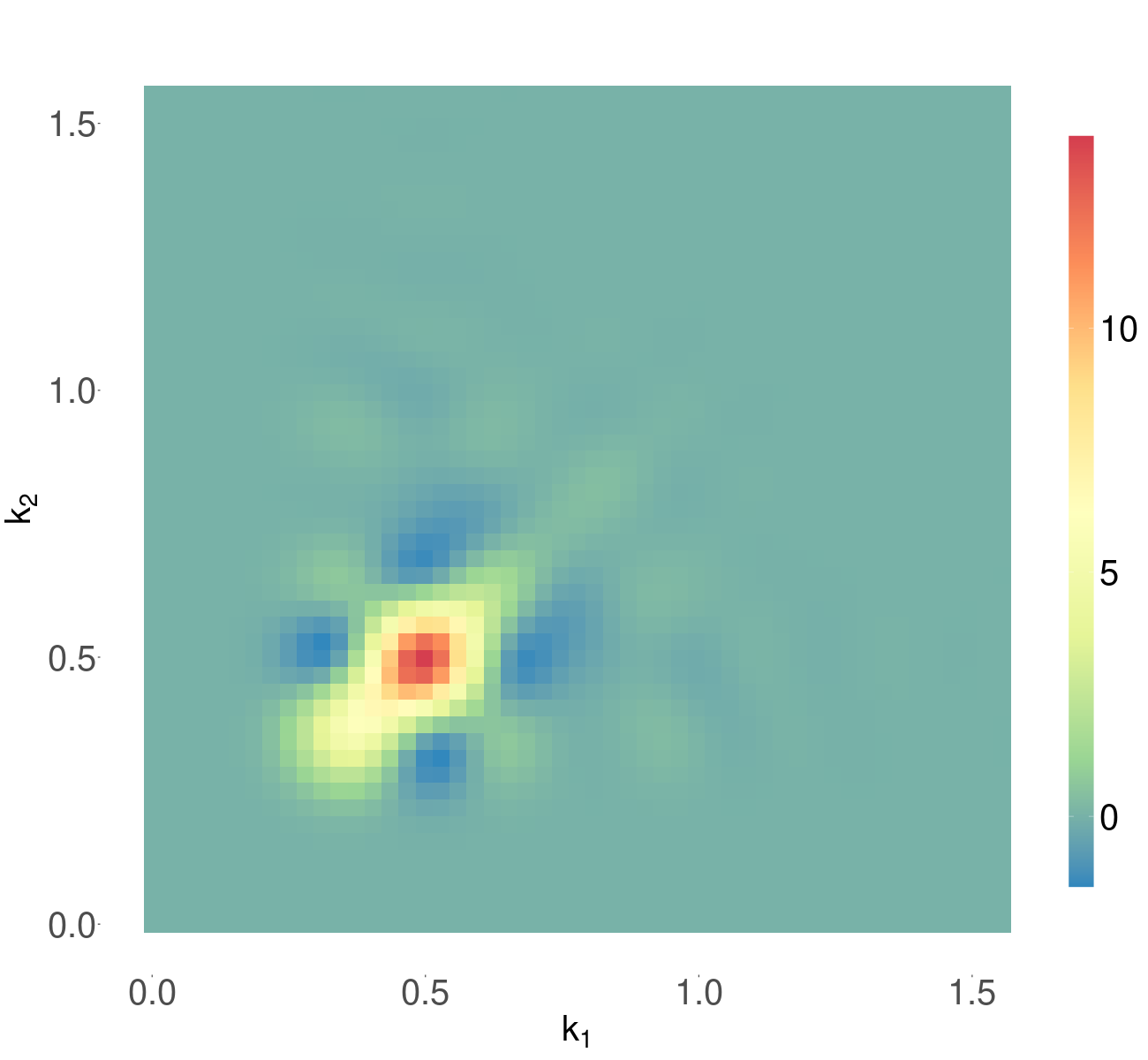}}\\
\subfloat[NRMSE of $\widetilde{\mC}_1$: 0.0821]{
\label{fig:autocorrelation_l1_noisy_denoise}
\includegraphics[width=0.31\textwidth]{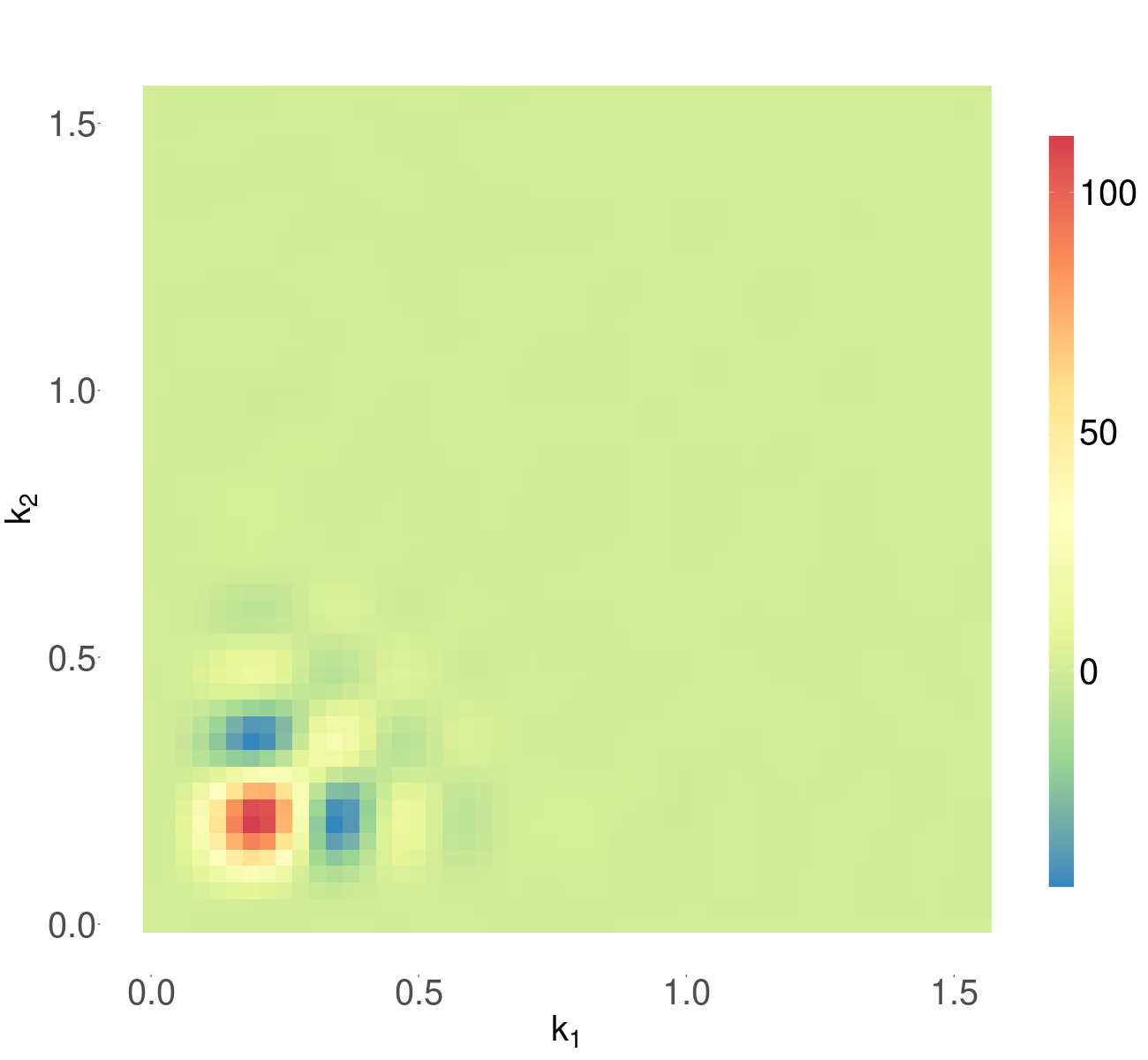}}
\subfloat[NRMSE of $\widetilde{\mC}_5$: 0.0197]{
\label{fig:autocorrelation_l5_noisy_denoise}
\includegraphics[width=0.31\textwidth]{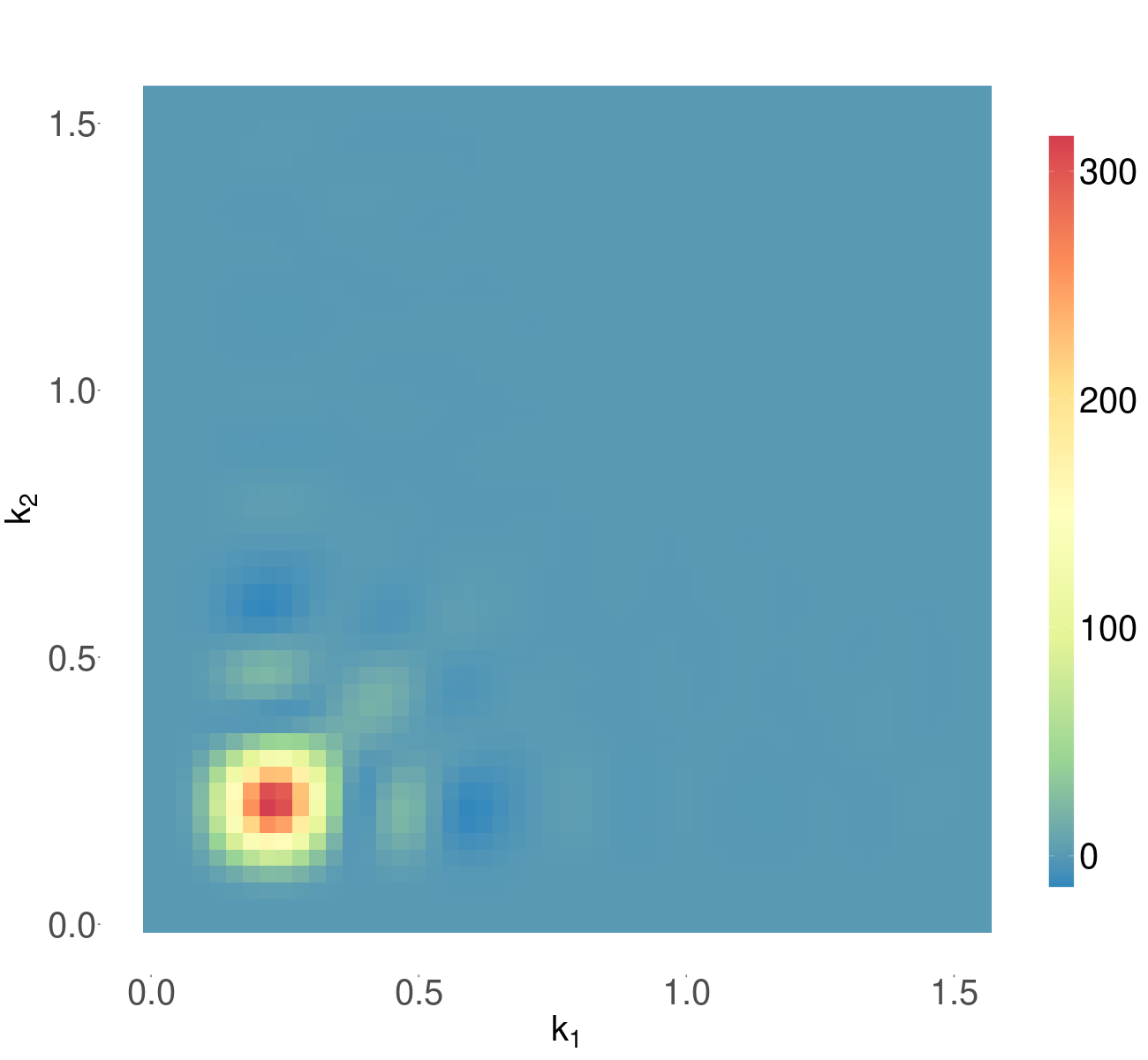}}
\subfloat[NRMSE of $\widetilde{\mC}_{10}$: 0.1034]{
\label{fig:autocorrelation_l10_noisy_denoise}
\includegraphics[width=0.31\textwidth]{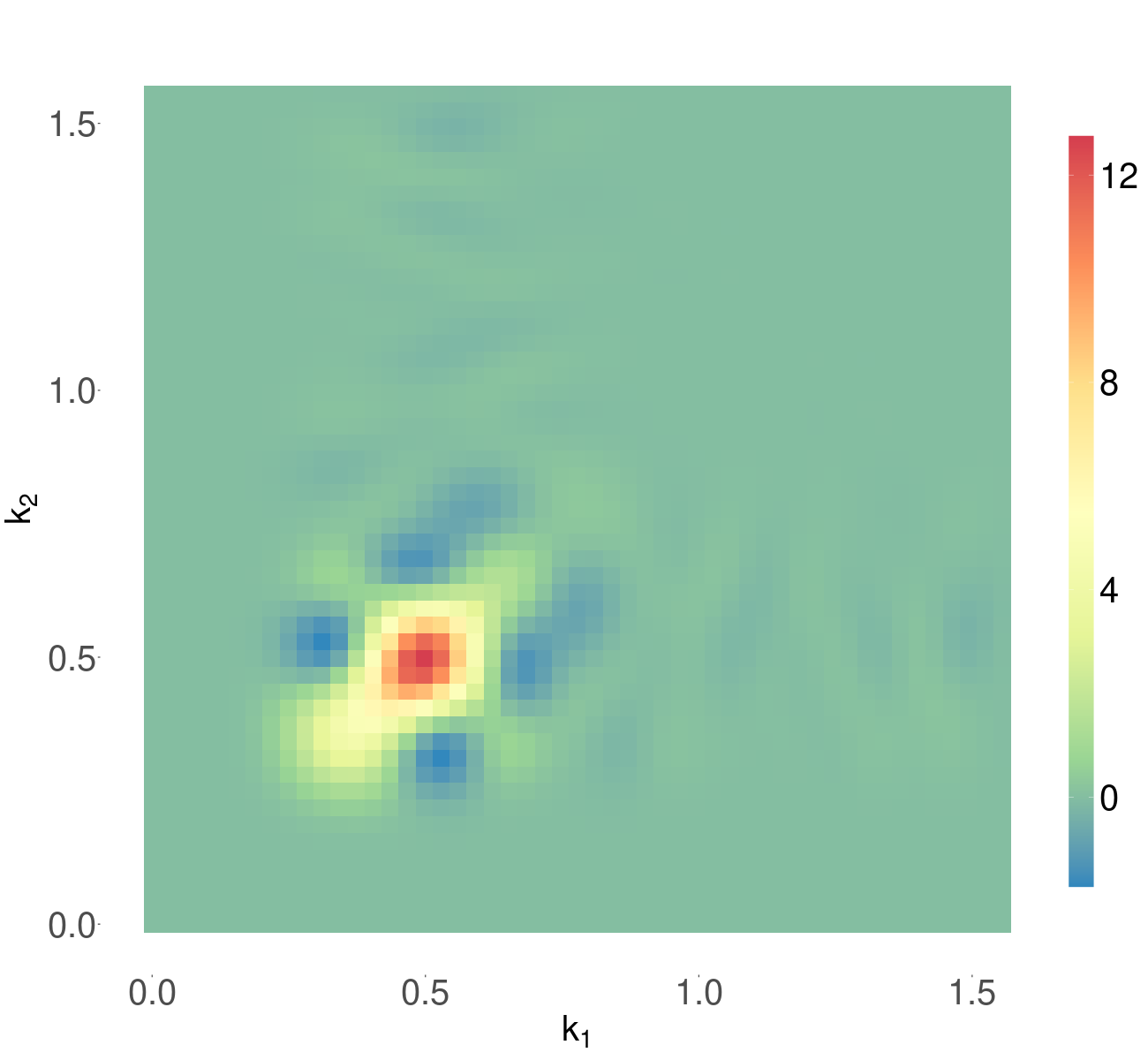}}
\caption{Fourier autocorrelation feature $\widetilde{\mC}_l$ is extracted from $N=10,000$ projection images with SNR=0.1. }
\label{fig:correlation_feature_noisy}
\end{figure*}

The linear radial features do not need to be debiased nor denoised. Using the ASPIRE package \cite{ASPIRE:2018}, we calculate the debiased and denoised autocorrelation features $\widetilde{\mC}_l$ via the fast steerable PCA \cite{zhao2016fast}. As shown in Fig. \ref{fig:correlation_feature_noisy}, the normalized-root-mean-squared-error (NRMSE) of the features is usually larger for higher spherical harmonic degrees.

\subsection{Reconstruction of 3D Density Maps}
\label{subsec:reconstruction_3d_density_maps}
The proposed OMR-SC approach leverages one denoised reference projection image as additional linear features. We denoise this one projection by multi-frequency vector diffusion maps \cite{MFVDM:Fan:2021}. Although the denoised images contain unknown bias and cannot be used
for computing MoM features, they can be used as linear ``measurements'' in our recovery formulation. As discussed in Section \ref{sec:parametric_density_map}, we can prune away the sampling locations $\{\boldsymbol\mu_d\}_d$ on the Cartesian grid that are not consistent with the reference projection to reduce the computational complexity. The regularization parameters $\lambda$ and $\xi$ in \ref{eq:fourier_objective_function} are both set to $100$ (by monitoring the loss on a training random density map generated under the same setting), and they are both set to $1$ in \ref{eq:spatial_objective_function}.

As discussed in Section \ref{sec:omr}, we can combine spatial radial features with different denoised images to compute different initializations for OMR-SC. Here we use 10 randomly selected projections to perform the reconstructions in parallel. In practice, we first perform \emph{ab initio} modeling via OMR-SC where a low-resolution density map is reconstructed by downsampling the projection images. Specifically, a $33\times 33 \times 33$ \emph{ab initio} model is computed for every density map.
We then refine the \emph{ab initio} models via OMR-SC again to get the high-resolution reconstruction. Among all the obtained solutions, we choose the one that minimizes the ``rotation-invariant'' MSE of autocorrelation features. On the other hand, the objective function of OMR-PM is not rotation-invariant and depends on the chosen set of projections \cite{Eithan2017}. Hence it could not be compared across different sets of projections to select the best solution.
As a result, we use all of the same 10 projections in OMR-PM to achieve the best performance. 

We also note that OMR-PM uses a different set of autocorrelation features constructed from the spherical Bessel expansion coefficients $a_{lms}$ of $A_{lm}(k)$ \cite{BhamreZS17}:
\begin{align}
    A_{lm}(k)=\sum_{s=1}^{S_l} a_{lms}\cdot j_{ls}(k)\,,
\end{align}
where $j_{ls}(k)$ is the normalized spherical Bessel function. We can see that $a_{lms}$ is connected to $A_{lm}(k)$ through a linear spherical Bessel transform, and our previous discussion of OMR-PM based on $A_{lm}(k)$ would still hold in this case. The proposed OMR-SC could also be adapted straightforwardly to use $a_{lms}$ as features, which produces similar performance to $A_{lm}(k)$.

A standard reconstruction quality metric in UVT is the Fourier shell correlation (FSC). It measures the similarity between two volumes ($\rho_1$, $\rho_2$) with respect to the spatial frequency, and is given by the normalized cross-correlation coefficient between two aligned volumes over corresponding spherical shells with radius $k$ in the frequency domain \cite{Harauz1986ExactFF}: 
\begin{align}
    \textnormal{FSC}(k)=\frac{\sum_{k_i=k}\widehat{\rho}_1(\vk_i)\widehat{\rho}_2(\vk_i)^*}{\sqrt{\sum_{k_i=k}|\widehat{\rho}_1(\vk_i)|^2}\cdot\sqrt{\sum_{k_i=k}|\widehat{\rho}_2(\vk_i)|^2}}\,,
\end{align}
where $\widehat{\rho}_1$ and $\widehat{\rho}_2$ are the Fourier transforms of the two aligned volumes, and $\vk_i$ corresponds to the $i$-th voxel in the frequency domain.
We apply a cutoff threshold of $0.5$ on the FSC curve to determine the volume resolution \cite{SPAHN:2004}. To evaluate how well the reconstruction matches the groundtruth globally, we calculate another standard metric, the correlation coefficient \cite{Afonine:2018}.

\subsubsection{Random Density Maps}
Table \ref{tab:compare_resolution} shows the resolutions (in voxel) of recovered random density maps using OMR-PM and OMR-SC.
Table \ref{tab:compare_cc} shows the corresponding correlation coefficients. We see that OMR-SC-F with Fourier autocorrelations generally performs better than OMR-SC-S with spatial autocorrelations. Without noise, OMR-SC-F performs better than OMR-PM on all the densities. With noise, OMR-SC-F performs significantly better than OMR-PM on all the densities. In both cases, OMR-SC-F is more robust across the different densities than OMR-PM. The information used by OMR-PM and OMR-SC is essentially the same. As discussed in Section \ref{subsec:compare_with_previous_omr}, the performance differences are due to the different problem formulations and their corresponding optimization procedures.

\begin{table}[tb]
\caption{Resolutions (in voxel) of recovered random density maps using the OMR-PM and OMR-SC approaches (FSC cutoff threshold = 0.5).}
\label{tab:compare_resolution}
\centering
\resizebox{\textwidth}{!}{
\begin{tabular}{llcccccccccc}
\toprule
& &D1 &D2 &D3 &D4 &D5 &D6 &D7 &D8 &D9 &D10\\ \midrule

&OMR-PM & 11.64 & 11.90 & 7.00 & 14.14 & 11.15 & 10.95 & 13.68 & 9.57 & 11.93 & 21.93 \\ 
&OMR-SC-S & 8.83 & 9.32 & 8.98 & 12.06 & 11.33 & 8.33 & 13.39 & 12.67 & 12.24 & 10.19 \\
\multirow{-3}{*}{Noiseless} &OMR-SC-F & \bf{6.70} & \bf{9.62} & \bf{5.90} & \bf{7.28} & \bf{5.94} & \bf{4.10} & \bf{10.63} & \bf{7.48} & \bf{8.46} & \bf{5.54}   \\ \midrule

&OMR-PM & 34.13 & 31.55 & 11.71 & 21.98 & 27.32 & 17.42 & 28.65 & 14.03 & 17.64 & 26.11  \\ 
&OMR-SC-S & 13.09 & 15.20 & 7.98 & 14.73 & 10.48 & 6.45 & 21.37 & 16.13 & 14.14 & 12.15 \\
\multirow{-3}{*}{Noisy} &OMR-SC-F & \bf{8.47} & \bf{10.80} & \bf{7.42} & \bf{10.98} & \bf{9.18} & \bf{6.01} & \bf{17.51} & \bf{8.94} & \bf{12.14} & \bf{8.45} \\ \bottomrule
\end{tabular}
}
\end{table}

\begin{table}[tb]
\caption{Correlation coefficients of recovered random density maps using the OMR-PM and OMR-SC approaches.}
\label{tab:compare_cc}
\centering
\resizebox{\textwidth}{!}{
\begin{tabular}{llcccccccccc}
\toprule
& &D1 &D2 &D3 &D4 &D5 &D6 &D7 &D8 &D9 &D10\\ \midrule

&OMR-PM & 0.81 & 0.69 & 0.90 & 0.76 & 0.75 & 0.61 & 0.79 & 0.86 & 0.84 & 0.47 \\ 
&OMR-SC-S & 0.78 & 0.75 & 0.83 & 0.75 & 0.76 & 0.73 & 0.80 & 0.75 & 0.79 & 0.80  \\
\multirow{-3}{*}{Noiseless} &OMR-SC-F & \bf{0.92} & \bf{0.84} & \bf{0.96} & \bf{0.92} & \bf{0.96} & \bf{0.96} & \bf{0.83} & \bf{0.91} & \bf{0.91} & \bf{0.91}  \\ \midrule

&OMR-PM & 0.39 & 0.45 & 0.78 & 0.56 & 0.43 & 0.50 & 0.58 & 0.77 & 0.67 & 0.38  \\ 
&OMR-SC-S & 0.68 & 0.62 & 0.85 & 0.74 & 0.75 & 0.80 & 0.62 & 0.68 & \bf{0.77} & 0.67    \\
\multirow{-3}{*}{Noisy} &OMR-SC-F & \bf{0.85} & \bf{0.71} & \bf{0.92} & \bf{0.82} & \bf{0.88} & \bf{0.89} & \bf{0.71} & \bf{0.85} & \bf{0.77} & \bf{0.86} \\ \bottomrule
\end{tabular}
}
\end{table}

\begin{figure*}[tb]
\centering
\subfloat[Density 1 (D1)]{
\label{fig:fsc_d1}
\includegraphics[width=0.45\textwidth]{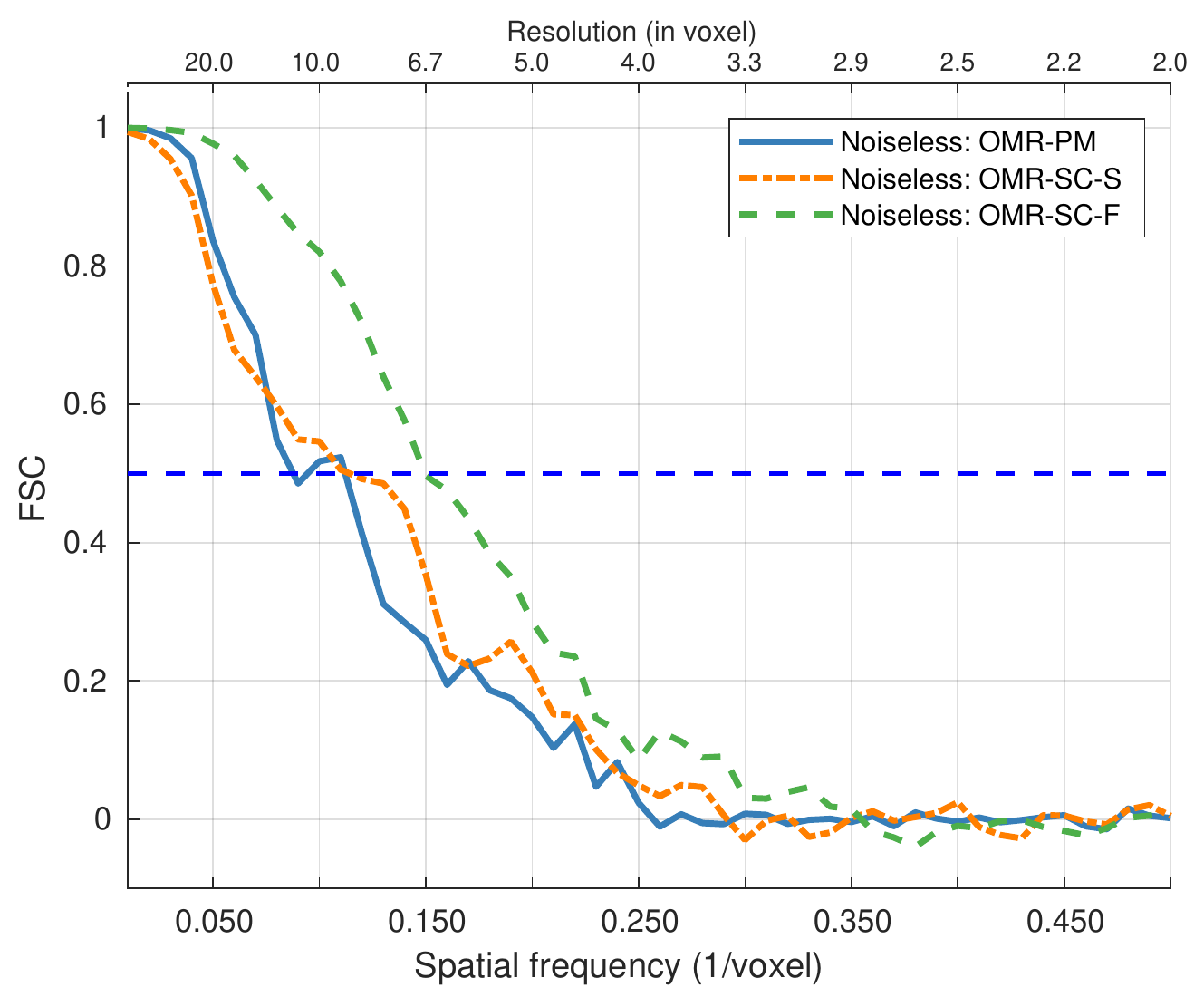}
\hspace{0.03\textwidth}
\includegraphics[width=0.45\textwidth]{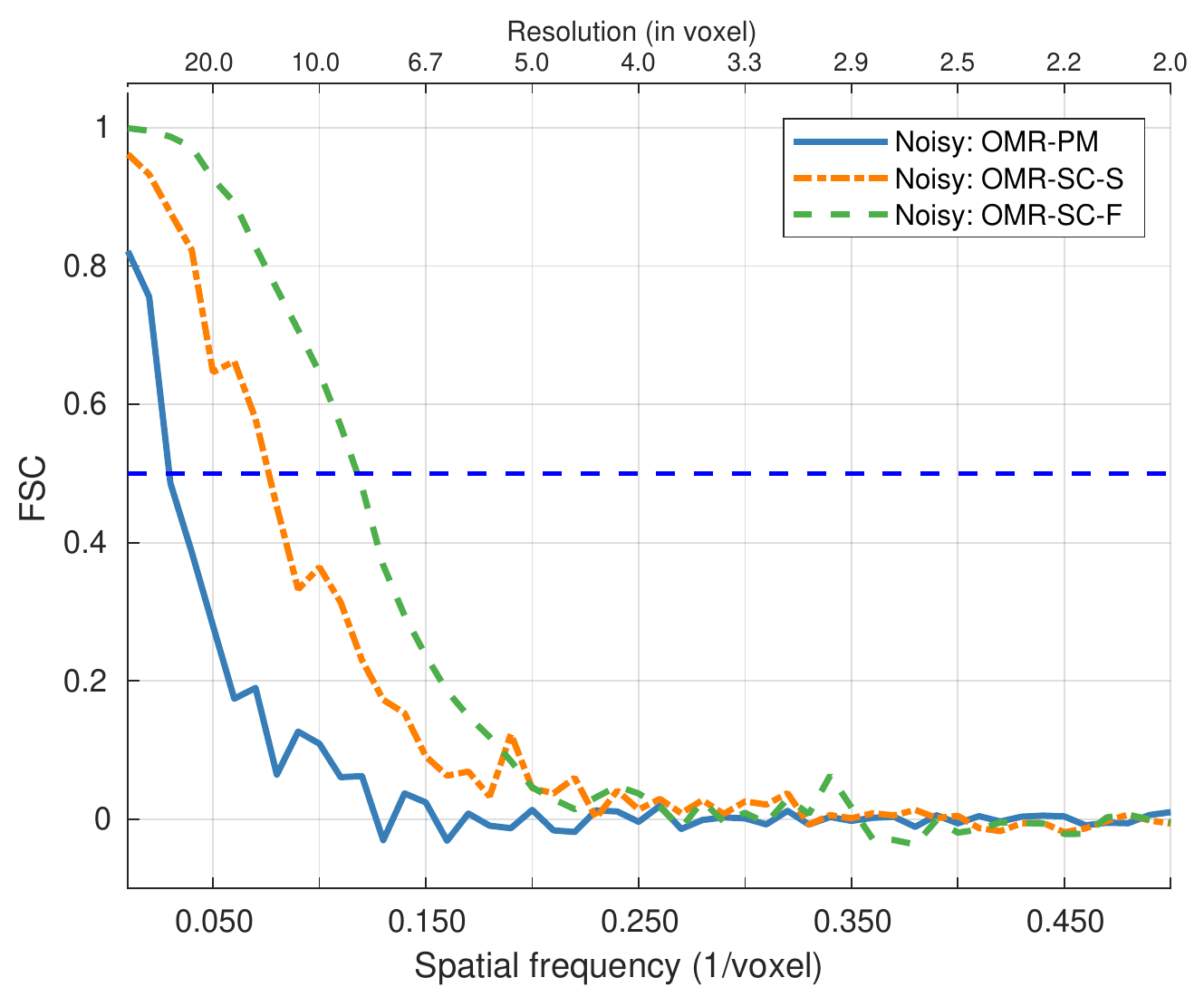}}\\
\subfloat[Density 10 (D10)]{
\label{fig:fsc_d10}
\includegraphics[width=0.45\textwidth]{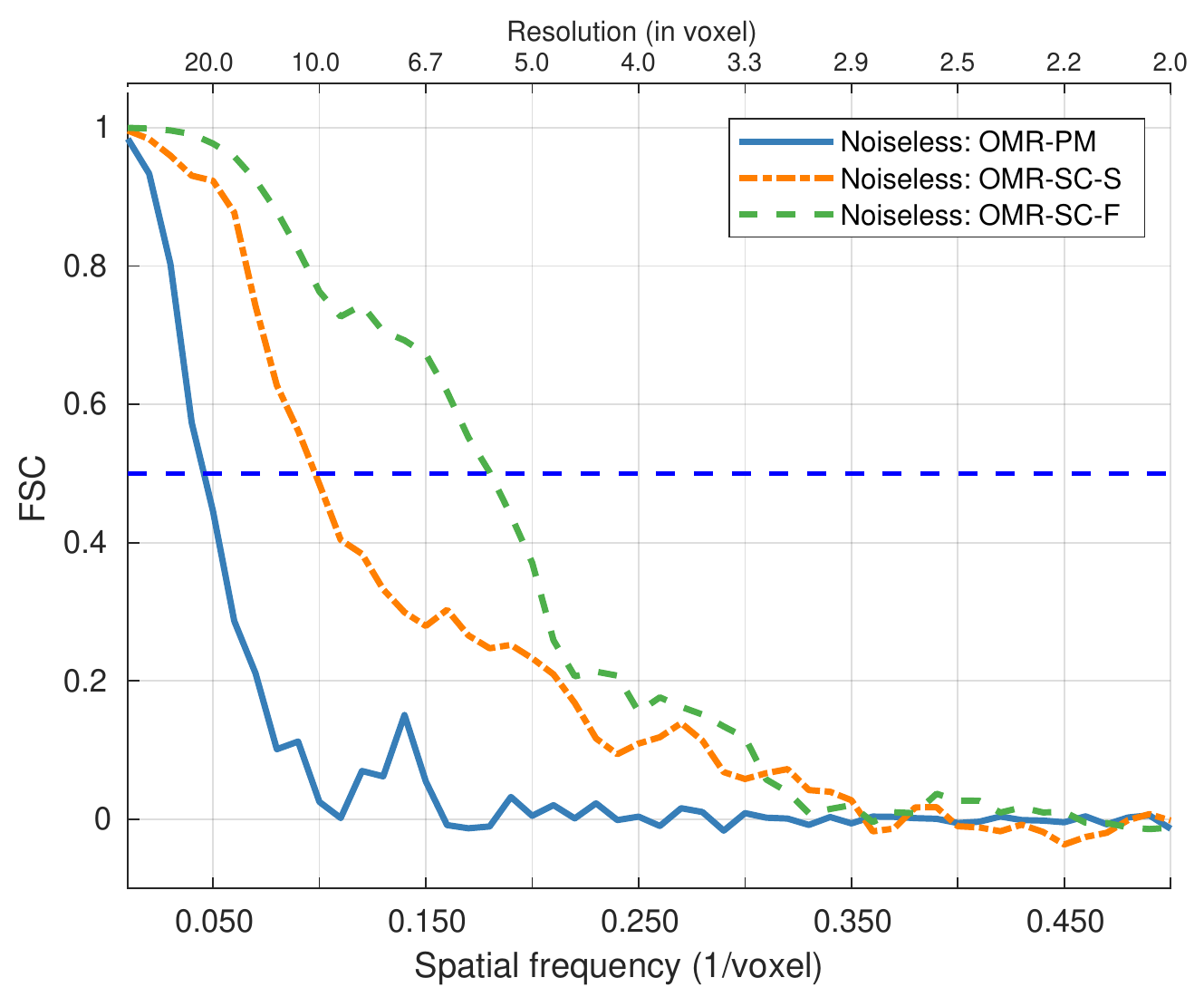}
\hspace{0.03\textwidth}
\includegraphics[width=0.45\textwidth]{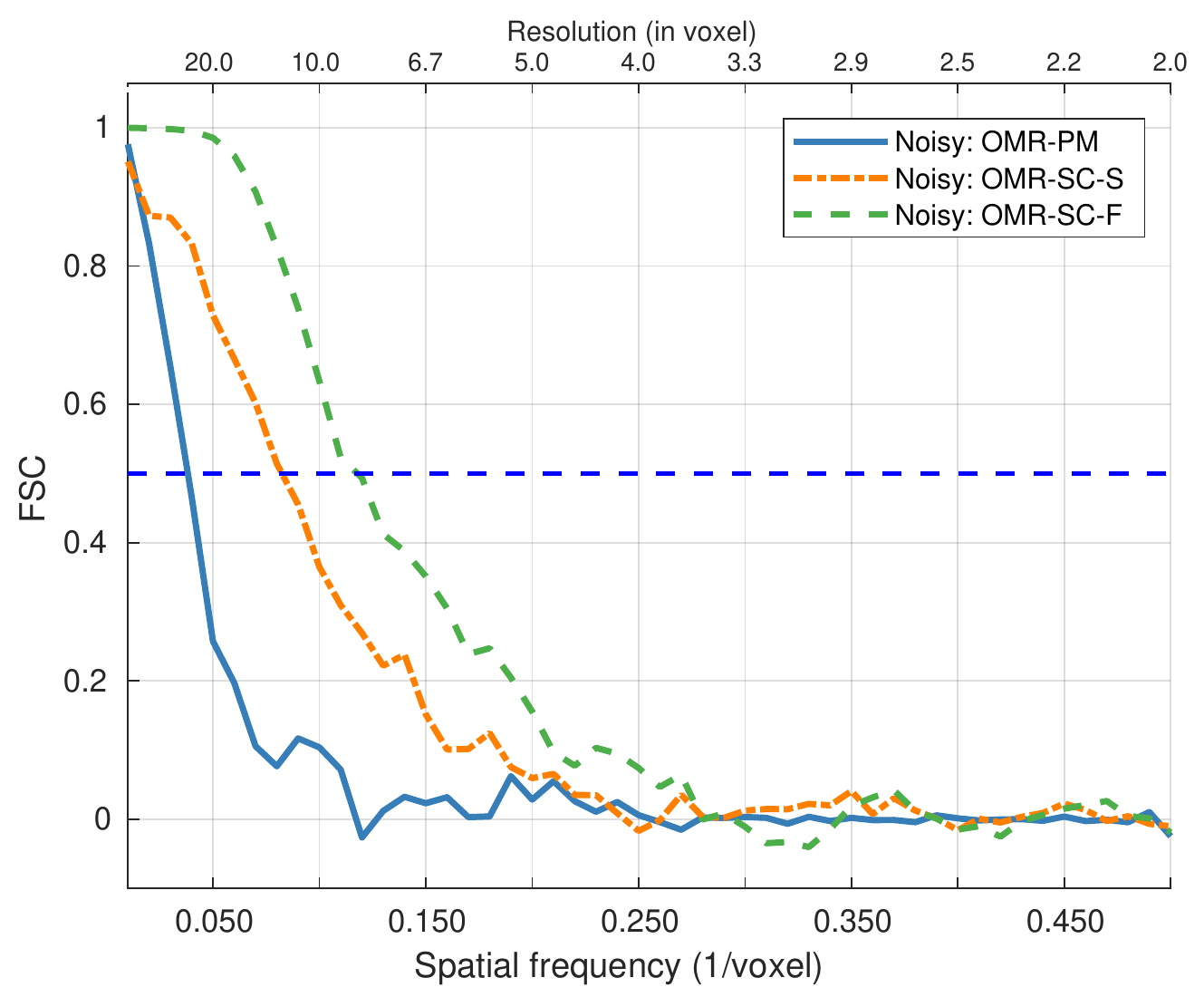}}
\caption{FSC curves of recovered random density maps D1 and D10 using the OMR-PM and OMR-SC approaches in the noiseless case and the noisy case (SNR=$0.1$). The cutoff-threshold of ``$0.5$'' is used to determine the resolution (in voxel). }
\label{fig:compare_fsc_d1_d10}
\end{figure*}

\begin{figure*}[tb]
\centering
\includegraphics[width=\textwidth]{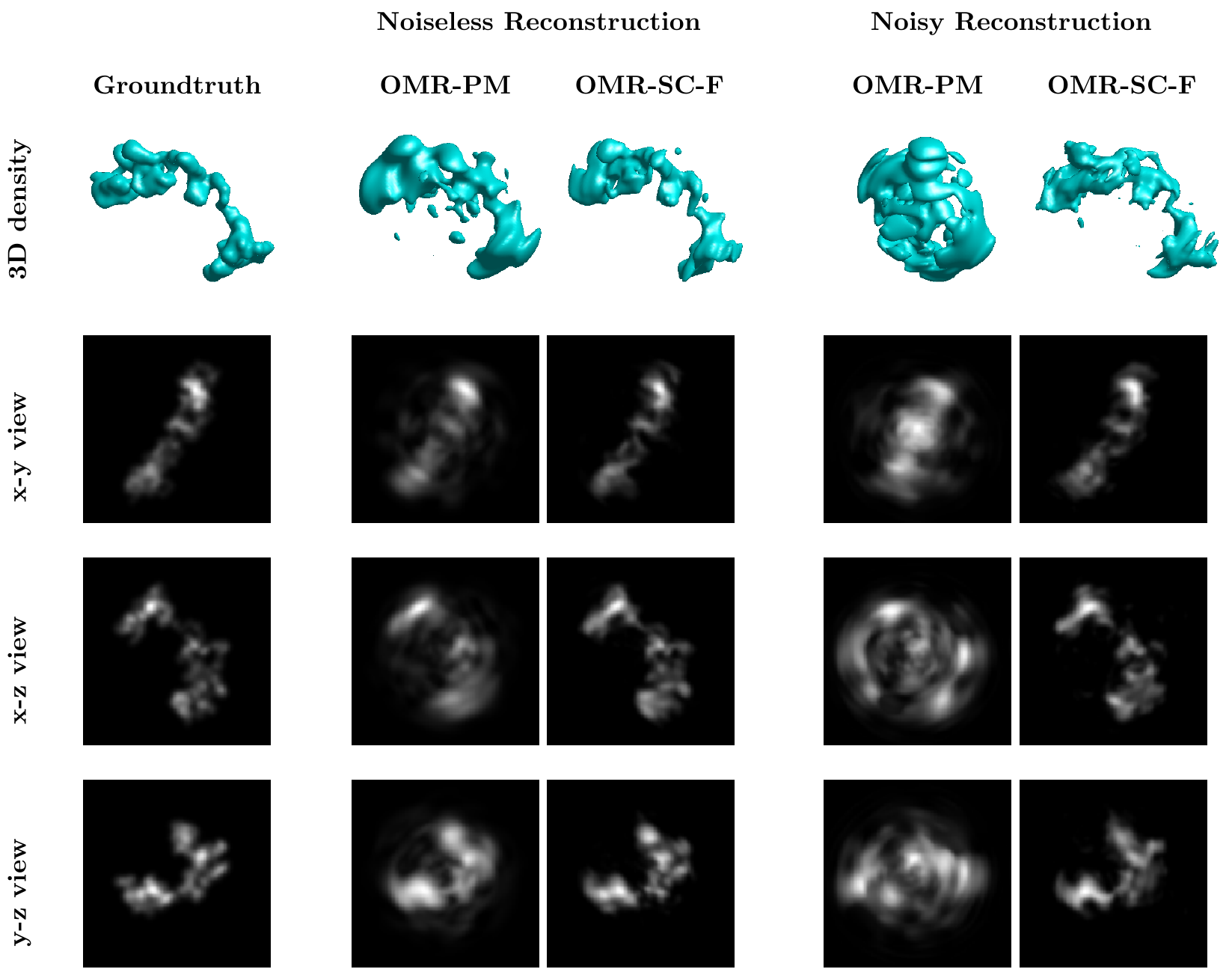}
\caption{Density 1 (D1): reconstructions using the OMR-PM and OMR-SC-F approaches in the noiseless case and the noisy case (SNR=$0.1$).}
\label{fig:compare_reconstruction_1st}
\end{figure*}

\begin{figure*}[tb]
\centering
\includegraphics[width=\textwidth]{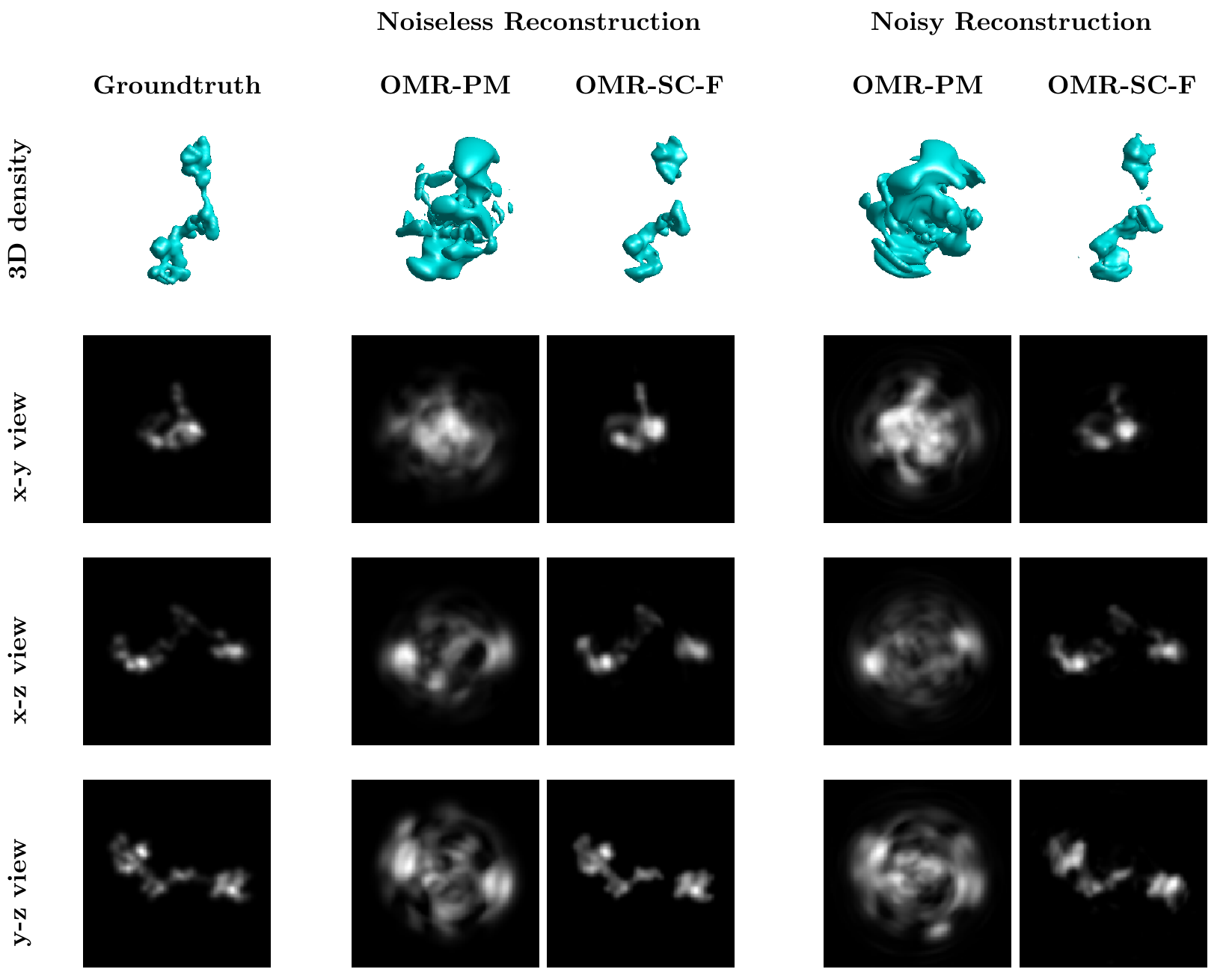}
\caption{Density 10 (D10): reconstructions using the OMR-PM and OMR-SC-F approaches in the noiseless case and the noisy case (SNR=$0.1$).}
\label{fig:compare_reconstruction_10th}
\end{figure*}

Using the first and tenth random densities as representative examples, we plot the FSC curves in Fig. \ref{fig:compare_fsc_d1_d10}, and show the 3D and 2D projection views of the reconstructed density maps in Fig. \ref{fig:compare_reconstruction_1st}-\ref{fig:compare_reconstruction_10th}. Additional figures showing the rest of the FSC curves and reconstructed density maps are given in the Supplementary Material. We can see that the OMR-PM reconstructions are generally much blurrier.
The quantitative correlation coefficients in Table \ref{tab:compare_cc} are also consistent with the visual observation that 3D structures are revealed better by OMR-SC reconstructions in general. We also note that there is always a mismatch between the parametric density map in \eqref{eq:discrete_representation} and the groundtruth, since the Gaussian mixtures in the groundtruth are generally not located on the sampling locations of the Cartesian grid and their variances are different from the one in the bump function \eqref{eq:gauss_basis}. Despite the model mismatch, OMR-SC is still able to recover the density map well.

Due to the nononvexity of the problem, the initialization directly affects the performance of OMR-SC. The proposed initialization scheme draws information from a reference image and spatial radial features. The performances of the \emph{ab initio} models produced by OMR-SC-F are given in Appendix \ref{app:sec:ab_initio_omr_sc_f}. Tables \ref{tab:compare_resolution_random}-\ref{tab:compare_cc_random} compare the performances of the OMR-SC-F approaches with the random (R) and proposed (P) initializations. We can see that the proposed initialization scheme generally leads to better performances except on the noisy recovery of the densities D9 and D10.

\begin{table}[htb]
\caption{Resolutions (in voxel) of recovered random density maps using the OMR-SC-F approaches with the random (R) and proposed (P) initialization schemes (FSC cutoff threshold = 0.5).}
\label{tab:compare_resolution_random}
\centering
\resizebox{\textwidth}{!}{
\begin{tabular}{llcccccccccc}
\toprule
& &D1 &D2 &D3 &D4 &D5 &D6 &D7 &D8 &D9 &D10\\ \midrule
&OMR-SC-F (R) & \bf{6.48} & 10.99 & 6.09 & 7.47 & 6.86 & 4.27 & 11.21 & 9.43 & 12.59 & 5.93  \\
\multirow{-2}{*}{Noiseless} &OMR-SC-F (P) & 6.70 & \bf{9.62} & \bf{5.90} & \bf{7.28} & \bf{5.94} & \bf{4.10} & \bf{10.63} & \bf{7.48} & \bf{8.46} & \bf{5.54}   \\ \midrule
&OMR-SC-F (R) & 9.72 & 18.73 & \bf{7.23} & 13.64 & \bf{8.83} & \bf{5.97} & \bf{16.39} & 14.51 & \bf{9.78} & \bf{7.29} \\
\multirow{-2}{*}{Noisy} &OMR-SC-F (P) & \bf{8.47} & \bf{10.80} & 7.42 & \bf{10.98} & 9.18 & 6.01 & 17.51 & \bf{8.94} & 12.14 & 8.45 \\ \bottomrule
\end{tabular}
}
\end{table}

\begin{table}[htb]
\caption{Correlation coefficients of recovered random density maps using the OMR-SC-F approaches with the random (R) and proposed (P) initialization schemes.}
\label{tab:compare_cc_random}
\centering
\resizebox{\textwidth}{!}{
\begin{tabular}{llcccccccccc}
\toprule
& &D1 &D2 &D3 &D4 &D5 &D6 &D7 &D8 &D9 &D10\\ \midrule
&OMR-SC-F (R) & 0.90 & 0.78 & 0.95 & 0.91 & 0.90 & 0.94 & 0.81 & 0.79 & 0.73 & 0.91    \\
\multirow{-2}{*}{Noiseless} &OMR-SC-F (P) & \bf{0.92} & \bf{0.84} & \bf{0.96} & \bf{0.92} & \bf{0.96} & \bf{0.96} & \bf{0.83} & \bf{0.91} & \bf{0.91} & \bf{0.91}  \\ \midrule
&OMR-SC-F (R) & 0.78 & 0.57 & 0.93 & 0.71 & 0.87 & 0.89 & 0.71 & 0.67 & \bf{0.87} & \bf{0.89}     \\
\multirow{-2}{*}{Noisy} &OMR-SC-F (P) & \bf{0.85} & \bf{0.71} & \bf{0.92} & \bf{0.82} & \bf{0.88} & \bf{0.89} & \bf{0.71} & \bf{0.85} & 0.77 & 0.86 \\ \bottomrule
\end{tabular}
}
\end{table}

The nonnegativity constraints are important to the robustness of OMR-SC, they ensure that the orthogonal matrices reach the consensus on a ``physical'' density map in the spatial domain. Tables \ref{tab:compare_resolution_nonnegative}-\ref{tab:compare_cc_nonnegative} compare the performances of the OMR-SC-F approaches with (w/) and without (w/o) the nonnegativity constraints. In terms of correlation coefficient, we can see that OMR-SC-F achieves  significantly better performances when the nonnegativity constraints are included.

\begin{table}[htb]
\caption{Resolutions (in voxel) of recovered random density maps using the OMR-SC-F approaches with (w/) and without (w/o) nonnegativity constraints (FSC cutoff threshold = 0.5).}
\label{tab:compare_resolution_nonnegative}
\centering
\resizebox{\textwidth}{!}{
\begin{tabular}{llcccccccccc}
\toprule
& &D1 &D2 &D3 &D4 &D5 &D6 &D7 &D8 &D9 &D10\\ \midrule
&OMR-SC-F(w/o) & 8.98 & \bf{9.54} & 11.16 & 10.67 & 8.01 & 6.27 & 15.24 & 9.03 & 16.16 & 6.54 \\
\multirow{-2}{*}{Noiseless} &OMR-SC-F(w/) & \bf{6.70} & 9.62 & \bf{5.90} & \bf{7.28} & \bf{5.94} & \bf{4.10} & \bf{10.63} & \bf{7.48} & \bf{8.46} & \bf{5.54}   \\ \midrule
&OMR-SC-F(w/o) & 9.68 & 13.30 & 8.54 & 68.03 & 21.32 & 10.19 & \bf{16.18} & 20.79 & 17.39 & 10.64 \\
\multirow{-2}{*}{Noisy} &OMR-SC-F(w/) & \bf{8.47} & \bf{10.80} & \bf{7.42} & \bf{10.98} & \bf{9.18} & \bf{6.01} & 17.51 & \bf{8.94} & \bf{12.14} & \bf{8.45} \\ \bottomrule
\end{tabular}
}
\end{table}

\begin{table}[htb]
\caption{Correlation coefficients of recovered random density maps using the OMR-SC-F approaches with (w/) and without (w/o) nonnegativity constraints.}
\label{tab:compare_cc_nonnegative}
\centering
\resizebox{\textwidth}{!}{
\begin{tabular}{llcccccccccc}
\toprule
& &D1 &D2 &D3 &D4 &D5 &D6 &D7 &D8 &D9 &D10\\ \midrule
&OMR-SC-F(w/o) & 0.78 & 0.75 & 0.76 & 0.85 & 0.86 & 0.75 & 0.78 & 0.83 & 0.68 & 0.87   \\
\multirow{-2}{*}{Noiseless} &OMR-SC-F(w/) & \bf{0.92} & \bf{0.84} & \bf{0.96} & \bf{0.92} & \bf{0.96} & \bf{0.96} & \bf{0.83} & \bf{0.91} & \bf{0.91} & \bf{0.91}  \\ \midrule
&OMR-SC-F(w/o) & 0.74 & 0.64 & 0.88 & 0.22 & 0.53 & 0.64 & 0.70 & 0.62 & 0.68 & 0.76     \\
\multirow{-2}{*}{Noisy} &OMR-SC-F(w/) & \bf{0.85} & \bf{0.71} & \bf{0.92} & \bf{0.82} & \bf{0.88} & \bf{0.89} & \bf{0.71} & \bf{0.85} & \bf{0.77} & \bf{0.86} \\ \bottomrule
\end{tabular}
}
\end{table}

\clearpage

\subsubsection{Protein Density Maps}

\begin{figure*}[tb]
\centering
\subfloat[CaS]{
\label{fig:fsc_emd25143}
\includegraphics[width=0.45\textwidth]{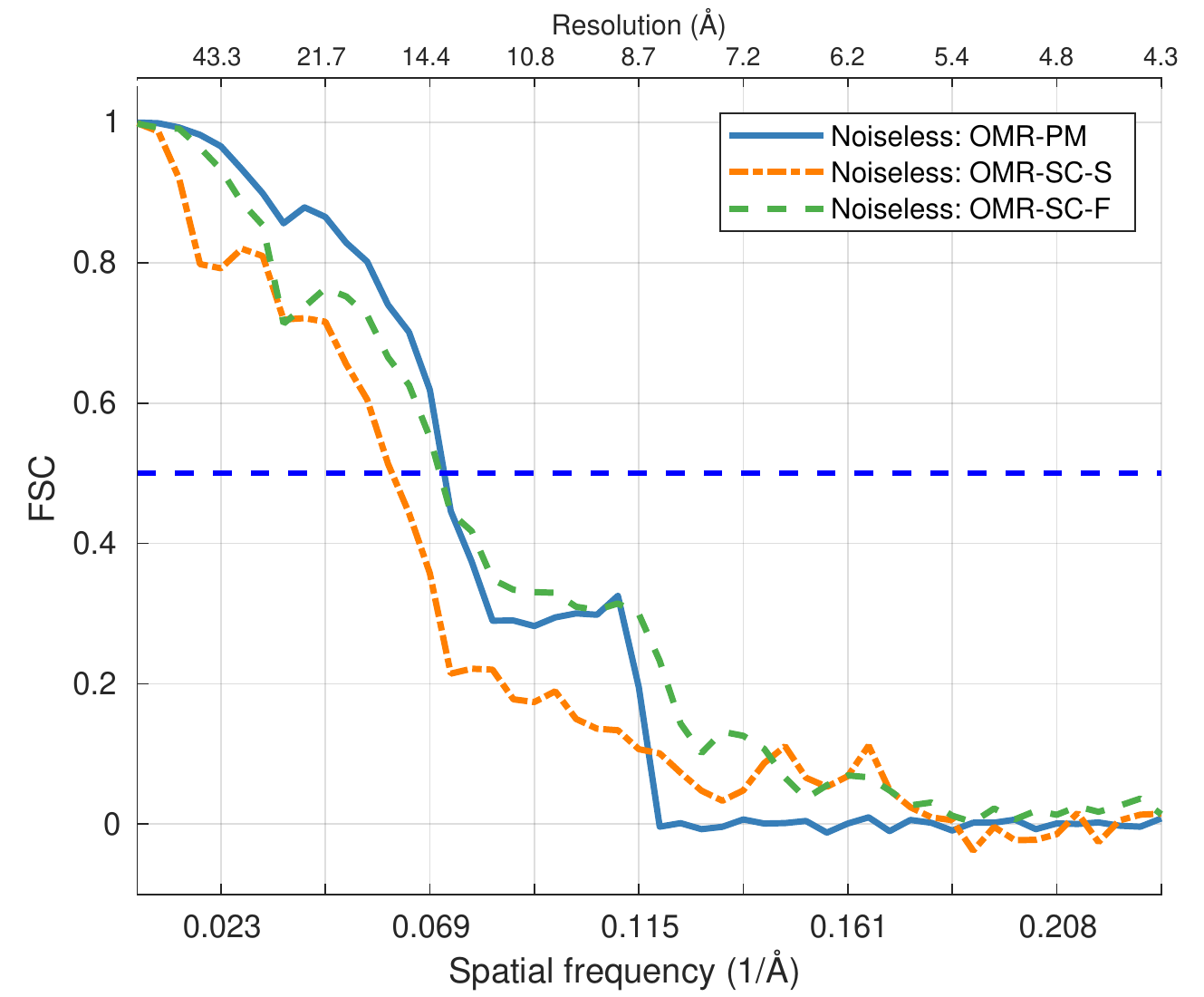}
\hspace{0.03\textwidth}
\includegraphics[width=0.45\textwidth]{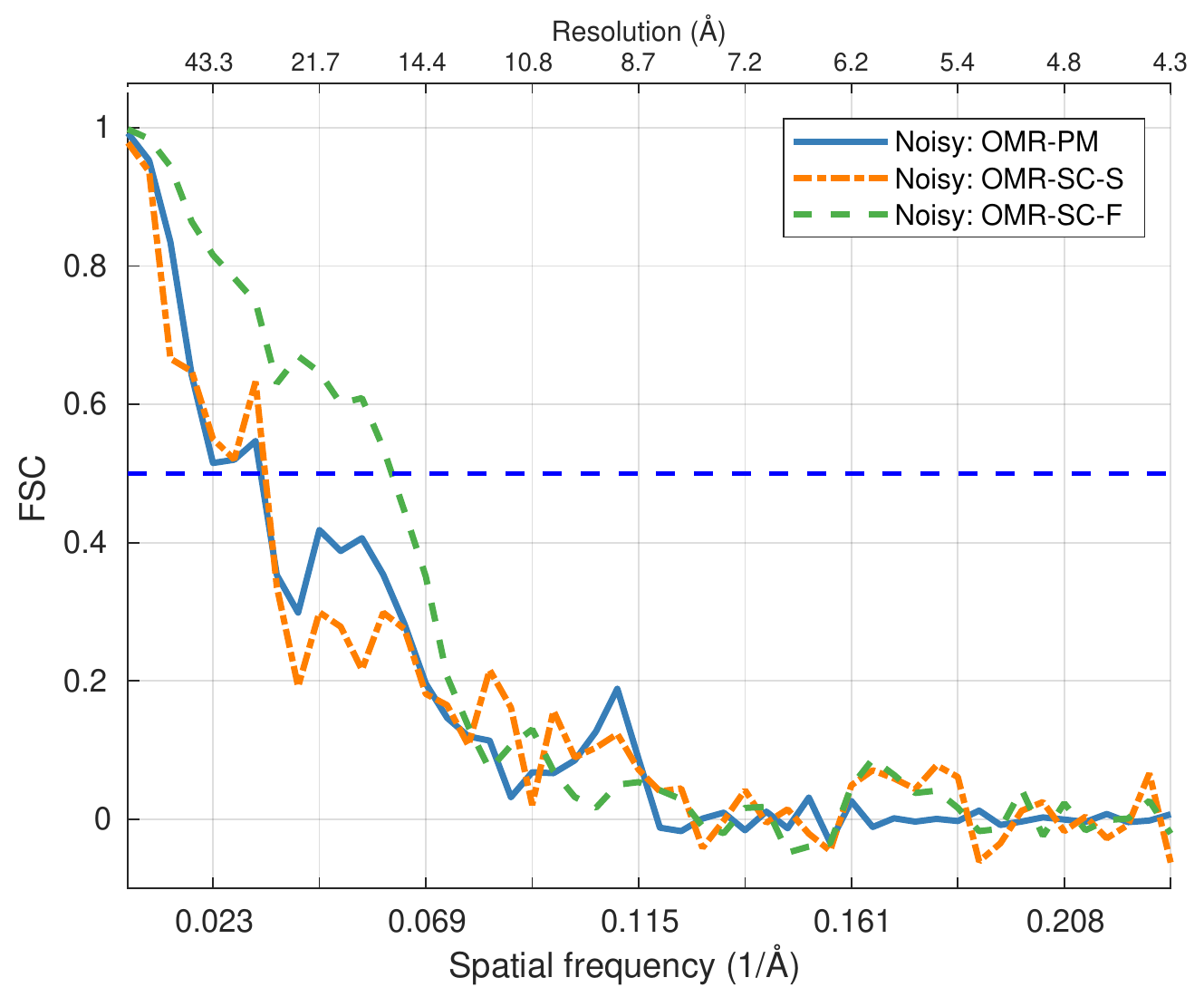}}\\
\subfloat[HJC]{
\label{fig:fsc_hjc}
\includegraphics[width=0.45\textwidth]{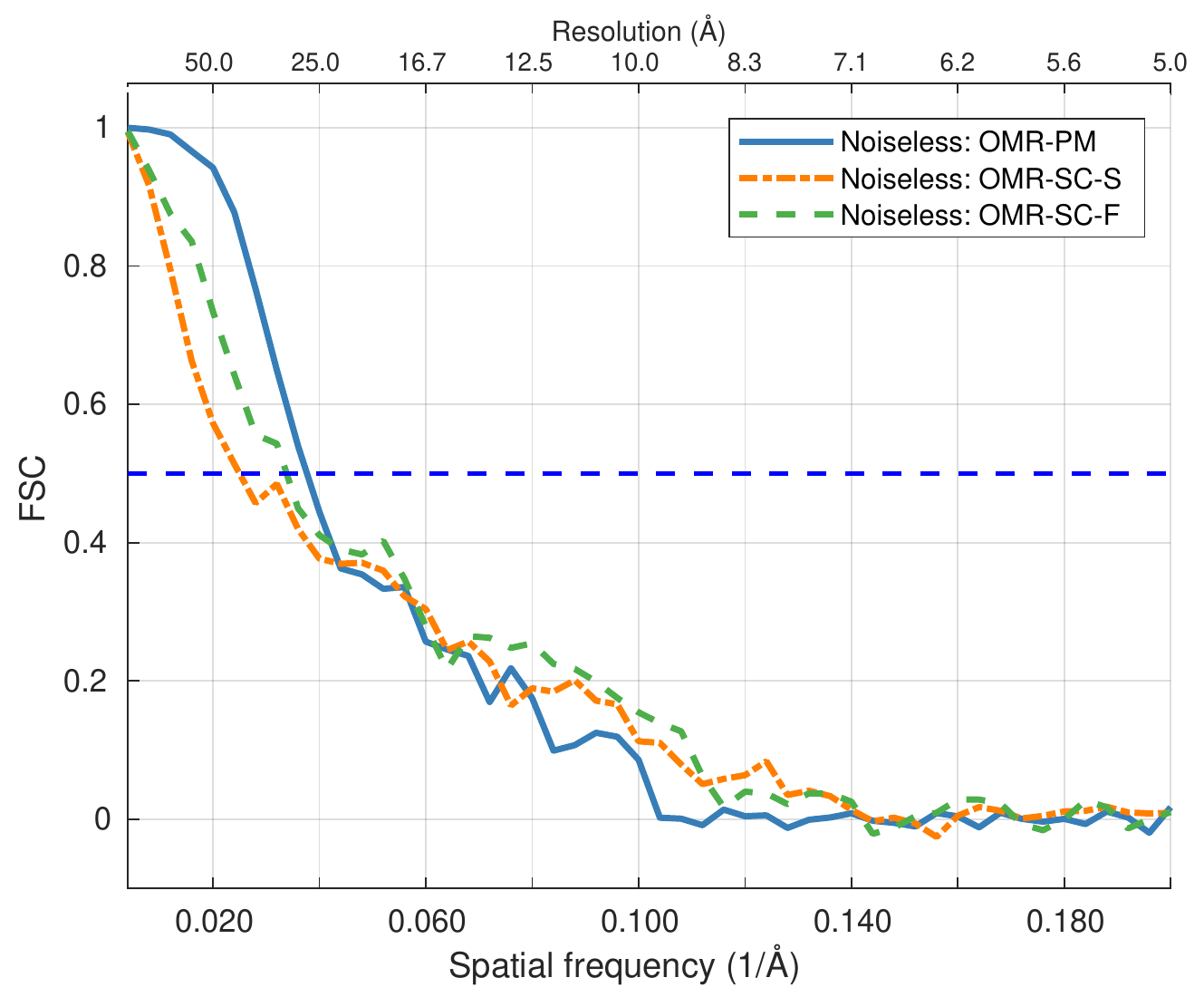}
\hspace{0.03\textwidth}
\includegraphics[width=0.45\textwidth]{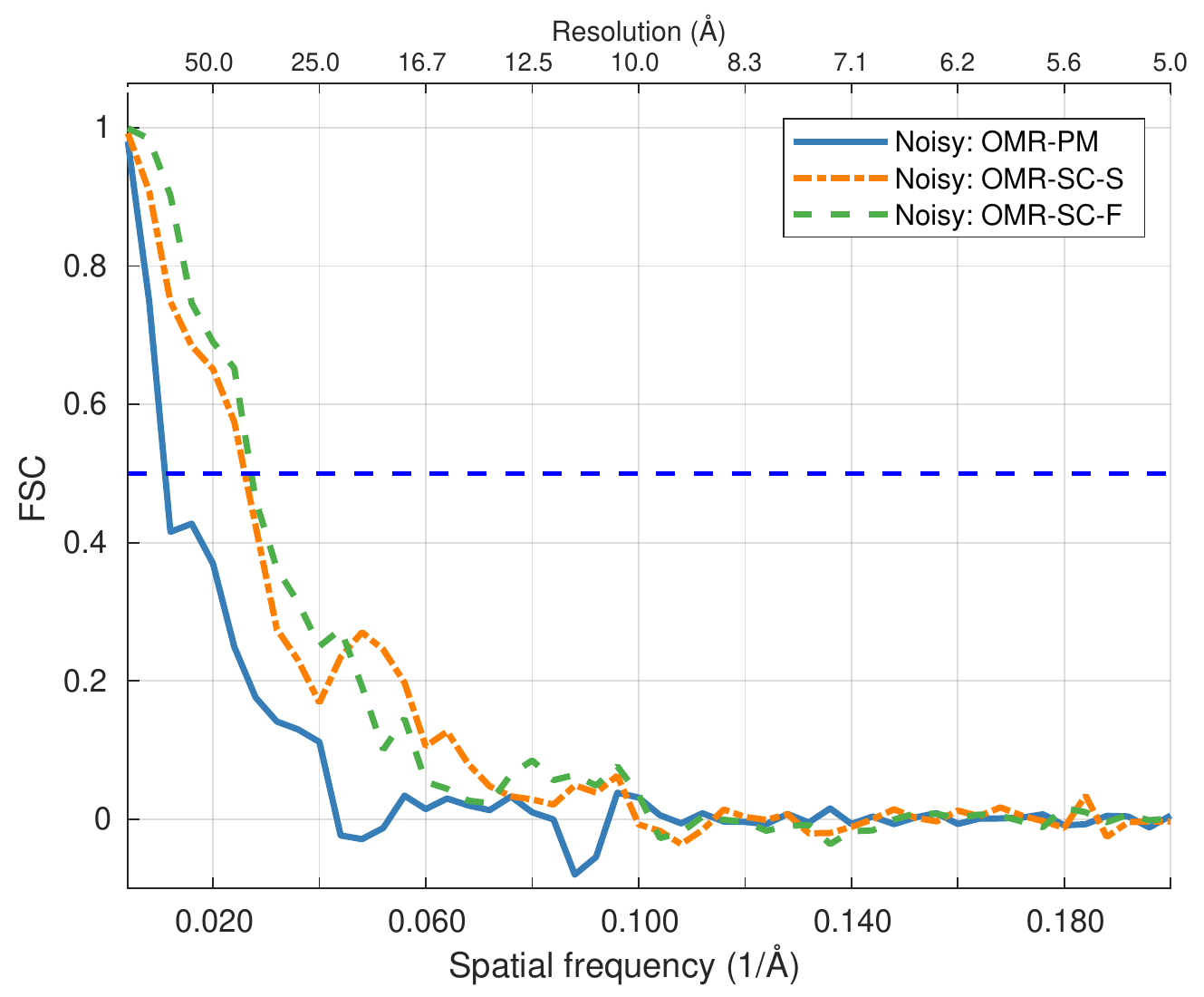}}\\
\subfloat[PTCH1]{
\label{fig:fsc_ptch1}
\includegraphics[width=0.45\textwidth]{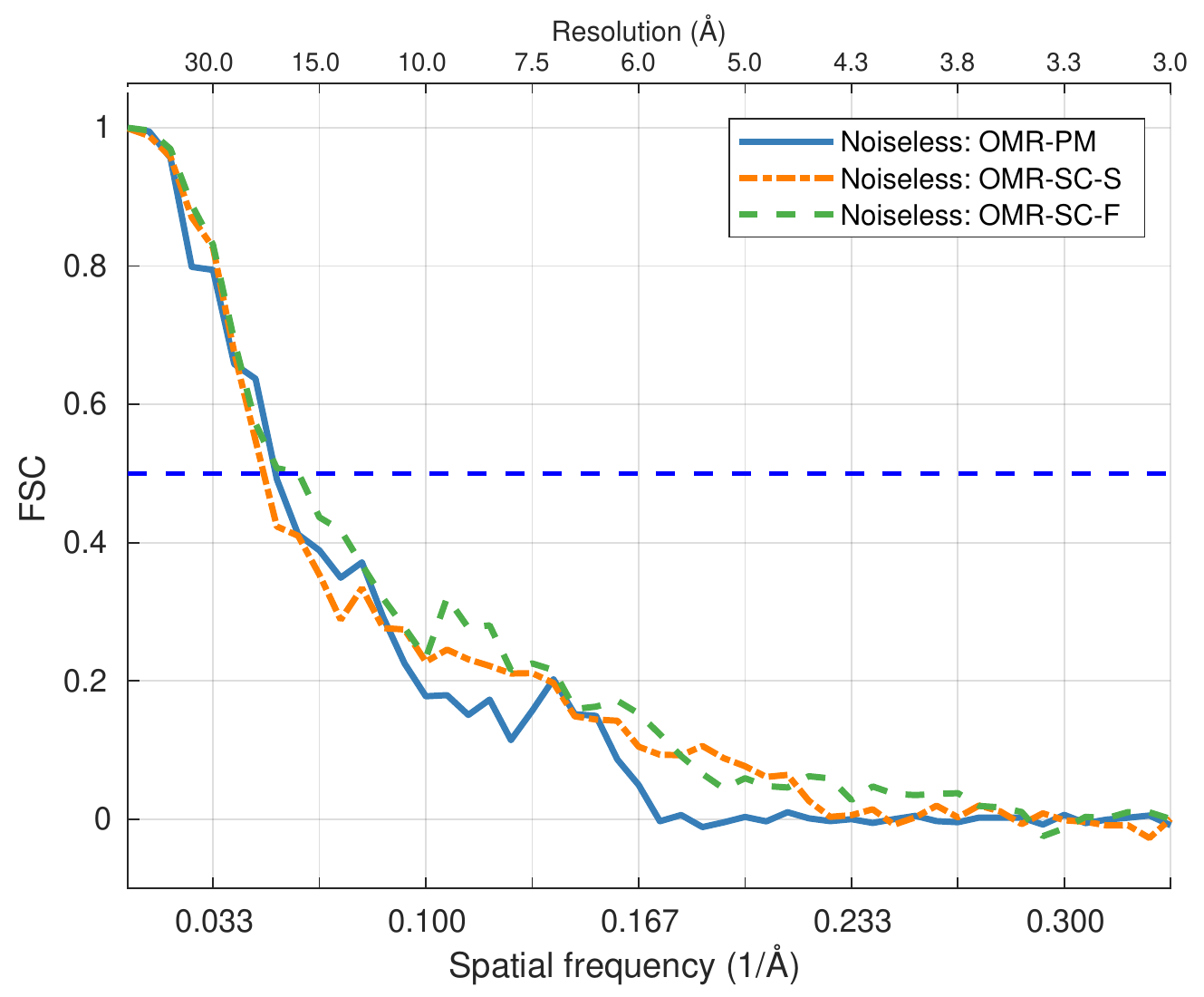}
\hspace{0.03\textwidth}
\label{fig:fsc_noisy_d1}
\includegraphics[width=0.45\textwidth]{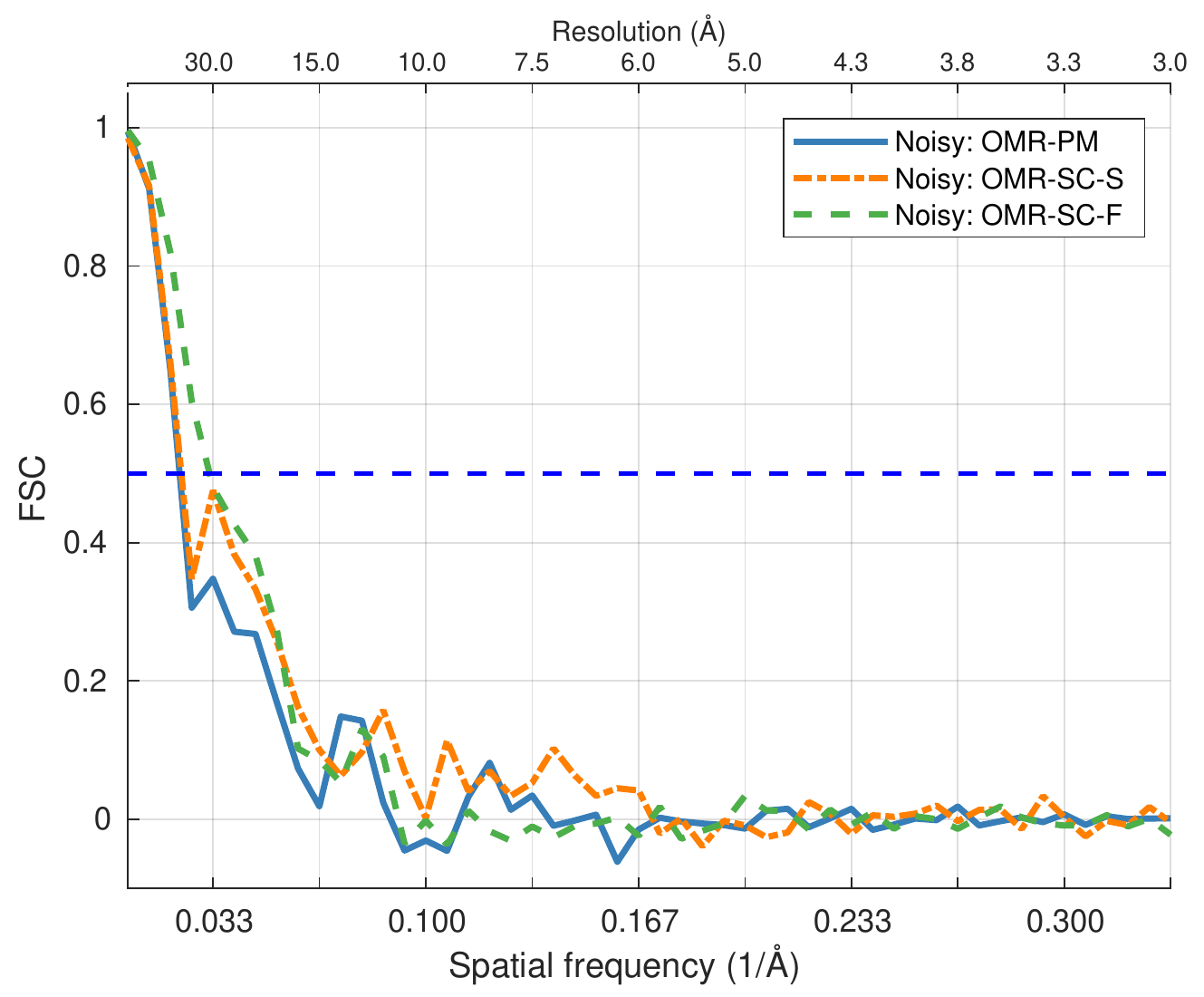}}
\caption{FSC curves of recovered protein density maps using the OMR-PM and OMR-SC approaches in the noiseless case and the noisy case (SNR=$0.1$). The cutoff-threshold of ``$0.5$'' is used to determine the resolution (in \r{A}). }
\label{fig:compare_fsc_emd_25143_hjc_ptch1}
\end{figure*}

\begin{table}[tb]
\caption{Resolutions (\r{A}) and correlation coefficients of recovered protein density maps using the OMR-PM and OMR-SC approaches (FSC cutoff threshold=0.5).}
\label{tab:compare_res_cc_emd_25143_hjc_ptch1}
\centering
\resizebox{\textwidth}{!}{
\begin{tabular}{llcccccccccc}
\toprule
& &\multicolumn{3}{c}{Resolution(\r{A})} &\multicolumn{3}{c}{Correlation Coefficient}\\ \cmidrule(lr){3-5} \cmidrule(lr){6-8}

& & CaS & HJC & PTCH1 & CaS & HJC & PTCH1 \\ \midrule

&OMR-PM & \bf{13.84} &\bf{26.59}  &18.92 &\bf{0.95} &\bf{0.88} &0.86  \\ 
&OMR-SC-S & 16.44 &40.26 & 20.46 &0.89 &0.74 &\bf{0.87} \\
\multirow{-3}{*}{Noiseless} &OMR-SC-F & 14.04 & 29.38 & \bf{16.52} &0.93 &0.78 &0.86 \\ \midrule

&OMR-PM & 29.49 &95.42 & 44.51 &0.77 &0.53 &0.69 \\ 
&OMR-SC-S & 28.78 &38.23 & 43.73 & 0.76 &0.73 &\bf{0.74} \\
\multirow{-3}{*}{Noisy} &OMR-SC-F & \bf{16.15} & \bf{36.82} &\bf{31.38} &\bf{0.87} & \bf{0.73} &0.73 \\ \bottomrule
\end{tabular}
}
\end{table}

\begin{figure*}[tb]
\centering
\includegraphics[width=\textwidth]{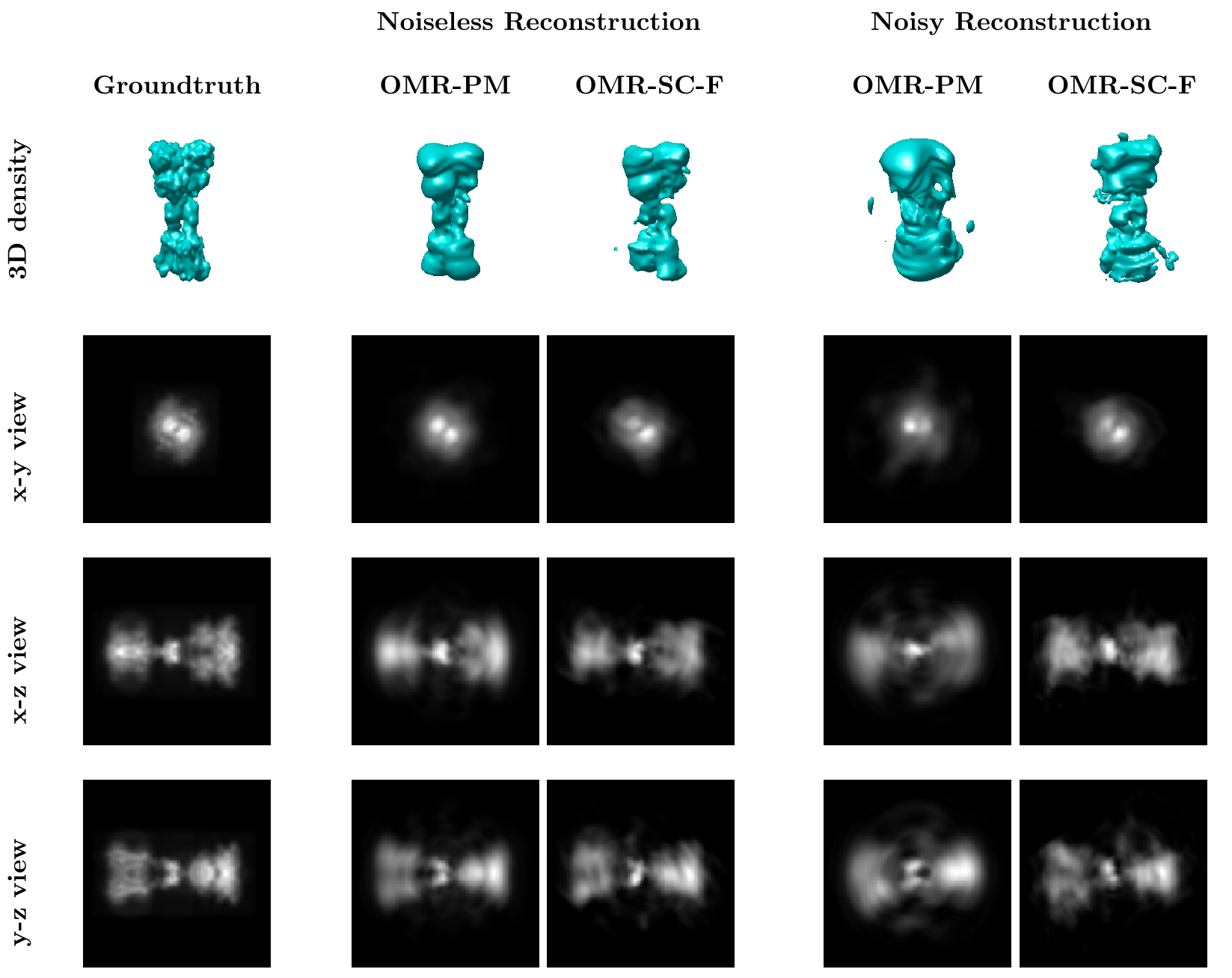}
\caption{Human calcium-sensing receptor (CaS): reconstructions using the OMR-PM and OMR-SC-F approaches in the noiseless case and the noisy case (SNR=$0.1$).}
\label{fig:compare_reconstruction_emd_25143}
\end{figure*}

\begin{figure*}[tb]
\centering
\includegraphics[width=\textwidth]{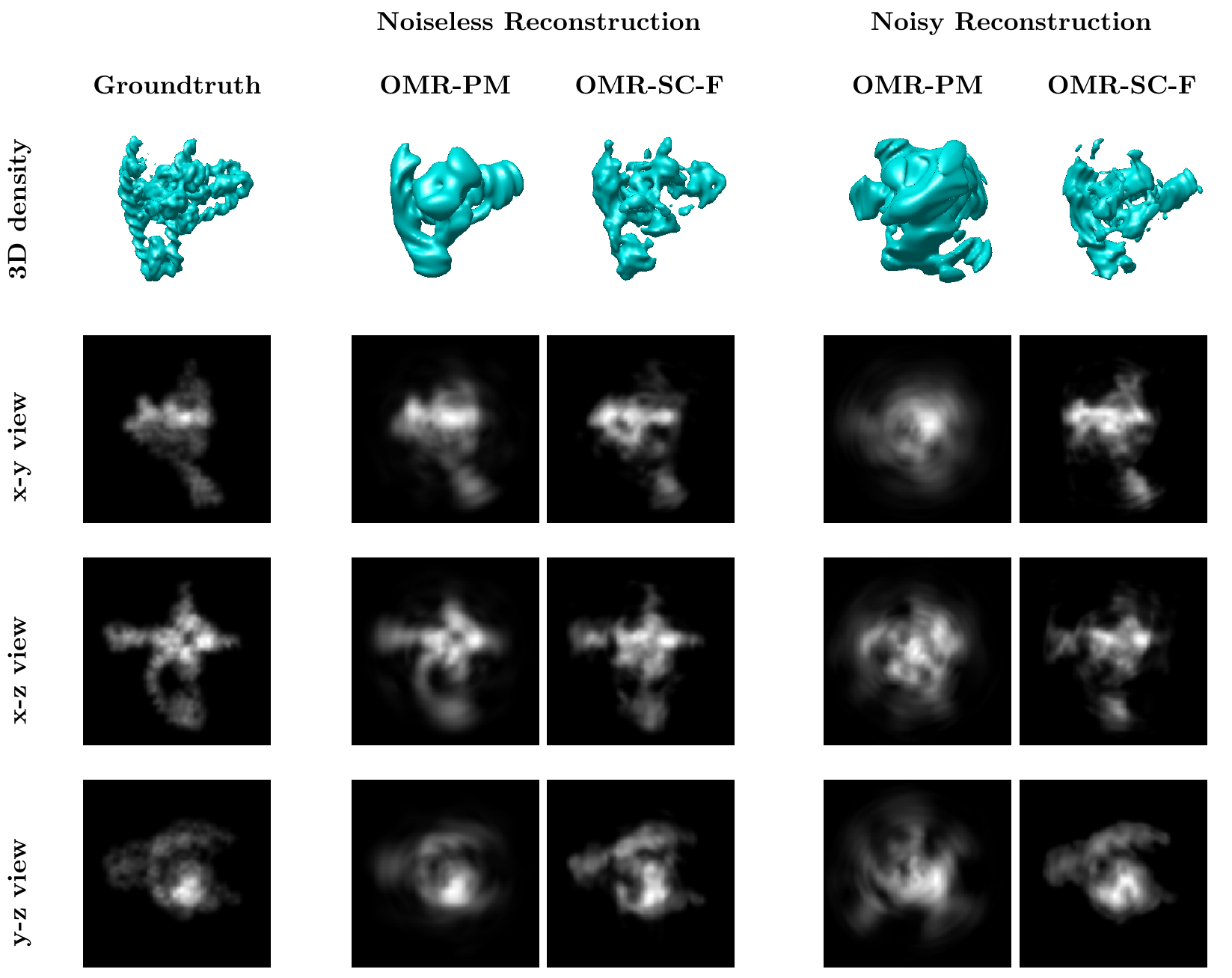}
\caption{Holliday junction complex (HJC): reconstructions using the OMR-PM and OMR-SC-F approaches in the noiseless case and the noisy case (SNR=$0.1$).}
\label{fig:compare_reconstruction_hjc}
\end{figure*}

\begin{figure*}[tb]
\centering
\includegraphics[width=\textwidth]{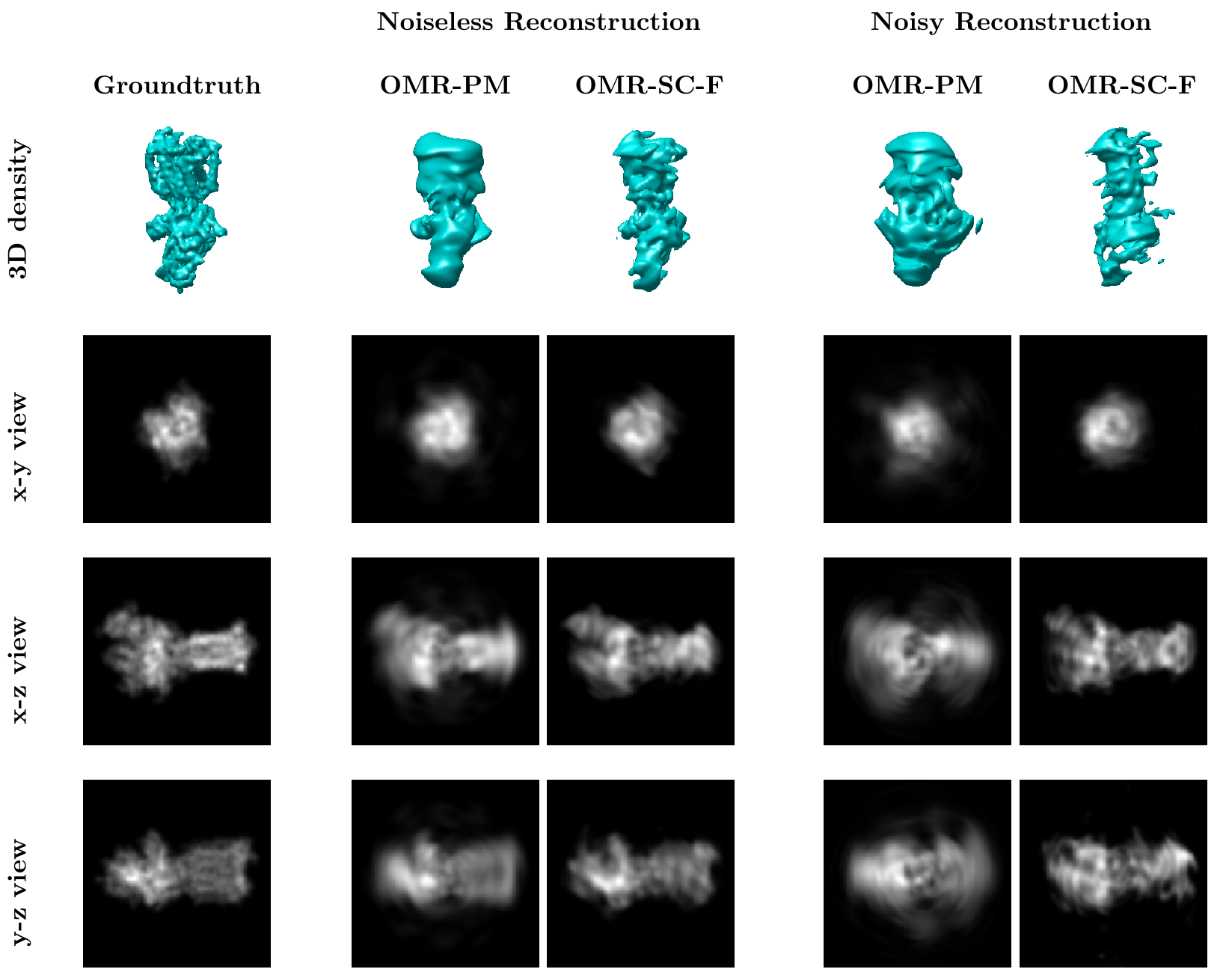}
\caption{Human patched 1 protein (PTCH1): reconstructions using the OMR-PM and OMR-SC-F approaches in the noiseless case and the noisy case (SNR=$0.1$).}
\label{fig:compare_reconstruction_ptch1}
\end{figure*}

Table \ref{tab:compare_res_cc_emd_25143_hjc_ptch1} shows the resolutions and correlation coefficients of recovered protein density maps, Fig. \ref{fig:compare_fsc_emd_25143_hjc_ptch1} shows the FSC curves, and Fig. \ref{fig:compare_reconstruction_emd_25143}-\ref{fig:compare_reconstruction_ptch1} show the recovered density maps in 3D and 2D projection views. 
In the noiseless case, we can see that OMR-PM performs better than OMR-SC on CaS and HJC, they perform equally well on PTCH1. However, in the noisy case OMR-PM becomes unstable, and OMR-SC performs much better on all three density maps. In particular, for the noisy reconstruction of PTCH1, the OMR-SC-S that uses spatial autocorrelations outperforms the OMR-SC-F that uses Fourier autocorrelations.

We also compare the performances of the OMR-SC-F approaches in Tables \ref{tab:compare_res_cc_emd_25143_hjc_ptch1_random}-\ref{tab:compare_res_cc_emd_25143_hjc_ptch1_nonnegative} where we either used random initializations or removed the nonnegativity constraints. The results also show that the proposed initialization and the nonnegativity constraints generally lead to better and more robust performances. In terms of correlation coefficient, although OMR-SC-F without the nonnegativity constraints performs better on the particular noisy recovery of PTCH1, it is not as robust as the OMR-SC-F with the constraints.

\begin{table}[tb]
\caption{Resolutions (\r{A}) and correlation coefficients of recovered protein density maps using the OMR-SC-F approaches with the random (R) and proposed (P) initialization schemes (FSC cutoff threshold=0.5).}
\label{tab:compare_res_cc_emd_25143_hjc_ptch1_random}
\centering
\resizebox{\textwidth}{!}{
\begin{tabular}{llcccccccccc}
\toprule
& &\multicolumn{3}{c}{Resolution(\r{A})} &\multicolumn{3}{c}{Correlation Coefficient}\\ \cmidrule(lr){3-5} \cmidrule(lr){6-8}

& & CaS & HJC & PTCH1 & CaS & HJC & PTCH1 \\ \midrule

&OMR-SC-F (R) & 15.99 &\bf{26.12} & 30.00 &0.85 &0.78 &0.74 \\
\multirow{-2}{*}{Noiseless} &OMR-SC-F (P) & \bf{14.04} & 29.38 & \bf{16.52} &\bf{0.93} &\bf{0.78} &\bf{0.86} \\ \midrule

&OMR-SC-F (R) & 31.78 &41.88 & \bf{28.20} & 0.75 &0.61 &0.71 \\
\multirow{-2}{*}{Noisy} &OMR-SC-F (P) & \bf{16.15} & \bf{36.82} &31.38 &\bf{0.87} & \bf{0.73} &\bf{0.73} \\ \bottomrule
\end{tabular}
}
\end{table}

\begin{table}[tb]
\caption{Resolutions (\r{A}) and correlation coefficients of recovered protein density maps using the OMR-SC-F approaches with (w/) and without (w/o) nonnegativity constraints (FSC cutoff threshold=0.5).}
\label{tab:compare_res_cc_emd_25143_hjc_ptch1_nonnegative}
\centering
\resizebox{\textwidth}{!}{
\begin{tabular}{llcccccccccc}
\toprule
& &\multicolumn{3}{c}{Resolution(\r{A})} &\multicolumn{3}{c}{Correlation Coefficient}\\ \cmidrule(lr){3-5} \cmidrule(lr){6-8}

& & CaS & HJC & PTCH1 & CaS & HJC & PTCH1 \\ \midrule

&OMR-SC-F(w/o) & 14.23 &39.62 & 24.27 &0.92 &0.67 &0.77 \\
\multirow{-3}{*}{Noiseless} &OMR-SC-F(w) & \bf{14.04} & \bf{29.38} & \bf{16.52} &\bf{0.93} &\bf{0.78} &\bf{0.86} \\ \midrule

&OMR-SC-F(w/o) & 17.11 &72.89 & 32.68 & 0.85 &0.58 &\bf{0.77} \\
\multirow{-3}{*}{Noisy} &OMR-SC-F(w/) & \bf{16.15} & \bf{36.82} &\bf{31.38} &\bf{0.87} & \bf{0.73} &0.73 \\ \bottomrule
\end{tabular}
}
\end{table}

\section{Conclusion}
\label{sec:conclusion}
In an effort to expand the applicability of the method of moments in unknown-view tomography, we proposed spatial radial and autocorrelation features that can be expressed as linear and quadratic functionals of the density map. 
Via a spherical Bessel transform, the spatial autocorrelations provide a closed-form link between the Fourier autocorrelations and the sought density.
Prior work noted that (under realistic assumptions) the autocorrelation features determine the spherical harmonic coefficients of the density map up to a set of unknown orthogonal matrices. But that prior work only attempted to recover the a priori uncoupled orthogonal matrices. Due in part to the functional forms of the used features, recovering the correct matrices that ``agree'' on a proper density estimate that satisfies the problem-specific constraints has been challenging. The challenge has been exacerbated by the nonconvexity of the involved optimization problems.

In this paper we addressed the first challenge by the newly-proposed closed-form spatial features, and we greatly alleviated the second, nonconvexity challenge, by designing an efficiently-computable initial density. We then formulated a joint recovery of the density map and the orthogonal matrices in an alternating fashion, constraining the density map to be 1) nonnegative, 2) compactly supported, 3) of correct total mass, and 4) consistent with the denoised reference projection. Experiments show that the proposed OMR with spatial consensus (OMR-SC) is more robust and performs much better in the the presence of noise than the previous state-of-the-art OMR by projection matching (OMR-PM).

The main drawback of the proposed OMR-SC is that it needs a large number of projection images to counter the noise and extract features of sufficient quality, especially at high spherical-harmonic degrees. A naive attempt to use off-the-shelf (or even bespoke) denoisers fails: since current image denoising methods only aim to improve image quality, they introduce typically unknown bias which makes them unsuitable for feature extraction. 
While we could improve the situation by denoising the extracted features, it is an interesting and important challenge to design schemes that simultaneously extract features and denoise projection images. Further, like prior work, OMR-SC relies on first- and second-order moments. Using higher-order moments is challenging since for a fixed SNR noise amplification and the related computational complexity grow exponentially with the moment order.

Another major challenge is the non-convexity introduced by the orthogonality constraints. We proposed an effective computational strategy that greatly reduces the impact of this challenge: a data-dependent, efficiently-computable initialization reminiscent of spectral initialization in phase retrieval with random measurements. Importantly, the reconstructions for multiple initializations can be carried out in parallel, and we empirically found that a small number of initialization ($\approx 10$) suffices to obtain a good reconstruction. Alas, unlike for phase retrieval with random measurements, existence of theoretical guarantees for deterministic moments used in our approach remains an open question.

\section*{Acknowledgement}
We are grateful to Yifeng Fan for his help in providing the code for denoising the images as linear projection features.

\clearpage

\begin{appendices}

\section{Noisy Feature Debiasing}
\label{app:sec:debias_features}
Let $P(x,y)$ denote the noiseless projection image of size $G\times G$, and the noisy image $S(x,y)=P(x,y)+\epsilon(x,y)$. We compute the frequency component $\widehat{S}(k_x,k_y)$ using the discrete (non-uniform) Fourier transform
\begin{align}
\label{eq:nufft}
\begin{split}
    \widehat{S}(k_x,k_y)&=\sum_x\sum_y \exp\left(-\mathfrak{i}\cdot[k_x\ k_y]\left[\begin{array}{c} x\\ y\end{array}\right]\right)\cdot \Big(P(x,y)+\epsilon(x,y)\Big)\\
    &=\vf_{\vk}^T(\vp+\boldsymbol\epsilon)\,,
\end{split}
\end{align}
where $\vp$ is the vectorized projection image $P(x,y)$, $\boldsymbol\epsilon$ is the vectorized noise $\epsilon(x,y)$, and $\vf_\vk$ is the vectorized exponential term to compute $\widehat{S}(\vk)$.

\subsection{Fourier Autocorrelation Feature}
\label{app:subsec:debias_autocorrelation_features}

When the number of projection images goes to infinity, i.e. $N\rightarrow\infty$, according to the strong law of large numbers, the sample average converges almost surely to the expected value. The average with respect to the projection images in \eqref{eq:C_k1_k2_psi} can thus be replaced by an expectation operator. The bias term $\zeta(k_1,k_2,\psi)$ is then
\begin{align}
\begin{split}
    \zeta(k_1,k_2,\psi)=&\mathbb{E}\left[\frac{1}{2\pi}\int_0^{2\pi}\vf_{\vk_1}^T(\vp+\boldsymbol\epsilon)\cdot \vf_{\vk_2}^*{}^T(\vp+\boldsymbol\epsilon)\ d\varphi_{\vk_1}\right]\\
    &-\mathbb{E}\left[\frac{1}{2\pi}\int_0^{2\pi}\vf_{\vk_1}^T\vp\cdot \vf_{\vk_2}^*{}^T\vp\ d\varphi_{\vk_1}\right]\,,
\end{split}
\end{align}
where $\varphi_{\vk_1} = \measuredangle \vk_1$ and $\varphi_{\vk_1} + \psi=\measuredangle \vk_2$. The noise $\boldsymbol\epsilon_n$ is modeled as additive white Gaussian noise with variance $\sigma_\epsilon^2$ \cite{CryoEM:Bendory:2020}. Since the AWGN noise $\boldsymbol\epsilon$ is independent from the projection image $\vp$ and $\mathbb{E}[\boldsymbol\epsilon]=\boldsymbol 0$, we have
\begin{align}
\label{eq:expectation_fk1_fk2}
\begin{split}
    &\mathbb{E}\left[\frac{1}{2\pi}\int_0^{2\pi}\widehat{\rho}(\vk_1)\cdot\widehat{\rho}^*(\vk_2)\ d\varphi_{\vk_1}\right]\\
    &=\frac{1}{2\pi}\int_0^{2\pi}\mathbb{E}\left[\vf_{\vk_1}^T\vp\cdot {\vf_{\vk_2}^*}^T\vp\right]\ d\varphi_{\vk_1} + \frac{1}{2\pi}\int_0^{2\pi}\mathbb{E}\left[\vf_{\vk_1}^T\vp\cdot{\vf_{\vk_2}^*}^T\boldsymbol\epsilon\right]\ d\varphi_{\vk_1} \\
    &\quad + \frac{1}{2\pi}\int_0^{2\pi}\mathbb{E}\left[\vf_{\vk_1}^T\boldsymbol\epsilon\cdot{\vf_{\vk_2}^*}^T\vp\right]\ d\varphi_{\vk_1} + \frac{1}{2\pi}\int_0^{2\pi}\mathbb{E}\left[\vf_{\vk_1}^T\boldsymbol\epsilon\cdot{\vf_{\vk_2}^*}^T\boldsymbol\epsilon\right]\ d\varphi_{\vk_1}\\
    &=\frac{1}{2\pi}\int_0^{2\pi}\mathbb{E}\left[\vf_{\vk_1}^T\vp\cdot {\vf_{\vk_2}^*}^T\vp\right]\ d\varphi_{\vk_1} + \frac{1}{2\pi}\int_0^{2\pi}\vf^T_{\vk_1}\mathbb{E}[\vp]\cdot \vf^T_{\vk_2}\mathbb{E}[\boldsymbol\epsilon]\ d\varphi_{\vk_1} \\
    &\quad+ \frac{1}{2\pi}\int_0^{2\pi}\vf^T_{\vk_1}\mathbb{E}[\boldsymbol\epsilon]\cdot \vf^T_{\vk_2}\mathbb{E}[\vp]\ d\varphi_{\vk_1} + \frac{1}{2\pi}\int_0^{2\pi}\mathbb{E}\left[\vf_{\vk_1}^T\boldsymbol\epsilon\cdot{\vf_{\vk_2}^*}^T\boldsymbol\epsilon\right]\ d\varphi_{\vk_1}\\
    &=\frac{1}{2\pi}\int_0^{2\pi}\mathbb{E}\left[\vf_{\vk_1}^T\vp\cdot {\vf_{\vk_2}^*}^T\vp\right]\ d\varphi_{\vk_1}+\frac{1}{2\pi}\int_0^{2\pi}\mathbb{E}\left[\vf_{\vk_1}^T\boldsymbol\epsilon\cdot{\vf_{\vk_2}^*}^T\boldsymbol\epsilon\right]\ d\varphi_{\vk_1}\,.
\end{split}
\end{align}
From the above \eqref{eq:expectation_fk1_fk2}, we can get the following bias term
\begin{align}
\begin{split}
    &\frac{1}{2\pi}\int_0^{2\pi}\mathbb{E}\left[\vf_{\vk_1}^T\boldsymbol\epsilon\cdot{\vf_{\vk_2}^*}^T\boldsymbol\epsilon\right]\ d\varphi_{\vk_1}\\
    &=\frac{1}{2\pi}\int_0^{2\pi}\mathbb{E}\left[\left(\sum_if_{\vk_1}(i)\cdot\epsilon_i\right)\cdot\left(\sum_jf_{\vk_2}^*(j)\cdot\epsilon_j\right)\right]\ d\varphi_{\vk_1}\\
    &=\frac{1}{2\pi}\int_0^{2\pi}\mathbb{E}\left[\sum_if_{\vk_1}(i)f_{\vk_2}^*(i)\cdot\epsilon_i^2+\sum_{i\neq j}f_{\vk_1}(i)f_{\vk_2}^*(j)\cdot\epsilon_i\epsilon_j\right]\ d\varphi_{\vk_1}\\
    &=\frac{1}{2\pi}\int_0^{2\pi}\sum_if_{\vk_1}(i)f_{\vk_2}^*(i)\cdot\mathbb{E}\left[\epsilon_i^2\right]+\sum_{i\neq j}f_{\vk_1}(i)f_{\vk_2}^*(j)\cdot\mathbb{E}\left[\epsilon_i\epsilon_j\right]\ d\varphi_{\vk_1}\\
    &=\mathbb{E}[\epsilon^2]\cdot\sum_i\frac{1}{2\pi}\int_0^{2\pi}f_{\vk_1}(i)f_{\vk_2}^*(i)\ d\varphi_{\vk_1}\\
    &=\mathbb{E}[\epsilon^2]\cdot\sum_x\sum_y\frac{1}{2\pi}\int_0^{2\pi}\exp\left(-\mathfrak{i}\cdot(\vk_1-\vk_2)^T\left[\begin{array}{c} x\\ y\end{array}\right]\right)\ d\varphi_{\vk_1}\,.
\end{split}
\end{align}
Since the noise $\epsilon$ is additive white Gaussian with variance $\sigma_\epsilon^2$, we have that $\mathbb{E}[\epsilon^2]=\sigma_\epsilon^2$. The bias term of the autocorrelation function is then
\begin{align}
    \zeta(k_1,k_2,\psi)=\sigma_\epsilon^2\cdot\sum_x\sum_y\frac{1}{2\pi}\int_0^{2\pi}\exp\left(-\mathfrak{i}\cdot(\vk_1-\vk_2)^T\left[\begin{array}{c} x\\ y\end{array}\right]\right)\ d\varphi_{\vk_1}\,.
\end{align}

\subsection{Fourier Radial Feature}
\label{app:subsec:debias_radial_features}
When the number of projections images goes to infinity, the average with respect to the projection images in \eqref{eq:radial_integration_feature_step_0} can also be replaced by the expectation operator. We then have
\begin{align}
\begin{split}
    \mathbb{E}\left[\frac{1}{2\pi}\int_0^{2\pi}\widehat{\rho}(\vk)\ d\varphi_{\vk}\right]&=\frac{1}{2\pi}\int_0^{2\pi}\mathbb{E}\left[\vf_{\vk}^T(\vp+\boldsymbol\epsilon)\right]\ d\varphi_{\vk} \\
    &=\frac{1}{2\pi}\int_0^{2\pi}\mathbb{E}\left[\vf_{\vk}^T\vp\right]\ d\varphi_{\vk} + \frac{1}{2\pi}\int_0^{2\pi}\sum_if_{\vk}(i)\mathbb{E}\left[\epsilon\right]\ d\varphi_{\vk} \\
    &=\frac{1}{2\pi}\int_0^{2\pi}\mathbb{E}\left[\vf_{\vk}^T\vp\right]\ d\varphi_{\vk}\,.
\end{split}
\end{align}
We can see that we do not need to perform debiasing on $M(k)$ since $\mathbb{E}[\epsilon]=0$\,.

\section{Autocorrelation Function in the Spatial Domain}
\label{app:sec:autocorrelation_fun_spatial}

Using \eqref{eq:B_lm_r}, the spatial autocorrelation feature $E_l(r_1,r_2)$ can be rewritten as
\begin{align}
\begin{split}
    E_l(r_1,r_2) &= \sum_{m=-l}^{l}B_{lm}(r_1)\cdot B_{lm}(r_2)\\
    &=\frac{1}{r_1^2r_2^2}\iiint\iiint \rho(\tilde\vr_1)\cdot\delta(r_1-\tilde{r}_1)\cdot\rho(\tilde\vr_2)\cdot\delta(r_2-\tilde{r}_2)\\
    &\quad\quad\quad\times\sum_{m=-l}^l Y_{lm}(\theta_{\tilde\vr_1},\varphi_{\tilde\vr_1})\cdot Y_{lm}(\theta_{\tilde\vr_2},\varphi_{\tilde\vr_2})\ d\tilde\vr_1 d\tilde\vr_2\\
    &=\frac{1}{r_1^2r_2^2}\iiint\iiint\rho(\tilde\vr_1)\cdot\delta(r_1-\tilde{r}_1)\cdot\rho(\vr_2)\cdot\delta(r_2-\tilde{r}_2)\\
    &\quad\quad\quad\times\frac{2l+1}{4\pi}P_l(\cos\psi_{\tilde\vr_1,\tilde\vr_2})\ d\tilde\vr_1d\tilde\vr_2\,,
\end{split}
\end{align}
where $\psi_{\tilde\vr_1,\tilde\vr_2}$ is the angle between $\tilde\vr_1$ and $\tilde\vr_2$. 

The autocorrelation function $E(r_1,r_2,\psi)$ can be computed from $E_l(r_1,r_2)$ using the completeness of Legendre polynomials as follows
\begin{align}
\begin{split}
    E(r_1,r_2,\psi) &= \sum_{l=0}^\infty E_l(r_1,r_2)\cdot P_l(\cos\psi)\\
    &=\frac{1}{r_1^2r_2^2}\iiint\iiint\rho(\tilde\vr_1)\cdot\delta(r_1-\tilde{r}_1)\cdot\rho(\tilde\vr_2)\cdot\delta(r_2-\tilde{r}_2)\\
    &\quad\quad\quad\times\frac{1}{2\pi}\sum_{l=0}^\infty \frac{2l+1}{2}P_l(\cos\psi_{\tilde\vr_1,\tilde\vr_2})\cdot P_l(\cos\psi)\ d\tilde\vr_1 d\tilde\vr_2\\
    &=\frac{1}{r_1^2r_2^2}\iiint\iiint\rho(\tilde\vr_1)\cdot\delta(r_1-\tilde{r}_1)\cdot\rho(\tilde\vr_2)\cdot\delta(r_2-\tilde{r}_2)\\
    &\quad\quad\quad\times\frac{1}{2\pi}\cdot\delta(\cos\psi_{\tilde\vr_1,\tilde\vr_2}-\cos\psi)\ d\tilde\vr_1 d\tilde\vr_2\\
    &=\frac{1}{2\pi r_1^2 r_2^2}\iiint\iiint\rho(\tilde\vr_1)\cdot\delta(r_1-\tilde{r}_1)\cdot\rho(\tilde\vr_2)\cdot\delta(r_2-\tilde{r}_2)\\
    &\quad\quad\quad\times\delta(\psi_{\tilde\vr_1,\tilde\vr_2}-\psi)\ d\tilde\vr_1 d\tilde\vr_2\,,
\end{split}
\end{align}
where $\psi\in[0,\pi]$.

\section{Spherical Harmonic Expansion of the Gaussian Basis Function}
\label{app:sec:she_gaussian}

The Gaussian basis function at the sampling location $\boldsymbol\mu_d$ is
\begin{align}
    h(\vr-\boldsymbol\mu_d)=\frac{1}{(2\pi)^{\frac{3}{2}}\sigma^3}\exp\left(-\frac{1}{2}\frac{\|\vr-\boldsymbol \mu_d\|_2^2}{\sigma^2}\right)\,.
\end{align}
It can be expanded in real spherical harmonics:
\begin{align}
    h(\vr-\boldsymbol\mu_d)=\sum_{l=0}^\infty\sum_{m=-l}^lg_{lm}(r,\boldsymbol\mu_d)\cdot Y_{lm}(\theta_\vr,\varphi_\vr)\,,
\end{align}
where $Y_{lm}(\theta_\vr,\varphi_\vr)$ is the real spherical harmonic function of degree $l$ and order $m$, and $g_{lm}(r,\boldsymbol\mu_d)$ is the expansion coefficient:
\begin{align}
\begin{split}
    &g_{lm}(r,\boldsymbol\mu_d)=\int_0^{2\pi}\int_{0}^\pi h(\vr-\boldsymbol\mu_d)\cdot Y_{lm}(\theta_\vr,\varphi_\vr)\ d\theta_\vr d\varphi_\vr\\
    &=\frac{1}{(2\pi)^{\frac{3}{2}}\sigma^3}\exp\left(-\frac{1}{2}\frac{r^2+\|\boldsymbol\mu_d\|_2^2}{\sigma^2}\right)\int_0^{2\pi}\int_0^{\pi}\exp\left(\frac{\vr^T\boldsymbol\mu_d}{\sigma^2}\right)\cdot Y_{lm}(\theta_\vr,\varphi_\vr)\sin\theta_\vr d\theta_\vr d\varphi_\vr\,.
\end{split}
\end{align}

\begin{enumerate}[label={\arabic*)}]
\item When $\boldsymbol\mu_d=\boldsymbol 0$, we have
\begin{align}
\begin{split}
    h(\vr-\boldsymbol 0)&=\left(\frac{1}{(2\pi)^{\frac{3}{2}}\sigma^3}\exp\left(-\frac{1}{2}\frac{r^2}{\sigma^2}\right)\sqrt{4\pi}\right)\cdot\sqrt{\frac{1}{4\pi}}\\
    &=g_{00}(r, \boldsymbol 0)\cdot Y_{00}(\theta_{\boldsymbol\mu_d},\varphi_{\boldsymbol\mu_d})\,,
\end{split}
\end{align}
where $Y_{00}(\theta_{\boldsymbol\mu_d},\varphi_{\boldsymbol\mu_d})=\sqrt{\frac{1}{4\pi}}$.

\item When $\boldsymbol\mu_d\neq \boldsymbol 0$, we can use the Funk-Hecke formula and get
\begin{align}
\label{eq:funk_hecke_1}
\begin{split}
    &\int_0^{2\pi}\int_0^{\pi}\exp\left(\frac{\vr^T\boldsymbol\mu_d}{\sigma^2}\right)\cdot Y_{lm}(\theta_\vr,\varphi_\vr)\sin\theta_\vr\ d\theta_\vr d\varphi_\vr\\
    &= \frac{4\pi\cdot c_l(r,\|\boldsymbol\mu_d\|_2)}{2l+1}\cdot Y_{lm}(\theta_{\boldsymbol\mu_d},\varphi_{\boldsymbol\mu_d})\,,
\end{split}
\end{align}
where $c_l(r,\|\boldsymbol\mu_d\|_2)$ is the Legendre series expansion coefficient of $\exp\left(\frac{\vr^T\boldsymbol\mu_d}{\sigma^2}\right)$. Let $\psi_{\vr,\boldsymbol\mu_d}$ denote the angle between $\vr$ and $\boldsymbol\mu_d$, and $\kappa=\cos\psi_{\vr,\boldsymbol\mu_d}$. We can compute $c_l(r,\|\boldsymbol\mu_d\|_2)$ as follows
\begin{align}
\label{eq:funk_hecke_2}
\begin{split}
    c_l(r,\|\boldsymbol\mu_d\|_2)&=\frac{2l+1}{2}\int_{-1}^{1}\exp\left(\frac{r\|\boldsymbol\mu_d\|_2\cdot \kappa}{\sigma^2}\right)\cdot P_l(\kappa) d\kappa\\
    &=\sqrt{\frac{\pi}{2r\|\boldsymbol\mu_d\|_2}}\cdot\sigma\cdot(2l+1)\cdot I_{l+\frac{1}{2}}\left(\frac{r\|\boldsymbol\mu_d\|_2}{\sigma^2}\right)\,,
\end{split}
\end{align}
where $P_l(\cdot)$ is the Legendre polynomial, and $I_{\nu}(\cdot)$ is the modified Bessel function of the first kind of order $\nu$. To avoid numeric overflow, we often use the scaled Bessel function $\widehat{I}_\nu(\cdot)$ instead.
\begin{align}
    \label{eq:scaled_bessel}
    \widehat{I}_{l+\frac{1}{2}}\left(\frac{r\|\boldsymbol\mu_d\|_2}{\sigma^2}\right) = \exp\left(-\frac{r\|\boldsymbol\mu_d\|_2}{\sigma^2}\right)\cdot I_{l+\frac{1}{2}}\left(\frac{r\|\boldsymbol\mu_d\|_2}{\sigma^2}\right)\,.
\end{align}

The expansion coefficient $g_{lm}(r,\boldsymbol\mu_d)$ is then
\begin{align}
\begin{split}
    g_{lm}(r,\boldsymbol\mu_d)=&\frac{1}{\sigma^2}\exp\left(-\frac{1}{2}\frac{(r-\|\boldsymbol\mu_d\|_2)^2}{\sigma^2}\right) \sqrt{\frac{1}{r\|\boldsymbol\mu_d\|_2}}\cdot \widehat{I}_{l+\frac{1}{2}}\left(\frac{r\|\boldsymbol\mu_d\|_2}{\sigma^2}\right)\\
    &\times Y_{lm}(\theta_{\boldsymbol\mu_d},\varphi_{\boldsymbol\mu_d})\,.
\end{split}
\end{align}
\end{enumerate}

\section{Linear Formulation of Spatial Radial Feature}
\label{app:sec:radial_integration_feature}
We first give the detailed derivation of \eqref{eq:radial_integration_feature_step_0}
\begin{align}
\begin{split}
M(k)&=\frac{1}{N}\sum_{n=1}^N\frac{1}{2\pi}\int_0^{2\pi}\widehat{S}_n(\vk)\ d\varphi\\
&\overset{N\rightarrow\infty}{=}\frac{\int_0^{2\pi}\int_0^{\pi}\widehat{\rho}(\vk)\cdot k^2\sin\theta_{\vk}d\varphi_{\vk}d\theta_{\vk}}{\int_0^{2\pi}\int_0^{\pi}k^2\sin\theta_{\vk}\ d\varphi_{\vk}d\theta_{\vk}}\\
&=\iiint\frac{\int_0^{2\pi}\int_0^{\pi}e^{-i\vk\cdot\vr}k^2\sin\theta_{\vk}\ d\theta_{\vk}d\varphi_{\vk}}{\int_0^{2\pi}\int_0^{\pi}k^2\sin\theta_{\vk}\ d\varphi_{\vk}d\theta_{\vk}}\cdot\rho(\vr)\ d\vr\\
&=\iiint\frac{\int_0^{2\pi}d\varphi_{\vk,\vr}\int_0^{\pi}e^{-ikr\cos\theta_{\vk,\vr}}k^2\sin\theta_{\vk,\vr}\ d\theta_{\vk,\vr}}{\int_0^{2\pi}d\varphi_{\vk,\vr}\int_0^{\pi}k^2\sin\theta_{\vk,\vr}\ d\theta_{\vk,\vr}}\cdot\rho(\vr)\ d\vr\\
&=\iiint\frac{\sin(kr)}{kr}\cdot\rho(\vr)\ d\vr\,,
\end{split}
\tag{\ref{eq:radial_integration_feature_step_0} revisited}
\end{align}
where $\varphi$ is the azimuth angle of $\vk$ in the Fourier slice $\widehat{S}_n$, and $(\theta_{\vk,\vr}, \varphi_{\vk,\vr})$ is the angular direction of $\vk$ with respect to $\vr$ when $\vr$ was selected as the pseudo $z$-axis during a change of coordinates.

As derived in \eqref{eq:radial_integration_linear}, we can compute the radial feature $W(r)$ as follows
\begin{align}
\label{eq:radial_integration_feature_step_1}
\begin{split}
    W(r)=&\iiint \rho(\tilde\vr)\delta(r-\tilde{r})\ d\tilde\vr\\
    =&\sum_{d=1}^D w_d\cdot\iiint \frac{1}{(2\pi)^{\frac{3}{2}}\sigma^3}\exp\left(-\frac{1}{2}\frac{\|\tilde\vr-\boldsymbol \mu_d\|_2^2}{\sigma^2}\right)\cdot\delta(r-\tilde{r})\ d\tilde\vr\\
    =&\sum_{d=1}^D w_d\cdot\frac{1}{(2\pi)^{\frac{3}{2}}\sigma^3}\int\exp\left(-\frac{1}{2}\frac{\tilde r^2+\|\boldsymbol\mu_d\|_2^2}{\sigma^2}\right)\\
    &\quad\quad\times\left(\int_0^{2\pi}\int_0^\pi\exp\left(\frac{\tilde\vr^T\boldsymbol\mu_d}{\sigma^2}\right)\sin\theta_{\tilde\vr}\ d\theta_{\tilde\vr} d\varphi_{\tilde\vr}\right)\cdot \tilde r^2\cdot\delta(r-\tilde{r})\ d\tilde r\,,
\end{split}
\end{align}
which can be simplified as
\begin{align}
     W(r)=\sum_{d=1}^Dw_d\cdot g(r,\boldsymbol\mu_d)\,.
\end{align}
\begin{enumerate}[label={\arabic*)}]
\item When $r=0$, we have
\begin{align}
    g(r,\boldsymbol\mu_d) = \frac{1}{(2\pi)^{\frac{3}{2}}\sigma^3}\exp\left(-\frac{1}{2}\frac{\|\boldsymbol\mu_d\|_2^2}{\sigma^2}\right)\,.
\end{align}
\item When $r\neq 0$ and $\|\boldsymbol\mu_d\|_2=0$, we have
\begin{align}
    g(r,\boldsymbol\mu_d) = 4\pi r^2\cdot\frac{1}{(2\pi)^{\frac{3}{2}}\sigma^3}\exp\left(-\frac{1}{2}\frac{r^2}{\sigma^2}\right)\,.
\end{align}
\item When $r\neq 0$ and $\|\boldsymbol\mu_d\|_2\neq 0$, by using \eqref{eq:funk_hecke_1},\eqref{eq:funk_hecke_2} and setting $l=0,m=0$, we can compute
\begin{align}
\label{eq:funk_hecke_3}
    \int_0^{2\pi}\int_0^\pi\exp\left(\frac{\tilde\vr^T\boldsymbol\mu_d}{\sigma^2}\right)\sin\theta_{\tilde\vr}\ d\theta_{\tilde\vr} d\varphi_{\tilde\vr} = 4\pi\cdot\sqrt{\frac{\pi}{2r\|\boldsymbol\mu_d\|_2}}\cdot\sigma\cdot I_{\frac{1}{2}}\left(\frac{r\|\boldsymbol\mu_d\|_2}{\sigma^2}\right)\,.
\end{align}
Plugging \eqref{eq:funk_hecke_3} into \eqref{eq:radial_integration_feature_step_1}, and using the scaled Bessel function in \eqref{eq:scaled_bessel}, we have
\begin{align}
    g(r,\boldsymbol\mu_d)=\frac{r^2}{\sigma^2}\exp\left(-\frac{1}{2}\frac{(r-\|\boldsymbol\mu_d\|_2)^2}{\sigma^2}\right)\sqrt{\frac{1}{r\|\boldsymbol\mu_d\|_2}}\cdot\widehat{I}_{\frac{1}{2}}\left(\frac{r\|\boldsymbol\mu_d\|_2}{\sigma^2}\right)\,.
\end{align}
\end{enumerate}
We then have the following measurement vector $\vg(r)$:
\begin{align}
\vg(r) = [g(r,\boldsymbol\mu_1)\quad g(r,\boldsymbol\mu_2)\quad\cdots\quad g(r,\boldsymbol\mu_D)]^T\,.
\end{align}

\section{Quadratic Formulation of Spatial Autocorrelation Feature}
\label{app:sec:quad_autocorrelation}
As derived in \eqref{eq:autocorrelation_1} of Section \ref{subsec:autocorrelation_features}, the autocorrelation feature $E_l(r_1,r_2)$ is given by
\begin{align}
\begin{split}
    E_l(r_1,r_2)&=\sum_{m=-l}^l\iiint\rho(\tilde\vr_1)\cdot\delta(r_1-\tilde{r}_1)\cdot Y_{lm}(\theta_{\tilde\vr_1},\varphi_{\tilde\vr_1})\sin\theta_{\tilde\vr_1}\ d\tilde r_1 d\theta_{\tilde\vr_1} d\varphi_{\tilde\vr_1}\\
    &\quad\quad\quad\quad\quad\times\iiint\rho(\tilde\vr_2)\cdot\delta(r_2-\tilde{r}_2)\cdot Y_{lm}(\theta_{\tilde\vr_2},\varphi_{\tilde\vr_2})\sin\theta_{\tilde \vr_2}\ d\tilde r_2d\theta_{\tilde\vr_2}\varphi_{\tilde\vr_2}\,. 
\end{split}
\end{align}

Using \eqref{eq:3d_density_sh_expansion}, we can compute the following integration
\begin{align}
\label{eq:rho_integration}
\begin{split}
    &B_l^m(r)=\iiint\rho(\tilde\vr)\cdot\delta(r-\tilde{r})\cdot Y_{lm}(\theta_{\tilde\vr},\varphi_{\tilde\vr})\sin\theta_{\tilde\vr}\ d\tilde r d\theta_{\tilde\vr} d\varphi_{\tilde\vr} \\
    &= \sum_{d=1}^Dw_d\cdot\sum_{l^\prime=0}^\infty\sum_{m^\prime=-l^\prime}^{l^\prime}g_{l^\prime m^\prime}(r,\boldsymbol\mu_d)\cdot \iint Y_{l^\prime m^\prime}(\theta_{\tilde\vr},\varphi_{\tilde\vr})\cdot Y_{lm}(\theta_{\tilde\vr},\varphi_{\tilde\vr})\sin\theta_{\tilde\vr}\ d\theta_{\tilde\vr} d\varphi_{\tilde\vr}\\
    &=\sum_{d=1}^Dw_d\cdot\sum_{l^\prime=0}^\infty\sum_{m^\prime=-l^\prime}^{l^\prime}g_{l^\prime m^\prime}(r,\boldsymbol\mu_d)\cdot\delta(l^\prime=l)\delta(m^\prime=m)\\
    &=\sum_{d=1}^Dw_d\cdot g_{lm}(r,\boldsymbol\mu_d)\\
    &=\vw^T\vg_{lm}(r)\,.
\end{split}
\end{align}

Using \eqref{eq:rho_integration}, we can get the quadratic formulation of $E_l(r_1,r_2)$:
\begin{align}
\begin{split}
    E_l(r_1,r_2)&=\sum_{m=-l}^l\vw^T\vg_{lm}(r_1)\cdot\vg_{lm}^T(r_2)\vw\\
    &=\vw^T\cdot\left(\sum_{m=-l}^l\vg_{lm}(r_1)\cdot\vg_{lm}^T(r_2)\right)\cdot\vw\,,
\end{split}  \tag{\ref{eq:autocorrelation_quad} revisited}
\end{align}
where $\vw=[w_1,\cdots,w_D]^T$, and the real measurement vector $\vg_{lm}(r)$ is 
\begin{align}
    \vg_{lm}(r)=\left[g_{lm}(r,\boldsymbol\mu_1)\quad g_{lm}(r,\boldsymbol\mu_2)\quad\cdots\quad g_{lm}(r,\boldsymbol\mu_D)\right]^T\,.
\end{align}

\section{Computational Complexity}
\label{app:sec:complexity}
Let $\Phi$ denote the number of selected GLQ points along the dimension of $\varphi$, $V$ denote the number of selected GLQ points along the dimension of $r$, $U$ denote the number of sampling points along the dimension of $k$, and $L$ denote the spherical harmonic bandwidth. There are a total of $N$ projection images, and the size of the projection image is $G\times G$. In general, the number of sampled points $\Phi,V,U$ scales linearly with $G$, i.e., $\mathcal{O}(G)$. For simplification, we shall give the complexity in terms of $L,N$ and $G$. The proposed OMR-SC approach consists of the following three steps with different computational complexities:
\begin{enumerate}[label={\arabic*)}]
\item \emph{Feature extraction.} 
The complexity of computing the Fourier-Bessel expansion $\widetilde{C}_{\mathrm{FB}}$ of $\widetilde{C}(k_1,k_2,\psi)$ using the fast steerable PCA and CWF is $\mathcal{O}(NG^3)$. The complexity of calculating $\{C_l(k_1,k_2)\}_l$ from $\widetilde{C}_{\mathrm{FB}}$ is $\mathcal{O}(LG^3)$. 
The complexity of performing Cholesky decompositions of $\{\mC_l\}_l$ is $\mathcal{O}(LU^3)=\mathcal{O}(LG^3)$. The overall complexity of extracting the Fourier autocorrelation features is $\mathcal{O}(NG^3+LG^3+LG^3)=\mathcal{O}(NG^3+LG^3)$.

The complexity of calculating the non-uniform FFT $\{\widehat{S}_n(k,\varphi)\}_n$ from $N$ projection images is $\mathcal{O}\left(NG^2\log G\right)$.  
The complexity of calculating $\{M(k)\}_k$ is $\mathcal{O}(NU\Phi)=\mathcal{O}(NG^2)$. The complexity of calculating $\{W(r)\}_r$ is $\mathcal{O}(VU)=\mathcal{O}(G^2)$. The overall complexity of extracting the spatial radial features is $\mathcal{O}\left(NG^2\log G+NG^2+G^2\right)=\mathcal{O}\left(NG^2\log G\right)$.

The overall complexity of feature extraction is as $\mathcal{O}\left(NG^3+LG^3+NG^2\log G\right)=\mathcal{O}\left(NG^3+LG^3\right)$.

\item \emph{Optimization of $\{\mO_l\}_{l=1}^L$.} The complexity of computing $\{\mB_l^T(\vw)\mQ_l\mF_l\}_l$ is $\mathcal{O}(L^2(LVL+VUL))=\mathcal{O}(L^3(LV+VU))=\mathcal{O}(L^3(LG+G^2))$. The complexity of performing SVD of $\{\mB_l(z)^*\mQ_l\mF_l\}_l$ is $\mathcal{O}(L^4)$. The complexity of computing $\{\mV_l\mU_l^T\}_l$ is $\mathcal{O}(L^4)$. The overall complexity of computing $\{\mO_l\}_{l=1}^L$ is thus $\mathcal{O}(L^3(LG+G^2))$.

\item \emph{Optimization of $\vw$.} The complexity of computing the vector $\vb_l(r,\vw)$ is $\mathcal{O}\big(lr^2\big)$. Since $V$ and $G$ are typically in the same order as $G$, i.e. $\mathcal{O}(V)=\mathcal{O}(U)=\mathcal{O}(G)$, the complexity of computing $\{\mB_l(\vw)\}_l$ is then $\mathcal{O}(L^2G^3)$. Computing $\{\mA_l\}_l$ via \eqref{eq:evaluate_A_matrix_form} only increases the complexity of evaluating features to $\mathcal{O}(L^2G^3+L^2G^2)$. The complexity of evaluating a radial feature $\vg(v)^T\vw$ is $\mathcal{O}(r^2)$, and the complexity of evaluating all radial features is $\mathcal{O}(G^3)$. The complexity of computing the gradient $\nabla f_2(\vw)$ is $\mathcal{O}(L^2G^3)$, and the projection operation $\mathcal{P}_\mathcal{S}(\cdot)$ has a complexity of $\mathcal{O}(G^3\log G)$. The overall complexity of computing $\vw$ is $\mathcal{O}(L^2G^3+G^3\log G)$.
\end{enumerate}

\section{OMR-SC Using Spatial Autocorrelations}
\label{app:sec:omr_sc_s}

For the spatial autocorrelation feature matrix $\mE_l$, we have for any orthogonal $\mO_l$ of size $(2l+1)\times(2l+1)$ that
\begin{align}
    \mE_l=\mB_l\mO_l^T\mO_l\mB_l^T=\mP_l\mP_l^T\,,
\end{align}
where $\mP_l=\mB_l\mO_l^T$ is obtained via the Cholesky decomposition of $\mE_l$. The formulation of OMR-SC using spatial autocorrelations (OMR-SC-S) is then
\begin{equation}
\tag{OMR-SC-S}
    \begin{aligned}
    \label{eq:spatial_objective_function}
    \minimize_{\vw,\{\mO_l\}_{l=0}^{L}}&\quad f(\vw,\mO_l):=\sum_{l=0}^L\left\|\mP_l\mO_l-\mB_l(\vw)\right\|_2^2\\
    &\quad\quad\quad\quad\quad\quad\quad+\lambda\cdot\sum_{v=1}^V\left(\vg(v)^T\vw-W(v)\right)^2\\
    &\quad\quad\quad\quad\quad\quad\quad+\xi\cdot\sum_{x,y}\big(P_{\vw}(x,y)-\overline{S}(x,y)\big)^2,\\
    \textnormal{subject to}&\quad 0\leq w_d\leq W_\rho,\\
    &\quad \sum_{d=1}^Dw_d=W_\rho\,,\\
    &\quad \mO_l^T \mO_l =\mO_l\mO_l^T =\mI, \quad l \in \{0, \ldots, L\}\,.
    \end{aligned}
\end{equation}
As summarized in Algorithm \ref{alg:OMR-SC-S}, we can compute $\{\mO_l\}_{l=1}^L$ and $\vw$ in a similar alternating fashion as in Section \ref{subsec:OMR-SC-F}.

\begin{algorithm}[tbp]
\caption{OMR with Spatial Consensus Using Spatial Autocorrelations}
\label{alg:OMR-SC-S}
\begin{algorithmic}[1]
\REQUIRE Denoised reference projections $\{\overline{S}_i\ |\ h=1,\cdots,H\}$, step size $\eta$, convergence threshold $\varsigma$.
\STATE Extract the spatial radial features $\{W(v)\}_v$ and Spatial autocorrelation features $\{\mE_l\}_l$.
\STATE Perform Cholesky decompositions of spatial autocorrelation matrices $\{\mE_l\}_l$.
\FOR{$h=\{1,\cdots,H\}$}
\STATE Compute the initialization $\vw_0(i)$ from the spatial radial features $\{W(v)\}_v$ and the $i$-th reference projection $\overline{S}_i$.
\FOR{$t=\{0,1,\cdots,T\}$}
    \STATE Fix $\vw_t(i)$, and update $\{\mO_l(i)\}_l$ with respect to $\vw_t(i)$ via singular value decomposition.
    \STATE Fix $\{\mO_l(i)\}_l$, and estimate $\vw_{t+1}(i)$ with respect to $\{\mO_l(i)\}_l$ via projected gradient descent.
    \IF {$\frac{\|\vw_{t+1}(i)-\vw_t(i)\|_2}{\|\vw_t(i)\|_2}<\varsigma$}
        \STATE Convergence is reached, set $\vw(i)=\vw_{t+1}(i)$ and \textbf{break}.
    \ENDIF
\ENDFOR
\STATE Save the $i$-th set of solutions $\{\vw(i), \{\mO_l(i)\}_l\}$.
\ENDFOR
\STATE Find the set of solutions that minimizes the MSE of autocorrelation features:
\begin{align}
    \tilde{i} = \arg\min_i\quad\sum_{l=0}^L\left\|\mP_l\mO_l(i)-\mB_l(\vw(i))\right\|_2^2
\end{align}
\STATE {\bfseries Return} $\tilde{\vw}=\vw(\tilde{i})$.
\end{algorithmic}
\end{algorithm}

\subsection{Performance Comparison of Spatial and Fourier Autocorrelations}
\label{app:sec:compare_spatial_fourier_autocorrelations}

\begin{figure*}[htbp]
\centering
\subfloat[Spatial expansion coefficients]{
\label{fig:spatial_feature_energy}
\includegraphics[width=0.45\textwidth]{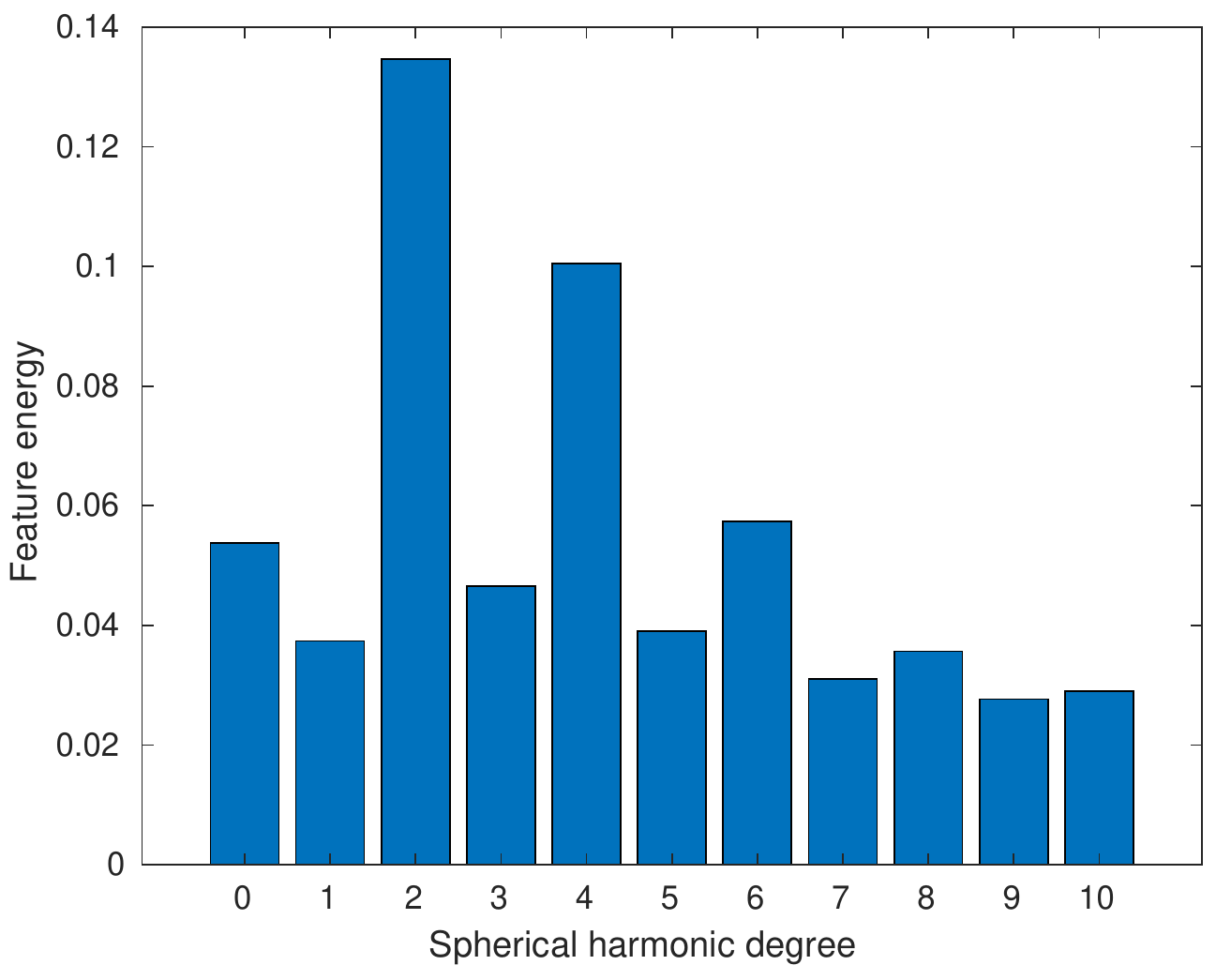}}
\subfloat[Fourier expansion coefficients]{
\label{fig:fourier_feature_energy}
\includegraphics[width=0.45\textwidth]{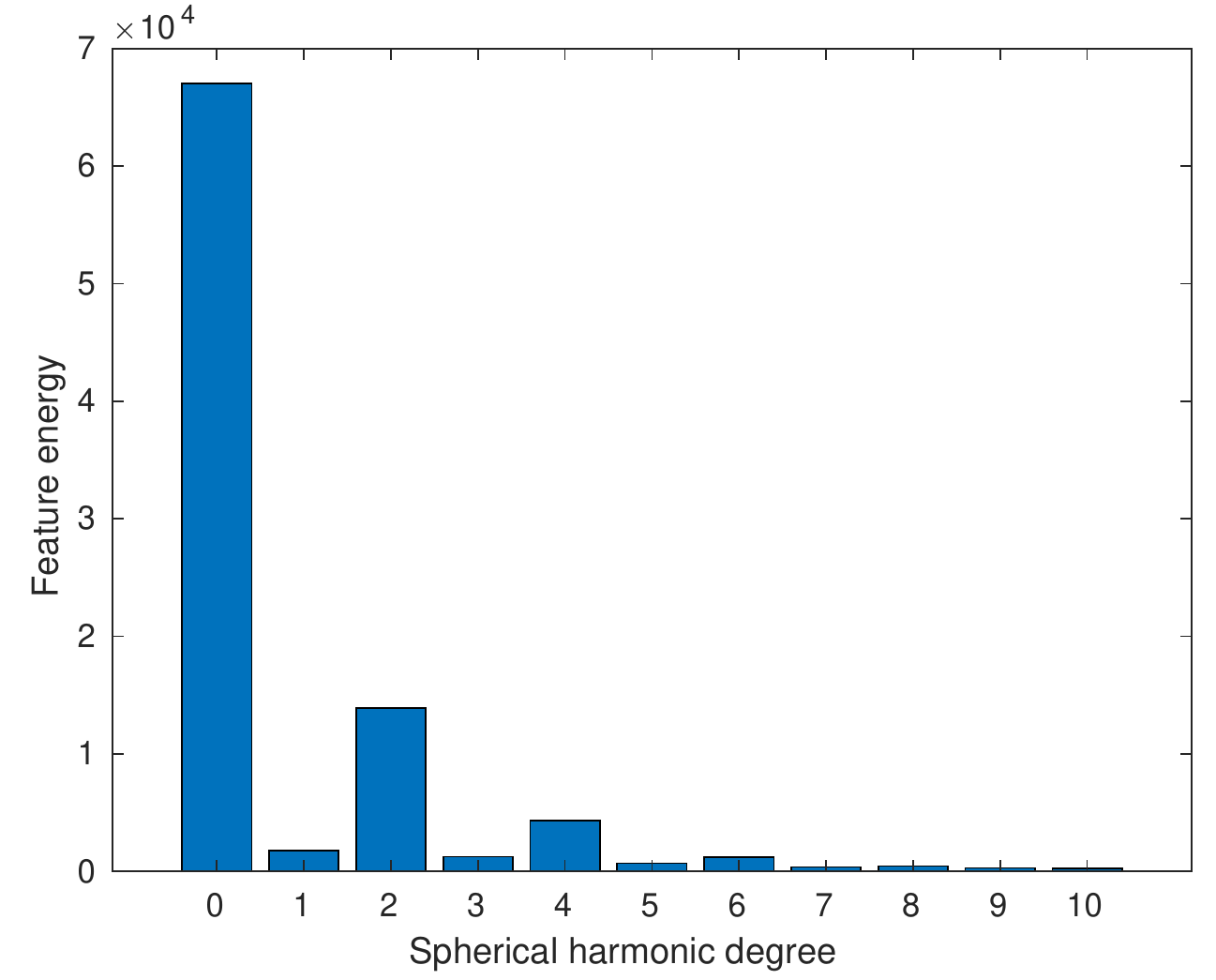}}
\caption{Energy distributions of spatial and Fourier spherical harmonic expansion coefficients.}
\label{fig:spatial_fourier_feature_energy}
\end{figure*}

\begin{figure*}[htbp]
\centering
\subfloat[OMR-SC-S]{
\label{fig:spatial_gradient_norm}
\includegraphics[width=0.45\textwidth]{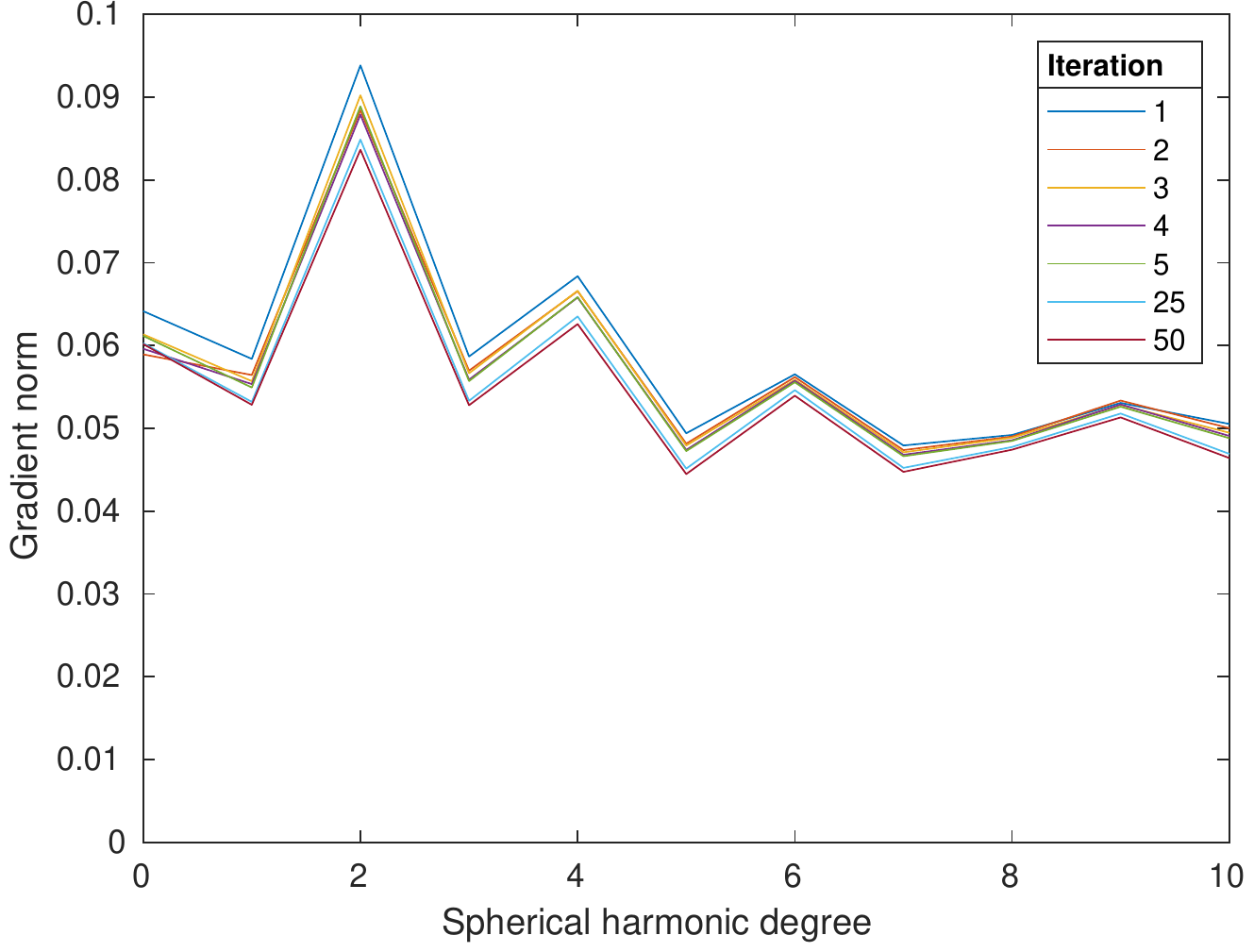}}
\subfloat[OMR-SC-F]{
\label{fig:fourier_gradient_norm}
\includegraphics[width=0.45\textwidth]{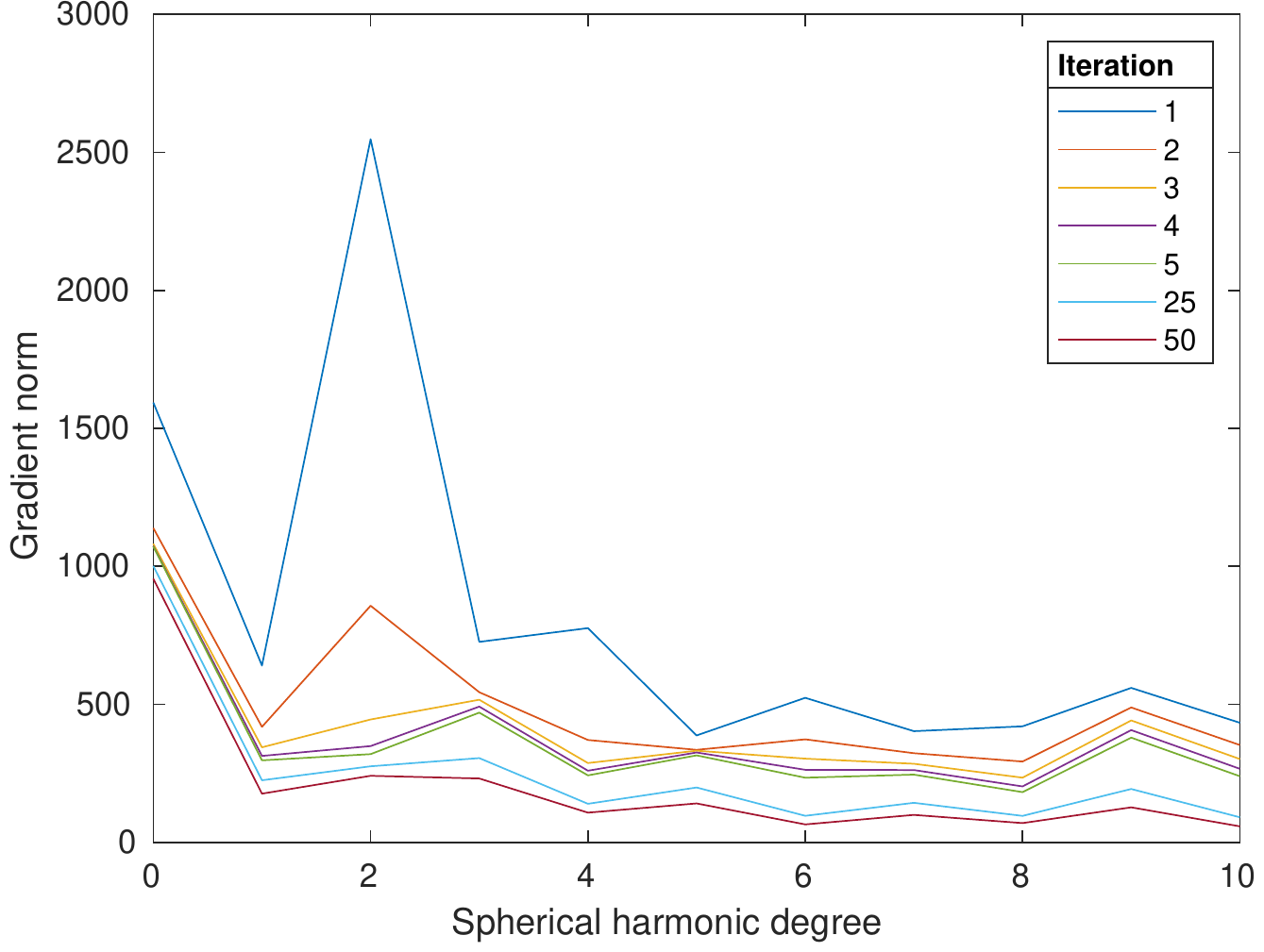}}
\caption{Gradient norms of OMR-SC-S and OMR-SC-F.}
\label{fig:spatial_fourier_gradient_norm}
\end{figure*}

As shown in Tables \ref{tab:compare_resolution}-\ref{tab:compare_cc}, OMR-SC with Fourier autocorrelations (OMR-SC-F) performs better than OMR-SC with spatial autocorrelations (OMR-SC-S). To find out possible causes for the performance differences, we take a closer look at the two features. Take the density ``D1'' for example, as shown in Fig. \ref{fig:spatial_fourier_feature_energy}, the energy of spatial spherical harmonic expansion coefficients is more evenly distributed across different degrees than that of Fourier spherical harmonic expansion coefficients, and the energy of Fourier expansion coefficients is mostly concentrated in the first few degrees. As a result, OMR-SC puts more effort into minimizing the errors from the lower spherical harmonic degrees when Fourier autocorrelations are used. This can be verified by comparing the gradient norms of OMR-SC-S and OMR-SC-F. As shown in Fig. \ref{fig:spatial_fourier_gradient_norm}, the gradient norms of OMR-SC-S are more evenly distributed, and the reductions of the gradient norms through the iterations are generally the same across different degrees. Whereas the gradient norms of OMR-SC-F are much larger in the first three degrees, and the reductions of the gradient norms through the iterations are more significant in the first three degrees, which leads to much stronger spherical-harmonic frequency marching effect. OMR-SC-F thus focuses more on recovering the low-resolution base structure of the density determined by the lower-degree spherical harmonic expansion coefficients, which proves to be beneficial during reconstruction.

\section{Ab initio Models from OMR-SC-F}
\label{app:sec:ab_initio_omr_sc_f}

We calculate the resolutions and correlation coefficients of the \emph{ab initio} models from OMR-SC-F by downsampling the projection images. As discussed in Section \ref{subsec:reconstruction_3d_density_maps}, a $33\times 33 \times 33$ \emph{ab initio} model is computed for every density map. The refinements are generally significantly better than the \emph{ab initio} models. In particular, due to the noisy features, the refinement of PTCH1 in the noisy case turns out to be worse than the \emph{ab initio} model.

\begin{table}[htb]
\caption{Resolutions (in voxel) of the \emph{ab initio} models and refinements produced by OMR-SC-F (FSC cutoff threshold = 0.5).}
\label{tab:compare_resolution_initialization}
\centering
\resizebox{\textwidth}{!}{
\begin{tabular}{llcccccccccc}
\toprule
&OMR-SC-F &D1 &D2 &D3 &D4 &D5 &D6 &D7 &D8 &D9 &D10\\ \midrule
&Ab initio & 20.88 & 30.77 & 11.83 & 12.21 & 15.29 & 13.70 & 32.36 & 13.97 & 14.20 & 13.76   \\
\multirow{-2}{*}{Noiseless} &Refinement & \bf{6.70} & \bf{9.62} & \bf{5.90} & \bf{7.28} & \bf{5.94} & \bf{4.10} & \bf{10.63} & \bf{7.48} & \bf{8.46} & \bf{5.54} \\ \midrule
&Ab initio & 12.76 & 22.88 & 11.35 & 29.67 & 24.39 & 10.91 & 42.19 & 17.57 & 16.86 & 11.51  \\
\multirow{-2}{*}{Noisy} &Refinement & \bf{8.47} & \bf{10.80} & \bf{7.42} & \bf{10.98} & \bf{9.18} & \bf{6.01} & \bf{17.51} & \bf{8.94} & \bf{12.14} & \bf{8.45} \\ \bottomrule
\end{tabular}
}
\end{table}

\begin{table}[htb]
\caption{Correlation coefficients of the \emph{ab initio} models and refinements produced by OMR-SC-F.}
\label{tab:compare_cc_initialization}
\centering
\resizebox{\textwidth}{!}{
\begin{tabular}{llcccccccccc}
\toprule
&OMR-SC-F &D1 &D2 &D3 &D4 &D5 &D6 &D7 &D8 &D9 &D10\\ \midrule
&Ab initio & 0.58 & 0.43 & 0.80 & 0.79 & 0.67 & 0.64 & 0.47 & 0.76 & 0.71 & 0.71   \\
\multirow{-2}{*}{Noiseless} &Refinement & \bf{0.92} & \bf{0.84} & \bf{0.96} & \bf{0.92} & \bf{0.96} & \bf{0.96} & \bf{0.83} & \bf{0.91} & \bf{0.91} & \bf{0.91}  \\ \midrule
&Ab initio & 0.74 & 0.52 & 0.80 & 0.40 & 0.49 & 0.72 & 0.46 & 0.69 & 0.75 & 0.76    \\
\multirow{-2}{*}{Noisy} &Refinement & \bf{0.85} & \bf{0.71} & \bf{0.92} & \bf{0.82} & \bf{0.88} & \bf{0.89} & \bf{0.71} & \bf{0.85} & \bf{0.77} & \bf{0.86} \\ \bottomrule
\end{tabular}
}
\end{table}

\begin{table}[htb]
\caption{Resolutions (\r{A}) and correlation coefficients of the \emph{ab initio} models and refinements produced by OMR-SC-F (FSC cutoff threshold=0.5). }
\label{tab:compare_res_cc_emd_25143_hjc_ptch1_initialization}
\centering
\resizebox{\textwidth}{!}{
\begin{tabular}{llcccccccccc}
\toprule
& &\multicolumn{3}{c}{Resolution(\r{A})} &\multicolumn{3}{c}{Correlation Coefficient}\\ \cmidrule(lr){3-5} \cmidrule(lr){6-8}
&OMR-SC-F & CaS & HJC & PTCH1 & CaS & HJC & PTCH1 \\ \midrule
&Ab initio & 25.77 &45.21 & 27.08 &0.87 &0.69 &0.80 \\
\multirow{-2}{*}{Noiseless} &Refinement & \bf{14.04} & \bf{29.38} & \bf{16.52} &\bf{0.93} &\bf{0.78} &\bf{0.86} \\ \midrule
&Ab initio & 50.64 &86.21 & \bf{30.67} & 0.74 &0.62 &\bf{0.79}  \\
\multirow{-2}{*}{Noisy} &Refinement & \bf{16.15} & \bf{36.82} &31.38 &\bf{0.87} & \bf{0.73} &0.73 \\ \bottomrule
\end{tabular}
}
\end{table}

\end{appendices}

\clearpage

\bibliographystyle{IEEEbib}
\bibliography{references}

\clearpage

\setcounter{section}{0}
\begin{center}
{\LARGE Supplementary Material}
\end{center}

\vspace{5em}

The Supplementary Material contains figures showing additional FSC curves and reconstructed density maps.


\renewcommand{\arraystretch}{1.2}
\begin{table}[htp]
\caption{Table of Notations.}
\label{tab:notations}
\centering
\begin{tabularx}{\textwidth}{lXl}
\toprule
Symbol 	&Description &Location \\ \midrule
$x$ & Spatial domain: the Cartesian coordinate in the $x$-direction & Sec. 2\\
$y$ & Spatial domain: the Cartesian coordinate in the $y$-direction & Sec. 2\\
$z$ & Spatial domain: the Cartesian coordinate in the $z$-direction & Sec. 2\\
$\vr$ & Spatial domain: the vector containing the Cartesian coordinates & Sec. 2  \\ 
$r$ & The $l_2$ norm of $\vr$: $r=\|\vr\|_2$ & Sec. 2 \\
$\theta_{\vr}$ & The polar angle of $\vr$ in the spherical coordinates & Sec. 3\\
$\varphi_{\vr}$ & The azimuthal angle of $\vr$ in the spherical coordinates & Sec. 3\\
$\rho(\vr)$ &The 3D density map in the spatial domain & Sec. 2 \\ 
$\mathfrak{i}$ & The imaginary unit & Sec. 2\\
$\widetilde{\rho}(\vr)$ &The recovered 3D density map in the spatial domain & Sec. 2\\
$N$ & The number of projection images & Sec. 2\\
$n$ & The index of the projection image & Sec. 2\\
$\mR_n$ & The $3\times 3$ Rotation matrix & Sec. 2\\
$\chi_n$ & A rotatation in the 3D rotation group $\mathrm{SO}(3)$ & Sec. 2\\
$P_n(x,y)$ & The noiseless projection image & Sec. 2\\
$\epsilon_n(x,y)$ & The additive noise added to the projection image & Sec. 2\\
$S_n(x,y)$ & The observed noisy projection image & Sec. 2\\
$k_x$ & Frequency domain: the Cartesian coordinate in the $x$-direction & Sec. 2\\
$k_y$ & Frequency domain: the Cartesian coordinate in the $y$-direction & Sec. 2\\
$k_z$ & Frequency domain: the Cartesian coordinate in the $z$-direction & Sec. 2\\
$\vk$ & Frequency domain: the vector containing the Cartesian coordinates & Sec. 2  \\ 
$k$ & The $l_2$ norm of $\vk$: $k=\|\vk\|_2$ & Sec. 2\\
$\theta_{\vk}$ & The polar angle of $\vk$ in the spherical coordinates & Sec. 2\\
$\varphi_{\vk}$ & The azimuthal angle of $\vk$ in the spherical coordinates & Sec. 2\\
$\widehat{\rho}(\vk)$ & The Fourier transform of the 3D density map $\rho(\vr)$ & Sec. 2\\
$\alpha$ & The angle between $\vk$ and $\vr$ & Sec. 2\\ \bottomrule
\end{tabularx}

\end{table}

\setcounter{table}{12}
\begin{table}[htp]
\caption{Table of Notations (continued).}
\label{tab:notations_2}
\centering
\begin{tabularx}{\textwidth}{lXl}
\toprule
Symbol 	&Description &Location \\ \midrule
$\psi$ & The angle between two frequency vectors $\vk_1$ and $\vk_2$ & Sec. 2\\
$\widehat{\epsilon}_n(k_x,k_y)$ & The Fourier transform of the additive noise $\epsilon_n(x,y)$ & Sec. 2\\
$\widehat{S}_n(k_x,k_y)$ & The Fourier transform of the projection image $S_n(x,y)$ & Sec. 2\\
$l$ & The spherical harmonic degree & Sec. 2\\
$L$ & The cutoff threshold on the spherical harmonic degree $l$: $l=0,1,\cdots,L$ & Sec. 4\\ 
$m$ & The spherical harmonic order under the degree $l$: $-l\leq m \leq l$ & Sec. 2\\
$A_{lm}(k)$ & The spherical harmonic expansion coefficients of $\widehat{\rho}(\vk)$ in the frequency domain & Sec. 2 \\
$Y_{lm}(\theta,\varphi)$ & The \emph{real} spherical harmonic function of degree $l$ and order $m$ & Sec. 2\\
$\delta(\cdot)$ & The Dirac impulse & Sec. 2\\
$P_l(\cdot)$ & The Legendre polynomial of degree $l$ & Sec. 2\\
$C_N(k_1,k_2,\psi)$ & The autocorrelation function computed from $N$ noisy projection images that contains a bias term $\zeta(k_1,k_2,\psi)$ & Sec. 2\\
$\zeta(k_1,k_2,\psi)$ & The bias term in the empirical autocorrelation function $C_N(k_1,k_2,\psi)$ & Sec. 2\\
$\widetilde{C}(k_1,k_2,\psi)$ & The ``debiased'' autocorrelation function computed from $N$ noisy projection images & Sec. 2\\
$C(k_1,k_2,\psi)$ & The asymptotic noiseless Fourier autocorrelation function & Sec. 2\\
$C_l(k_1,k_2)$ & The Fourier autocorrelation feature at $(k_1,k_2)$ with degree $l$ & Sec. 2\\
$U$ & The number of sampling points in the domain of the frequency radius $k$ & Sec. 2\\
$\mC_l$ & The $U\times U$ Fourier autocorrelation feature matrix of rank $(2l+1)$, it contains all the Fourier autocorrelation features with degree $l$ & Sec. 2\\
$\mA_l$ & The $U\times (2l+1)$ coefficient matrix that contains all the spherical harmonic coefficients $A_{lm}(k)$ with degree $l$, each column corresponds to an order $m$ & Sec. 2\\
$\mF_l$ & The matrix returned by performing Cholesky decomposition of $\mC_l=\mF_l\mF_l^*$ & Sec. 2\\
$\mO_l$ & The orthogonal matrix that rotates $\mF_l$ back to $\mA_l$: $\mF_l\mO_l = \mA_l$ & Sec. 2\\
$M_N(k)$ & The estimated first-order moment of the Fourier transform $\widehat{\rho}(\vr)$ from $N$ projection images & Sec. 3\\
$M(k)$ & The asymptotic first-order moment of the Fourier transform $\widehat{\rho}(\vr)$ & Sec. 3\\ \bottomrule
\end{tabularx}

\end{table}

\setcounter{table}{12}
\begin{table}[htp]
\caption{Table of Notations (continued).}
\label{tab:notations_3}
\centering
\begin{tabularx}{\textwidth}{lXl}
\toprule
Symbol 	&Description &Location \\ \midrule
$W(r)$ & The spatial radial feature, i.e. the integration of the density $\rho(\vr)$ on the sphere with radius $r$ & Sec. 3\\
$W_\rho$ & The total mass of the density $\rho(\vr)$ & Sec. 3\\
$B_{lm}(r)$ & The spherical harmonic expansion coefficients of $\rho(\vr)$ in the spatial domain & Sec. 3\\
$E(r_1,r_2,\psi)$ & The spatial autocorrelation function & Sec. 3\\
$E_l(r_1,r_2)$ & The spatial autocorrelation feature at $(r_1,r_2)$ with degree $l$ & Sec. 3\\
$V$ & The number of sampling points in the domain of the spatial radius $r$ & Sec. 3\\
$\mE_l$ & The $V\times V$ spatial autocorrelation feature matrix of rank $(2l+1)$, it contains all the spatial autocorrelation features with degree $l$ & Sec. 3\\
$\mB_l$ & The $V\times (2l+1)$ coefficient matrix that contains all the spherical harmonic coefficients $B_{lm}(r)$ with degree $l$ & Sec. 3\\
$j_l(\cdot)$ & The spherical Bessel function of order $l$ & Sec. 3\\
$Y_l^m(\theta,\varphi)$ & The \emph{complex} spherical harmonic function of degree $l$ and order $m$\\
$A_l^m(k)$ & The \emph{complex} spherical harmonic expansion coefficients of $\widehat{\rho}(\vk)$ in the frequency domain & Sec. 3\\
$\{q_1,\cdots,q_V\}$ & The Gauss--Legendre quadrature weights associated with the sampling locations $r\in\{v_1,\cdots,v_V\}$ in the spatial domain & Sec. 3\\
$\mQ_l$ & The $V\times U$ matrix that transforms $\mB_l$ to $\mA_l$: $\mA_l=\mQ_l^*\mB_l$ & Sec. 3\\
$G$ & The size of the Cartesian sampling grid along one direction in the spatial domain & Sec. 4\\
$d$ & The index of the sampling location on the Cartesian grid & Sec. 4\\
$\boldsymbol\mu_d$ & The coordinate of the $d$-th sampling location & Sec. 4\\
$h(\cdot)$ & A nonnegative bump function associated with the sampling grid. In this paper, it is chosen to be an isotropic Gaussian function & Sec. 4\\
$\tau$ & The width parameter of the isotropic Gaussian bump function $g(\cdot)$ & Sec. 4\\
$w_d$ & The weight corresponding to the $d$-th sampling location & Sec. 4\\
$\overline{S}(x,y)$ & The denoised reference projection image & Sec. 4\\
$\overline{s}$ & A threshold chosen to filter out those small perturbations in the denoised image $\overline{S}(x,y)$ & Sec. 4\\
$\mathcal{M}$ & The set of sampling grid points that are consistent with the reference projection $\overline{S}(x,y)$ and have radii less than $\frac{G-1}{2}$ & Sec. 4\\ \bottomrule
\end{tabularx}

\end{table}

\setcounter{table}{12}
\begin{table}[htp]
\caption{Table of Notations (continued).}
\label{tab:notations_4}
\centering
\begin{tabularx}{\textwidth}{lXl}
\toprule
Symbol 	&Description &Location \\ \midrule
$D$ & The total number of sampling locations in the set $\mathcal{M}$ & Sec. 4\\
$\vw$ & The vector containing the weights of all the sampling locations in the set $\mathcal{M}$ & Sec. 4\\
$h(\vr-\boldsymbol\mu_d)$ & The bump function at the $d$-th sampling location & Sec. 4\\
$g_{lm}(r,\boldsymbol\mu_d)$ & The spherical harmonic expansion coefficient of $h(\vr-\boldsymbol\mu_d)$ & Sec. 4 \\
$\vg_{lm}(r)$ & A vector containing the spherical harmonic expansion coefficients of all the bump functions: $\vg_{lm}(r)=[g_{lm}(r,\boldsymbol\mu_1)\quad g_{lm}(r,\boldsymbol\mu_2) \quad \cdots \quad g_{lm}(r,\boldsymbol\mu_D)]^T$ & Sec. 4\\
$g(r,\boldsymbol\mu_d)$ & The integration of the bump function $h(\vr-\boldsymbol\mu_d)$ on the sphere with radius $r$ & Sec. 4\\
$\vg(r)$ & A vector used to compute the spatial radial feature $W(r)$: $\vg(r) = [g(r,\boldsymbol\mu_1)\quad g(r,\boldsymbol\mu_2)\quad\cdots\quad g(r,\boldsymbol\mu_D)]^T$ & Sec. 4\\
$\mI$ & The identity matrix & Sec. 4\\
$I_{\nu}(\cdot)$ & The modified Bessel function of the first kind of order $\nu$ & App. C\\
$\widehat{I}_{\nu}(\cdot)$ & The scaled Bessel function $I_{\nu}(\cdot)$ & App. C\\
$b_l(v,\vw)$ & The row vector of spherical harmonic coefficient matrix $\mB_l$ & Sec. 5\\
$\mB_l(\cdot)$ & A linear operator applied on the weight vector $\vw$ such that $\mB_l = \mB_l(\vw)$ & Sec. 5\\
$P_{\vw}(x,y)$ & The projection of the density map represented by $\vw$ along the $z$-direction & Sec. 5\\
$\lambda$ & The regularization parameter corresponding to the mean-squared-error (MSE) of radial features & Sec. 5\\
$\xi$ & The regularization parameter corresponding to the mean-squared-error (MSE) of projection features & Sec. 5\\
$\mathcal{S}$ & A convex set defined by the nenegative (box) and summation constraints & Sec. 5\\
$\mathcal{P}_{\mathcal{S}}(\cdot)$ & A linear projection of an vector $\vw$ onto the convex set $\mathcal{S}$ & Sec. 5\\
$\mV,\mU$ & The two matrix are obtained by performing the singular value decomposition of $\mB_l^T(\vw)\mQ_l\mF_l$ & Sec. 5\\
$\eta$ & The step size of the projected gradient descent & Sec. 5\\
$j_{ls}(k)$ & The normalized spherical Bessel function & Sec. 6\\
$a_{lms}$ & The spherical Bessel expansion coefficients of $A_{lm}(k)$ & Sec. 6\\
& place holder & \\ \bottomrule
\end{tabularx}

\end{table}

\clearpage

\begin{figure*}[tb]
\centering
\subfloat[Noiseless recovery]{
\label{fig:fsc_noiseless_d2}
\includegraphics[width=0.45\textwidth]{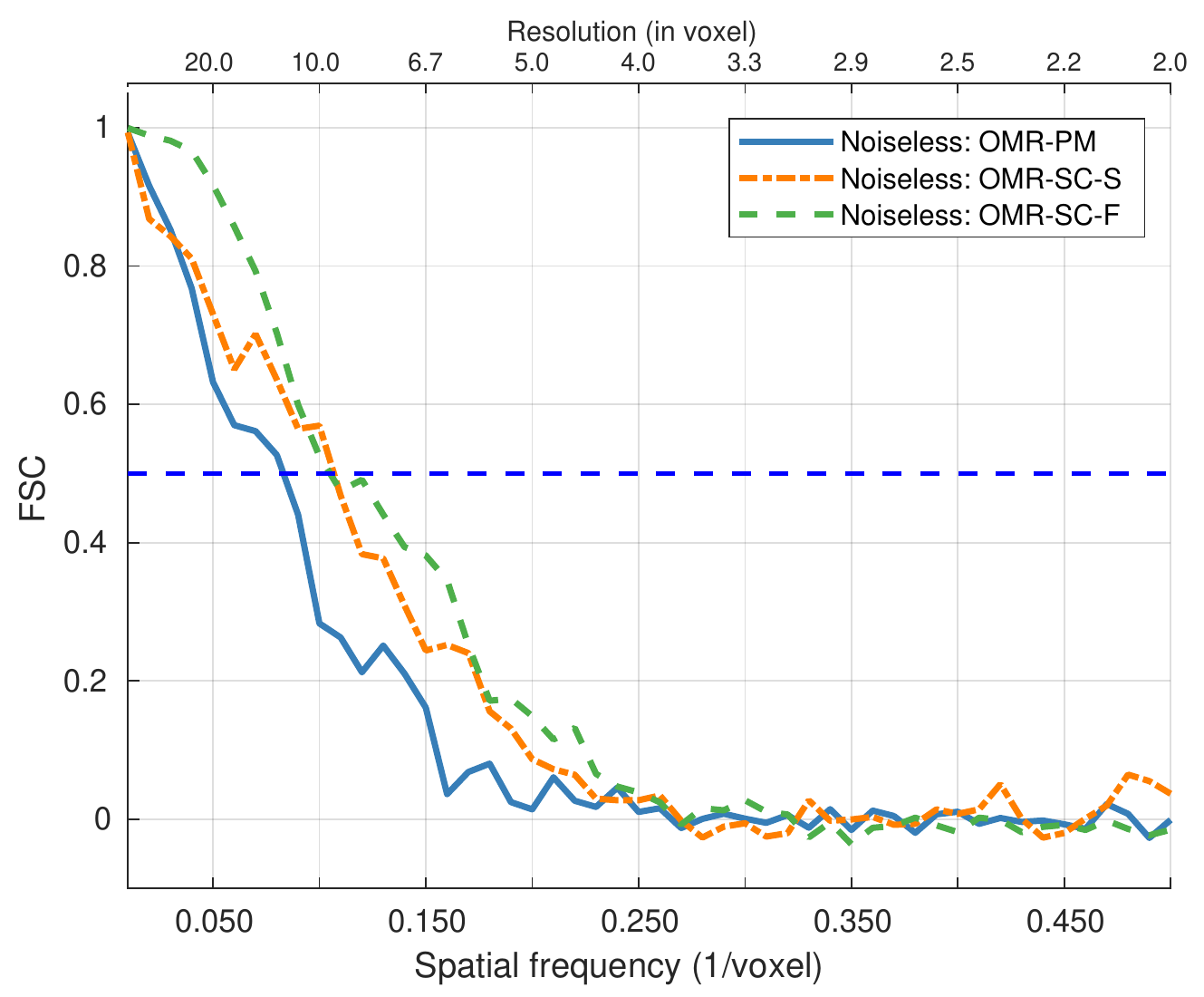}}
\hspace{0.5cm}
\subfloat[Noisy recovery]{
\label{fig:fsc_noisy_d2}
\includegraphics[width=0.45\textwidth]{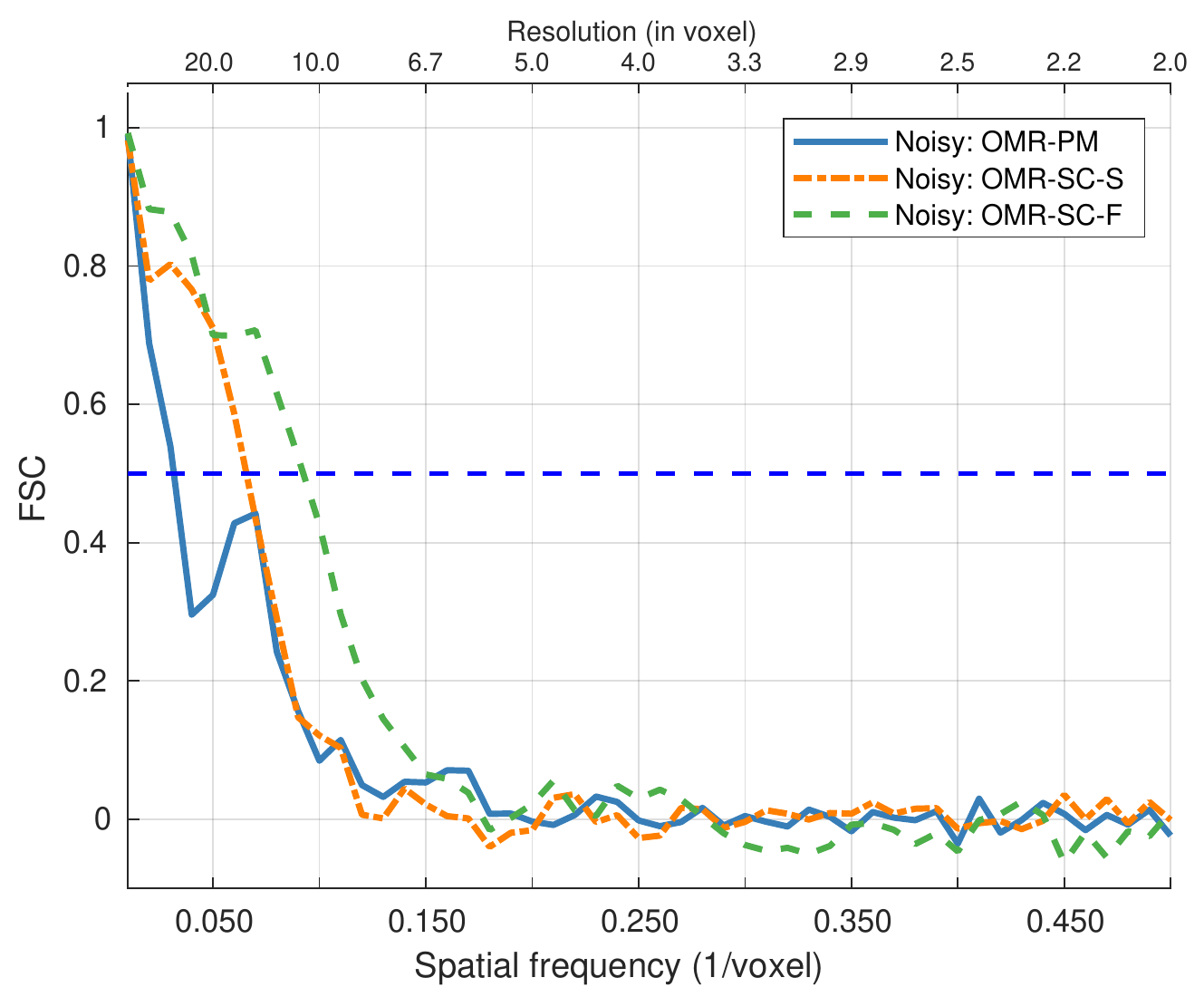}}
\caption{Density 2 (D2): FSC curves of recovered random density maps using the OMR-PM and OMR-SC approaches in the noiseless case and the noisy case (SNR=$0.1$). The cutoff-threshold of ``$0.5$'' is used to determine the resolution (in voxel). }
\label{fig:compare_fsc_d2}
\end{figure*}

\begin{figure*}[tb]
\centering
\subfloat[Noiseless recovery]{
\label{fig:fsc_noiseless_d3}
\includegraphics[width=0.45\textwidth]{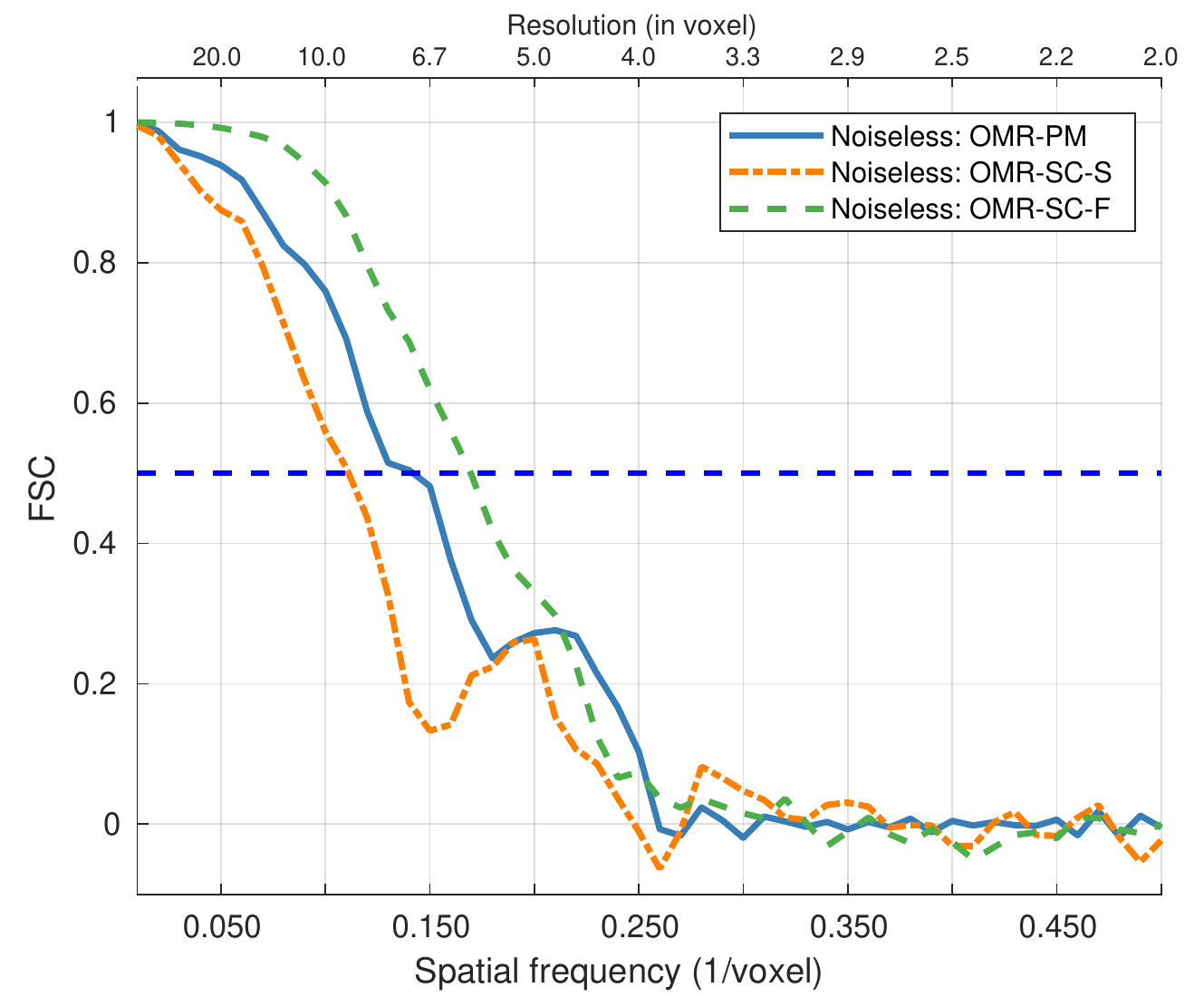}}
\hspace{0.5cm}
\subfloat[Noisy recovery]{
\label{fig:fsc_noisy_d3}
\includegraphics[width=0.45\textwidth]{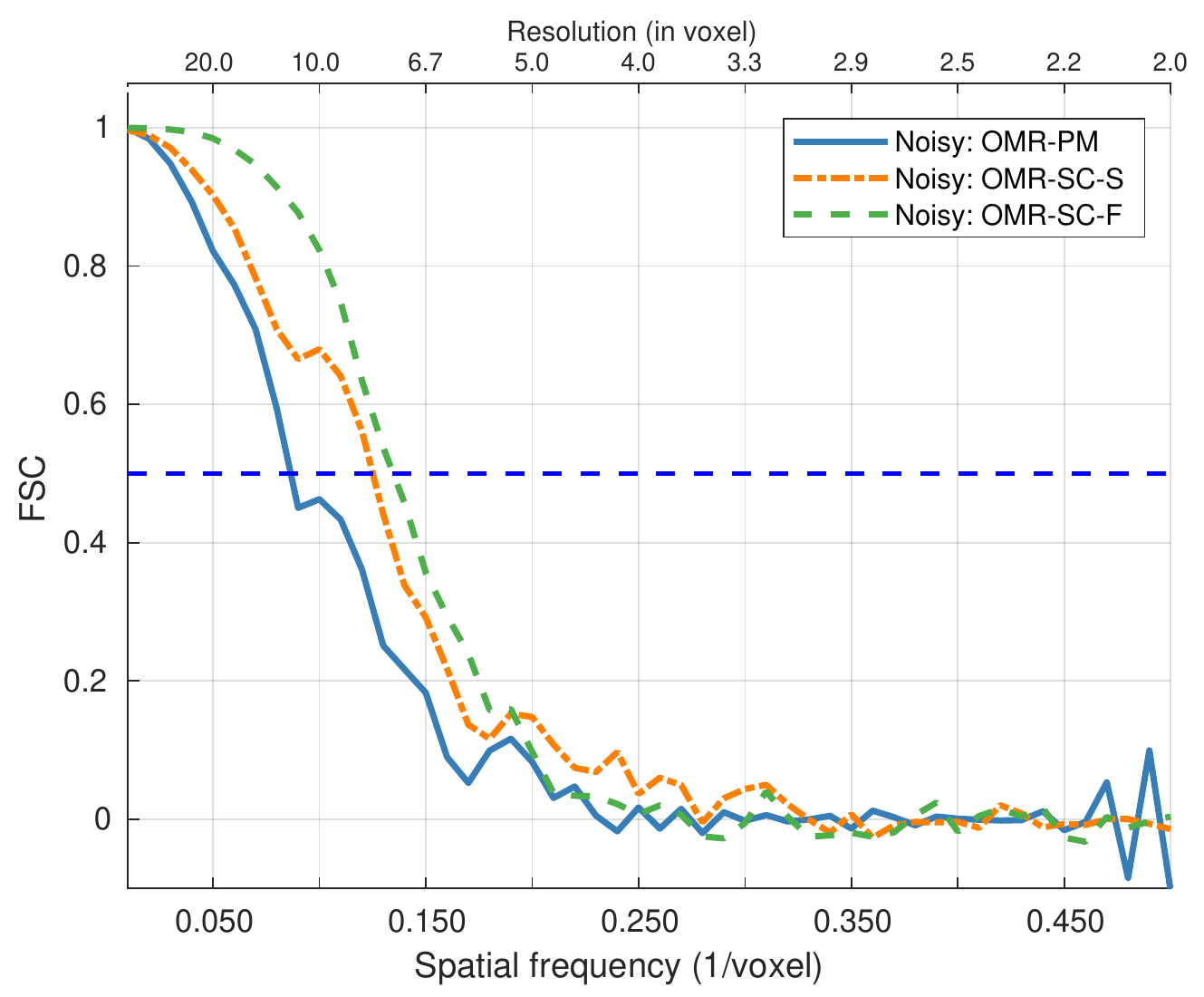}}
\caption{Density 3 (D3): FSC curves of recovered random density maps using the OMR-PM and OMR-SC approaches in the noiseless case and the noisy case (SNR=$0.1$). The cutoff-threshold of ``$0.5$'' is used to determine the resolution (in voxel). }
\label{fig:compare_fsc_d3}
\end{figure*}

\begin{figure*}[tb]
\centering
\subfloat[Noiseless recovery]{
\label{fig:fsc_noiseless_d4}
\includegraphics[width=0.45\textwidth]{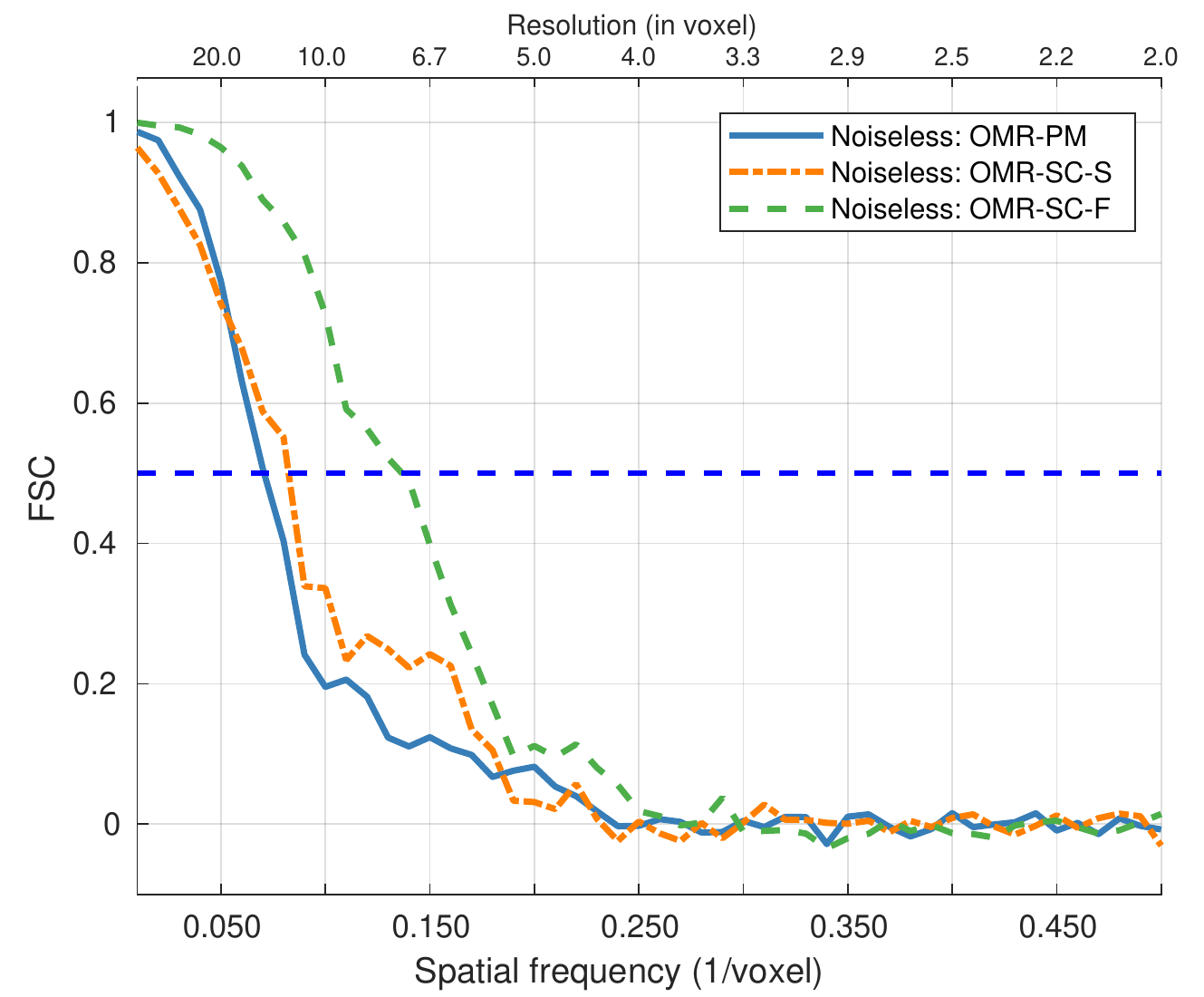}}
\hspace{0.5cm}
\subfloat[Noisy recovery]{
\label{fig:fsc_noisy_d4}
\includegraphics[width=0.45\textwidth]{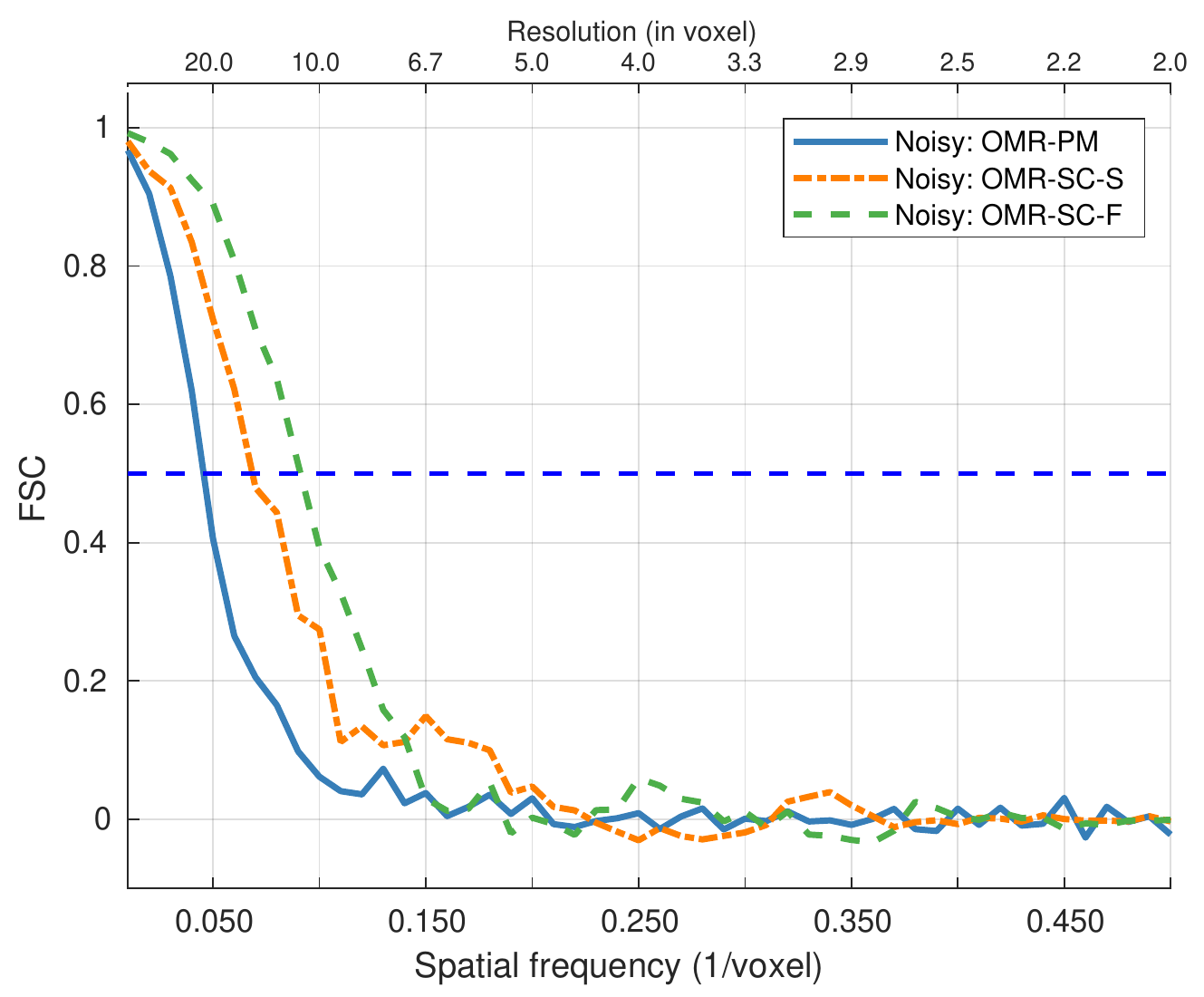}}
\caption{Density 4 (D4): FSC curves of recovered random density maps using the OMR-PM and OMR-SC approaches in the noiseless case and the noisy case (SNR=$0.1$). The cutoff-threshold of ``$0.5$'' is used to determine the resolution (in voxel). }
\label{fig:compare_fsc_d4}
\end{figure*}

\begin{figure*}[tb]
\centering
\subfloat[Noiseless recovery]{
\label{fig:fsc_noiseless_d5}
\includegraphics[width=0.45\textwidth]{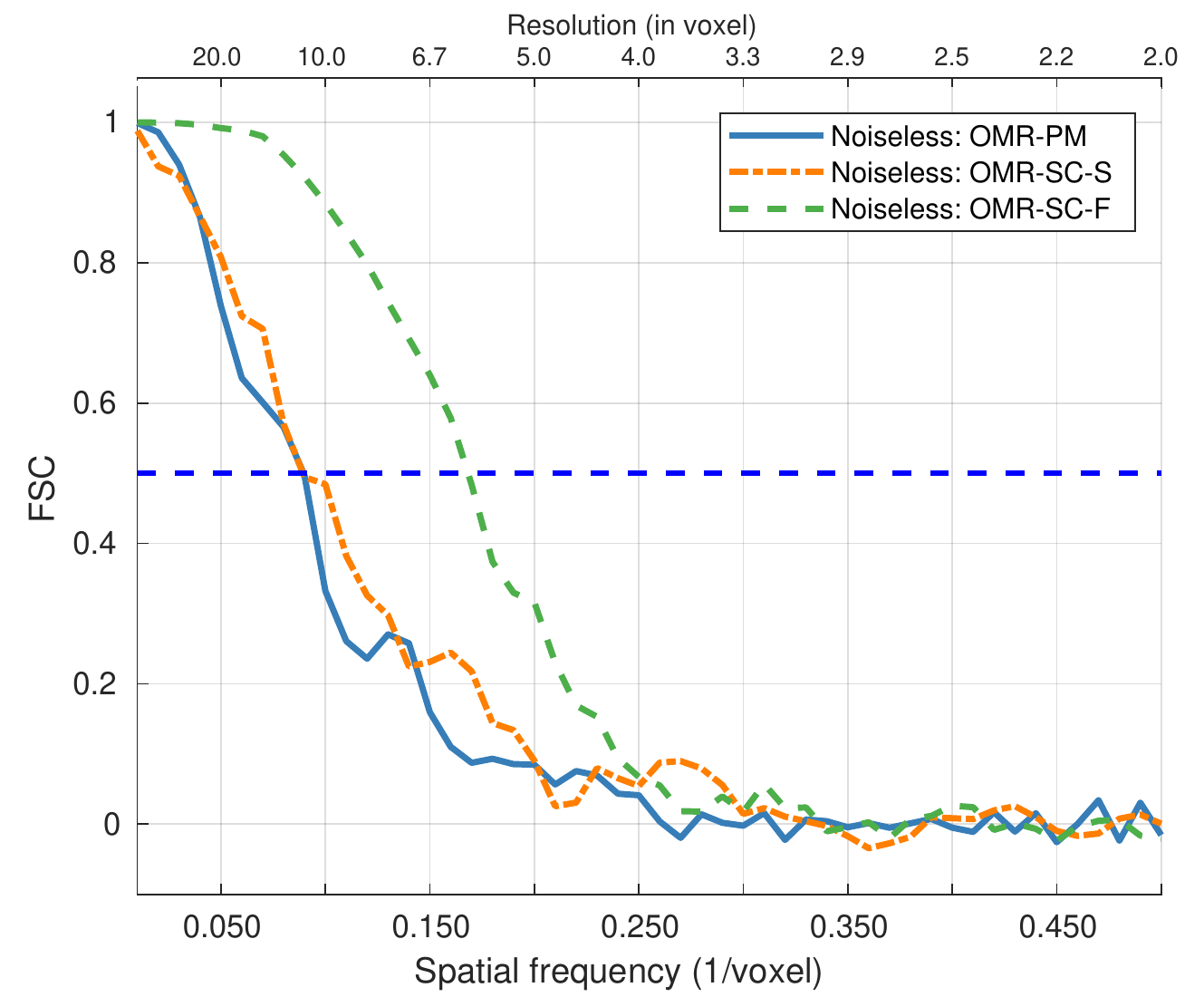}}
\hspace{0.5cm}
\subfloat[Noisy recovery]{
\label{fig:fsc_noisy_d5}
\includegraphics[width=0.45\textwidth]{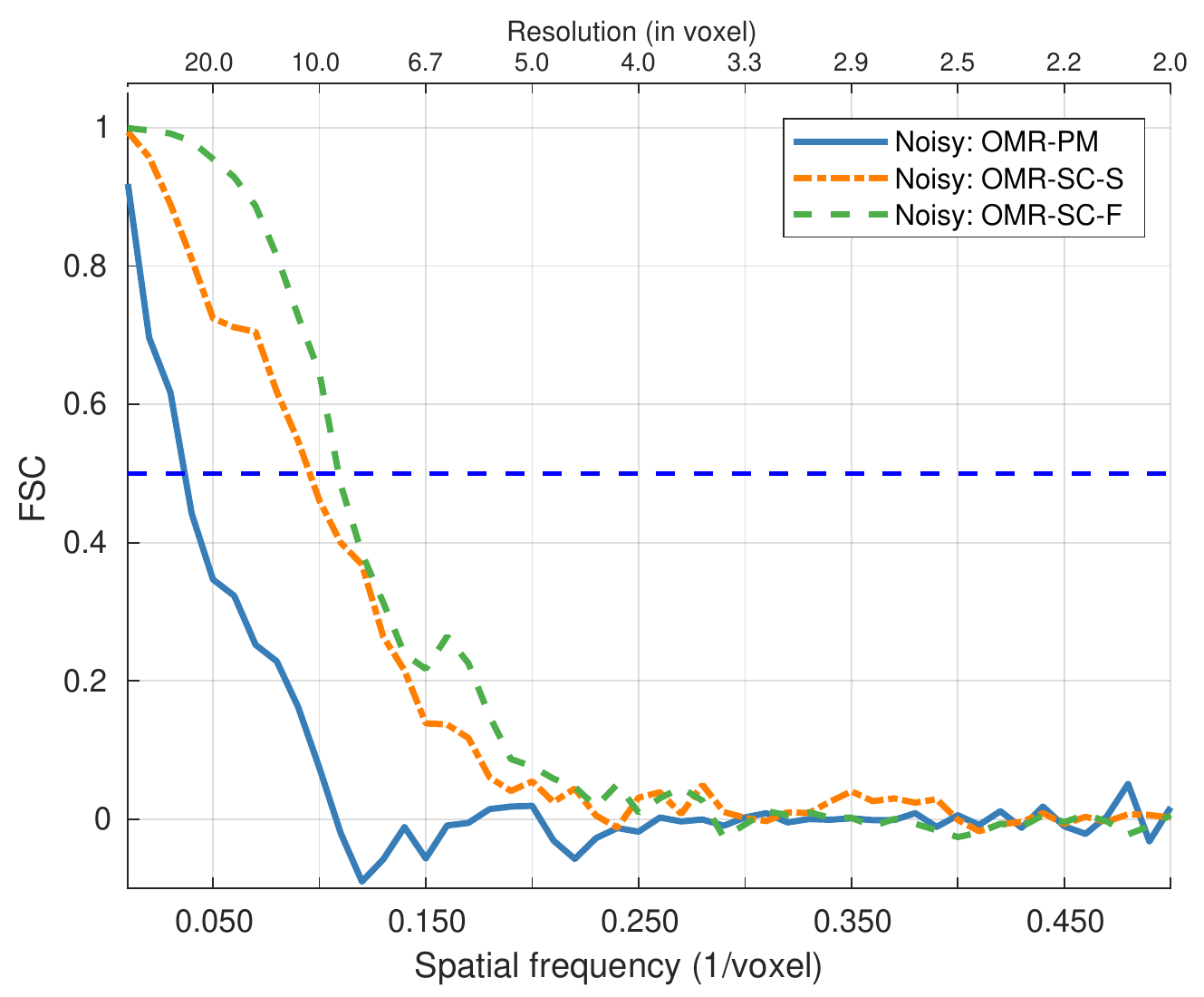}}
\caption{Density 5 (D5): FSC curves of recovered random density maps using the OMR-PM and OMR-SC approaches in the noiseless case and the noisy case (SNR=$0.1$). The cutoff-threshold of ``$0.5$'' is used to determine the resolution (in voxel). }
\label{fig:compare_fsc_d5}
\end{figure*}

\begin{figure*}[tb]
\centering
\subfloat[Noiseless recovery]{
\label{fig:fsc_noiseless_d6}
\includegraphics[width=0.45\textwidth]{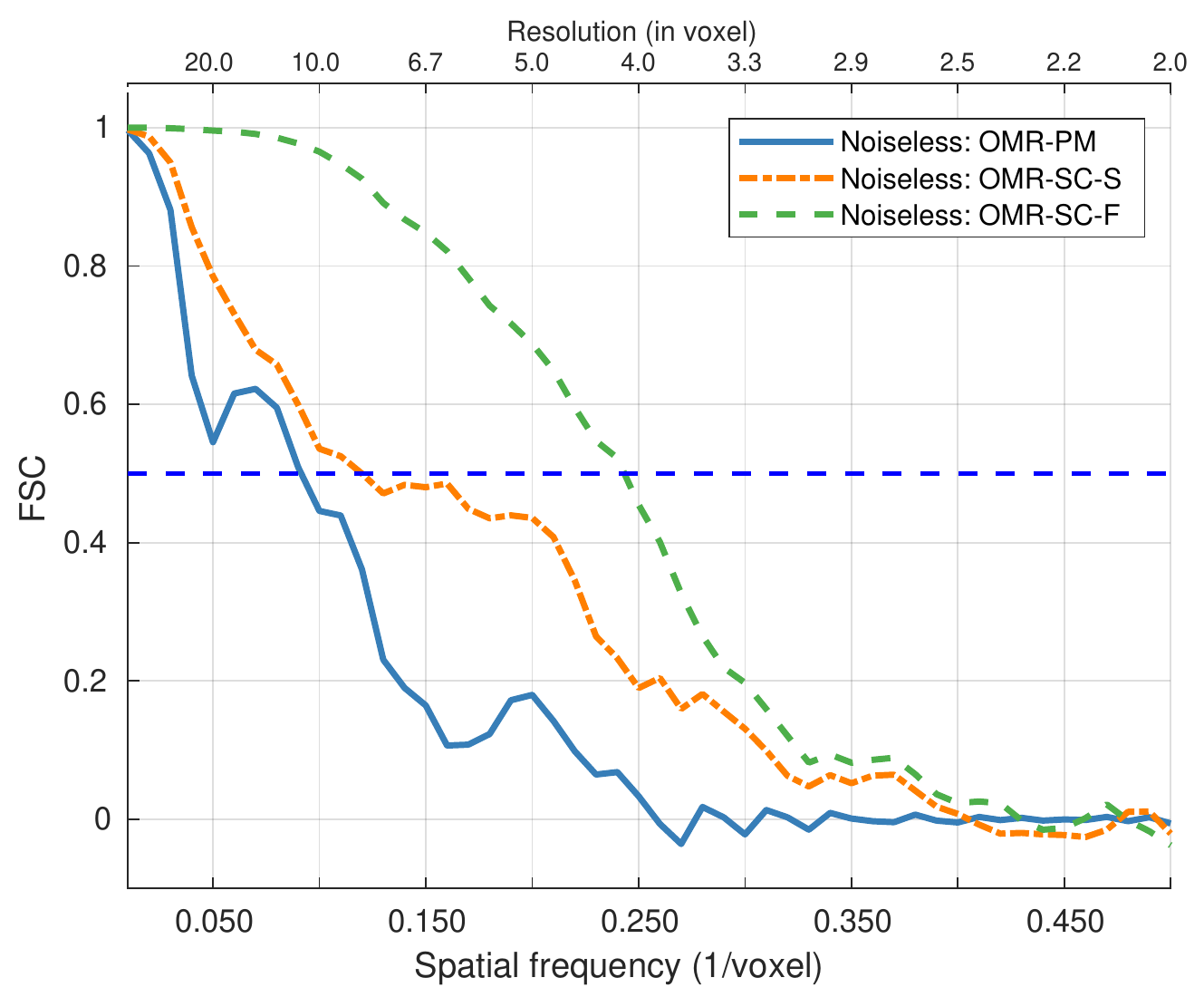}}
\hspace{0.5cm}
\subfloat[Noisy recovery]{
\label{fig:fsc_noisy_d6}
\includegraphics[width=0.45\textwidth]{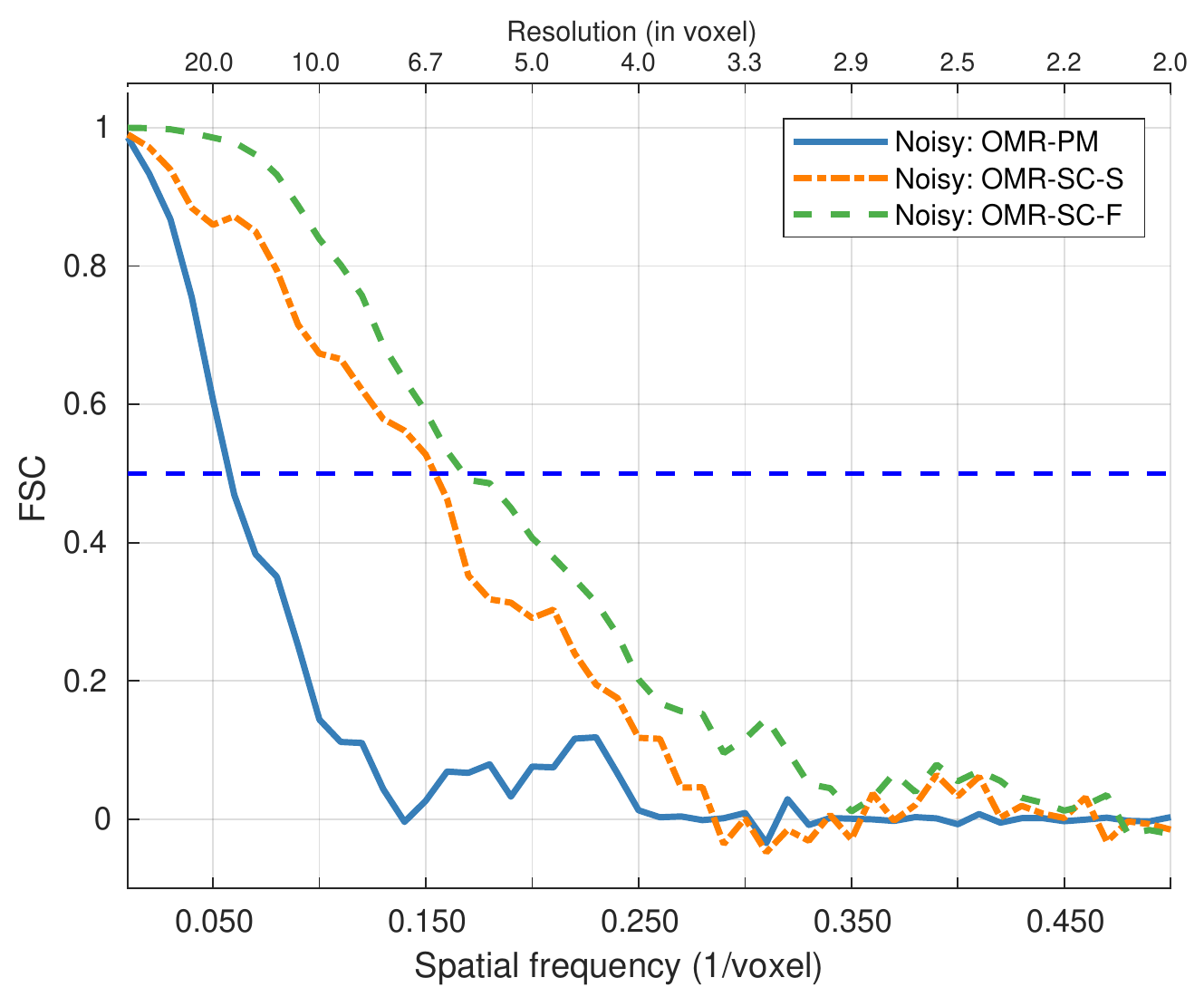}}
\caption{Density 6 (D6): FSC curves of recovered random density maps using the OMR-PM and OMR-SC approaches in the noiseless case and the noisy case (SNR=$0.1$). The cutoff-threshold of ``$0.5$'' is used to determine the resolution (in voxel). }
\label{fig:compare_fsc_d6}
\end{figure*}

\begin{figure*}[tb]
\centering
\subfloat[Noiseless recovery]{
\label{fig:fsc_noiseless_d7}
\includegraphics[width=0.45\textwidth]{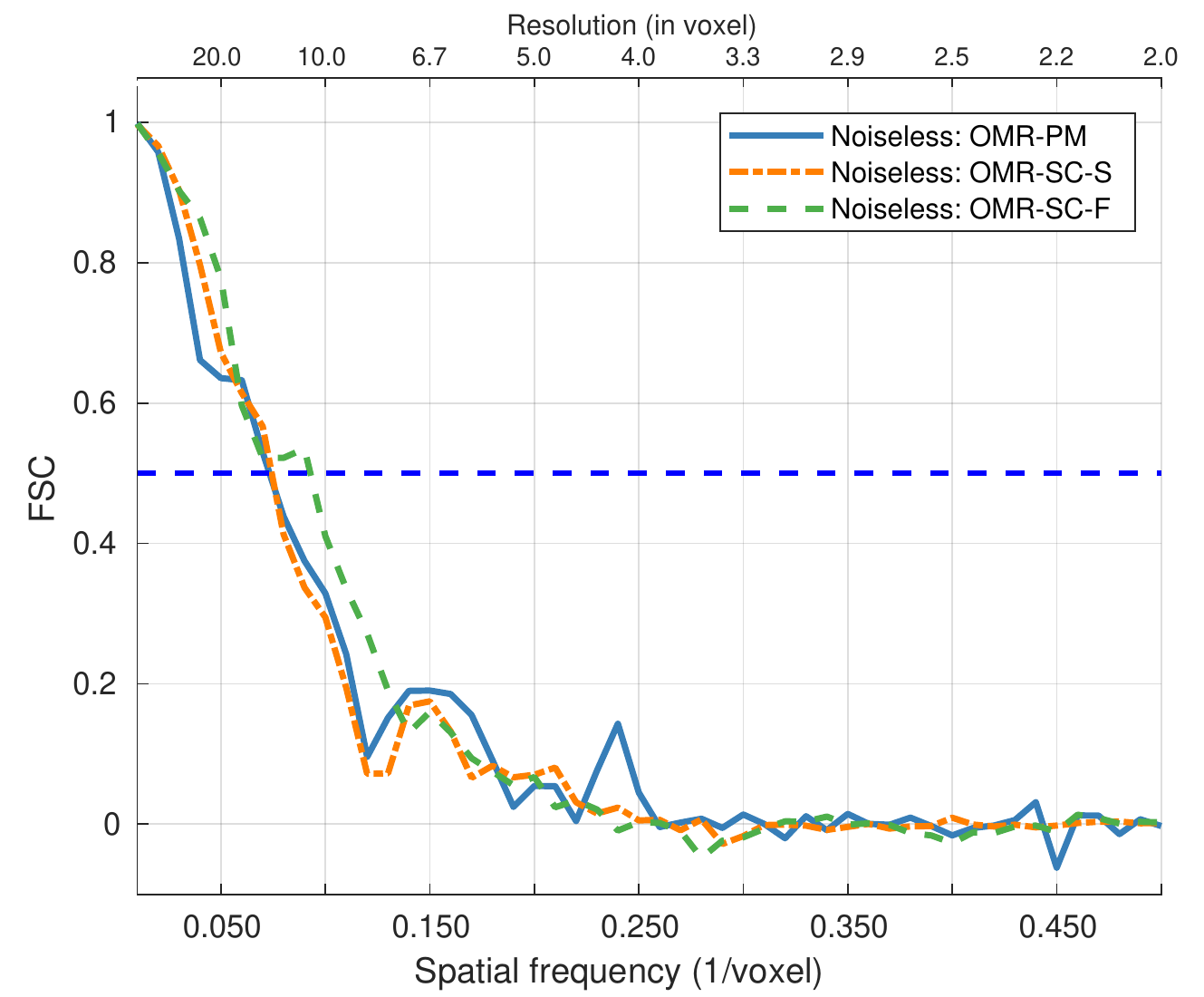}}
\hspace{0.5cm}
\subfloat[Noisy recovery]{
\label{fig:fsc_noisy_d7}
\includegraphics[width=0.45\textwidth]{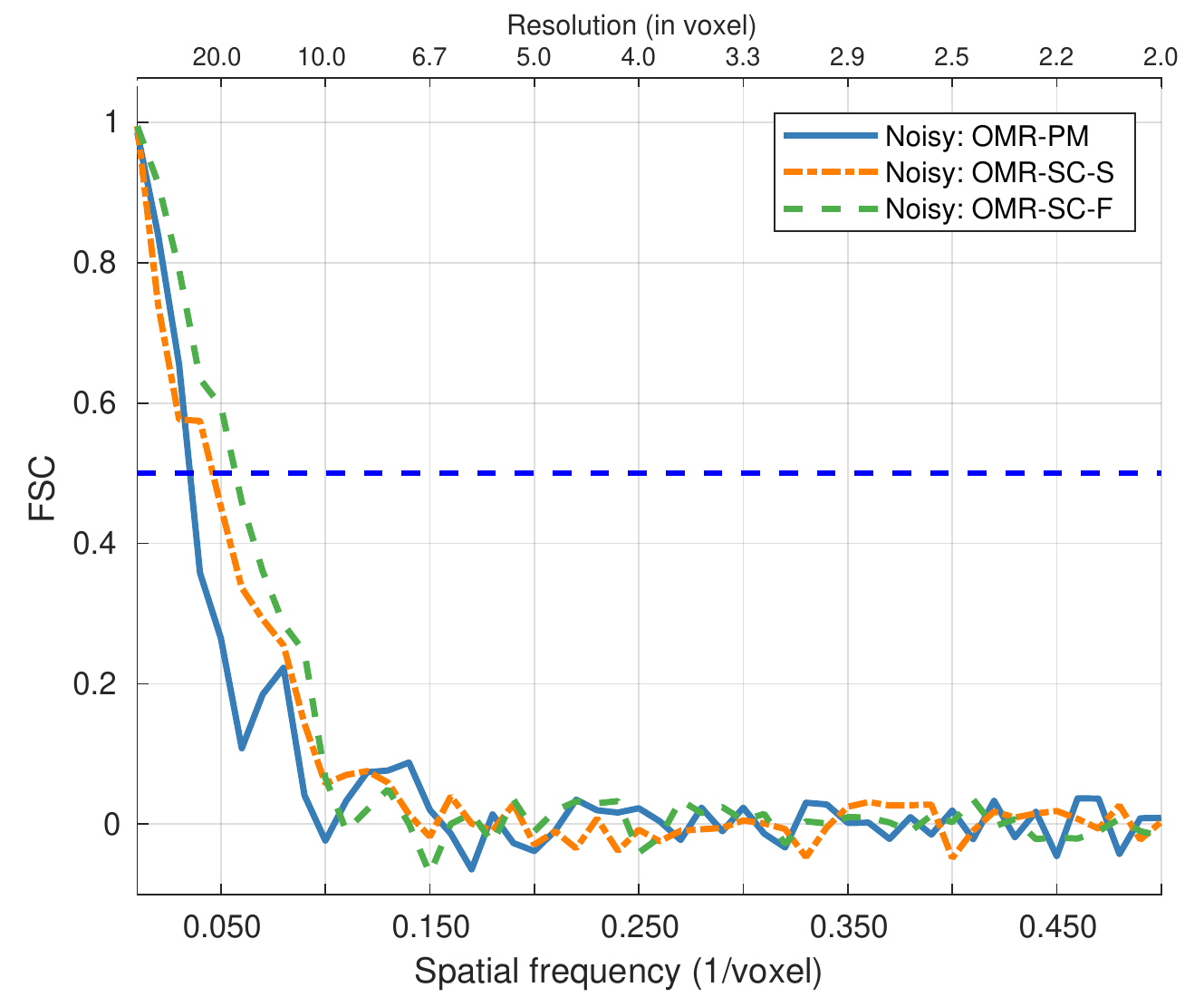}}
\caption{Density 7 (D7): FSC curves of recovered random density maps using the OMR-PM and OMR-SC approaches in the noiseless case and the noisy case (SNR=$0.1$). The cutoff-threshold of ``$0.5$'' is used to determine the resolution (in voxel). }
\label{fig:compare_fsc_d7}
\end{figure*}

\begin{figure*}[tb]
\centering
\subfloat[Noiseless recovery]{
\label{fig:fsc_noiseless_d8}
\includegraphics[width=0.45\textwidth]{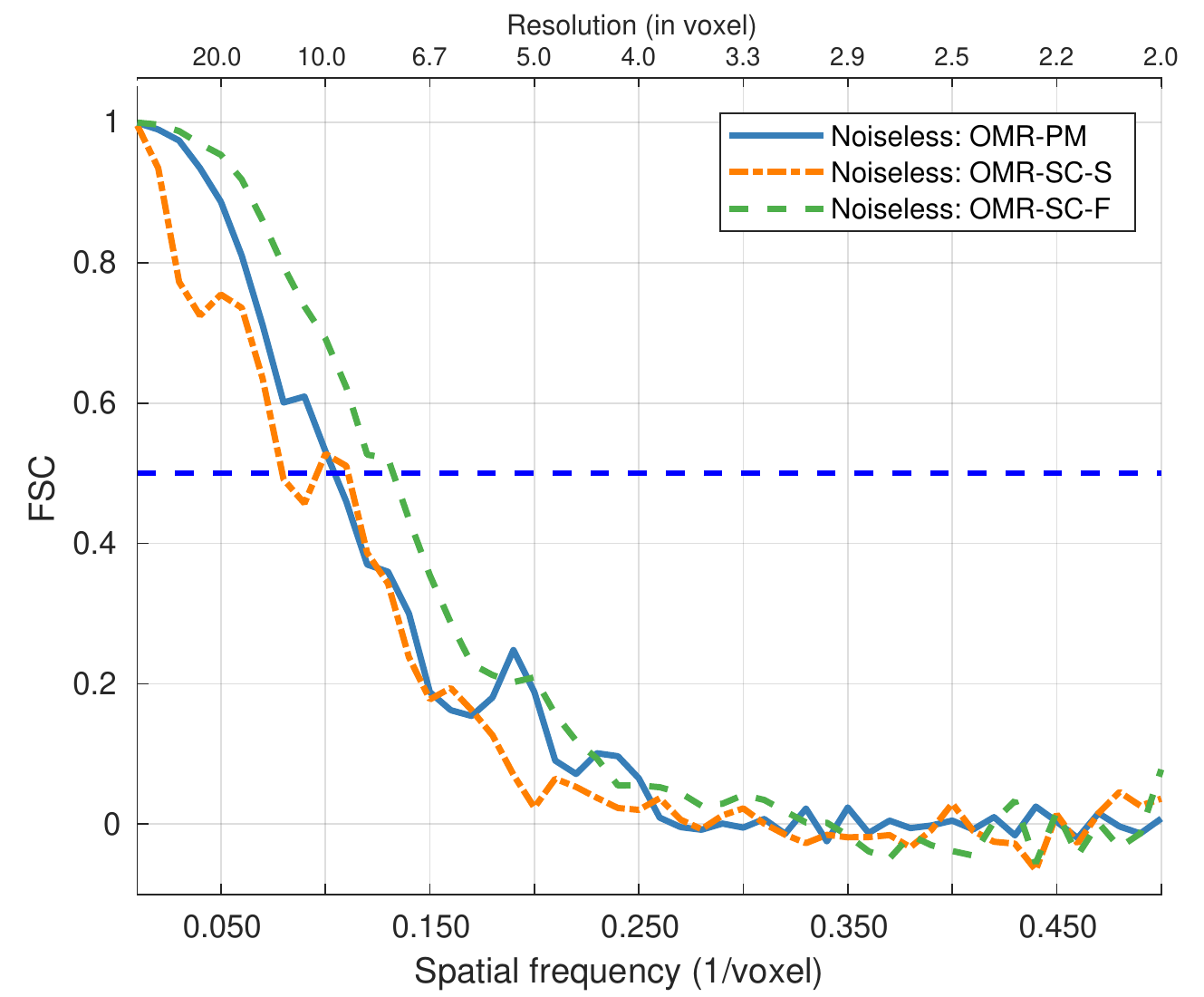}}
\hspace{0.5cm}
\subfloat[Noisy recovery]{
\label{fig:fsc_noisy_d8}
\includegraphics[width=0.45\textwidth]{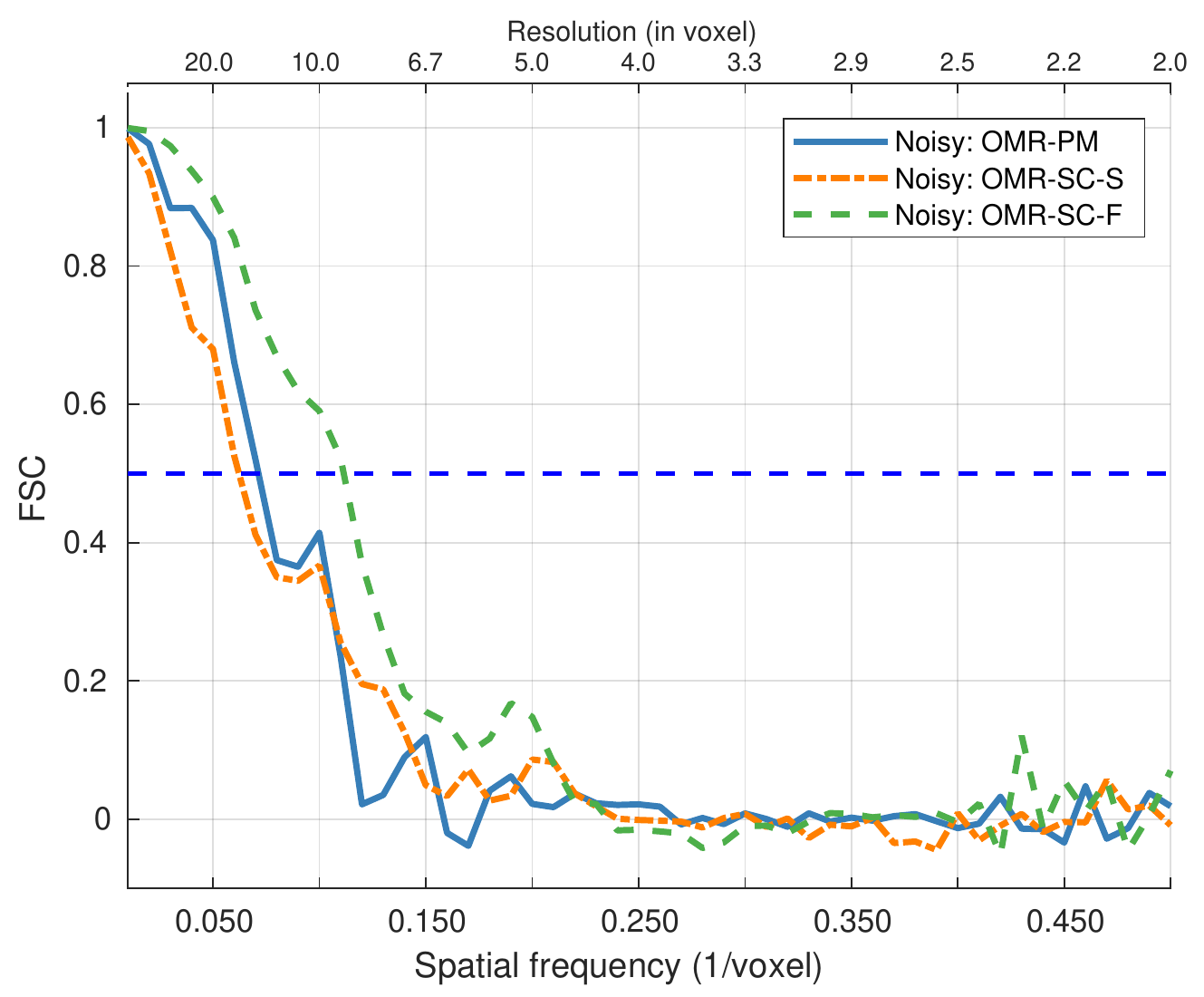}}
\caption{Density 8 (D8): FSC curves of recovered random density maps using the OMR-PM and OMR-SC approaches in the noiseless case and the noisy case (SNR=$0.1$). The cutoff-threshold of ``$0.5$'' is used to determine the resolution (in voxel). }
\label{fig:compare_fsc_d8}
\end{figure*}

\begin{figure*}[tb]
\centering
\subfloat[Noiseless recovery]{
\label{fig:fsc_noiseless_d9}
\includegraphics[width=0.45\textwidth]{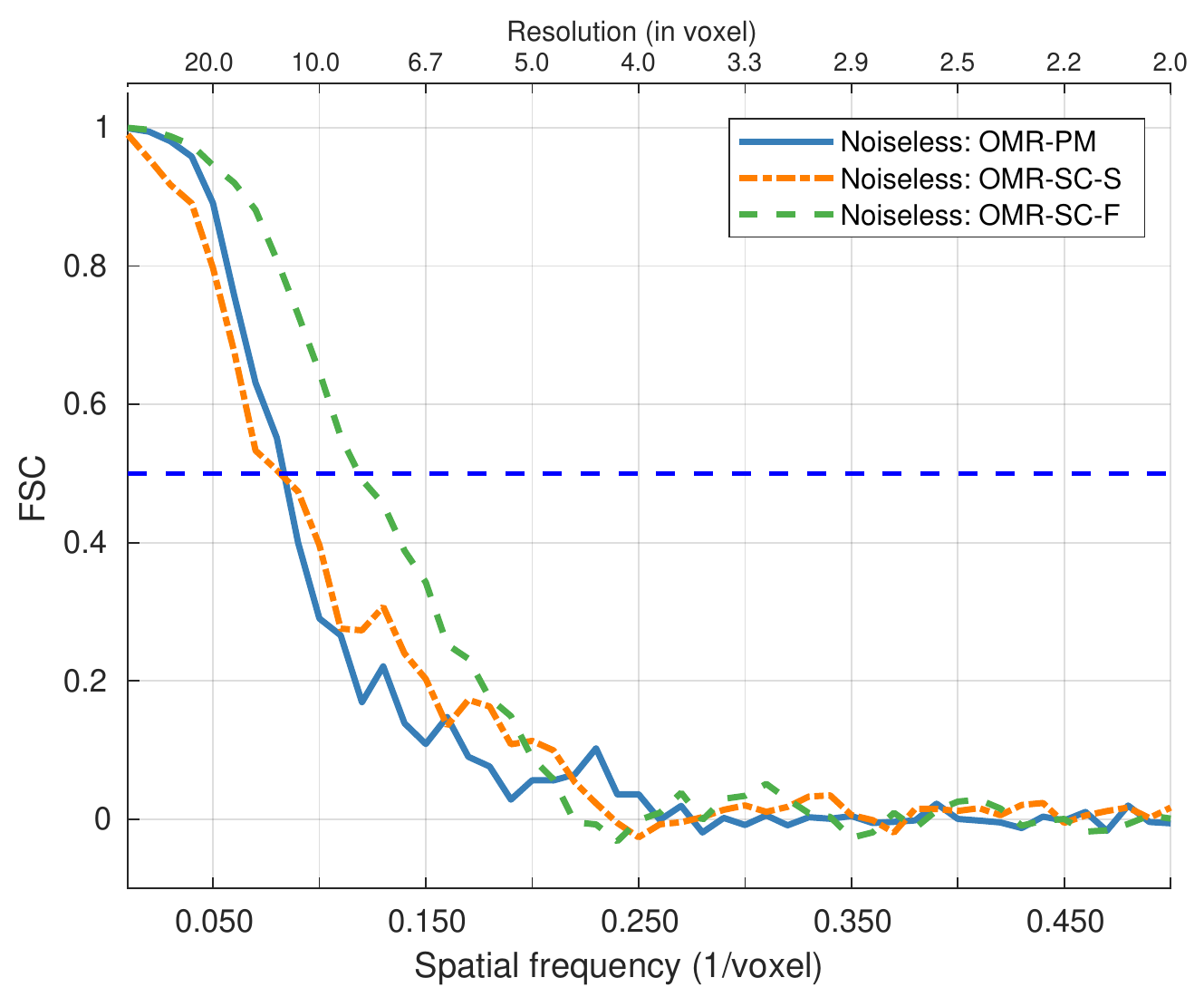}}
\hspace{0.5cm}
\subfloat[Noisy recovery]{
\label{fig:fsc_noisy_d9}
\includegraphics[width=0.45\textwidth]{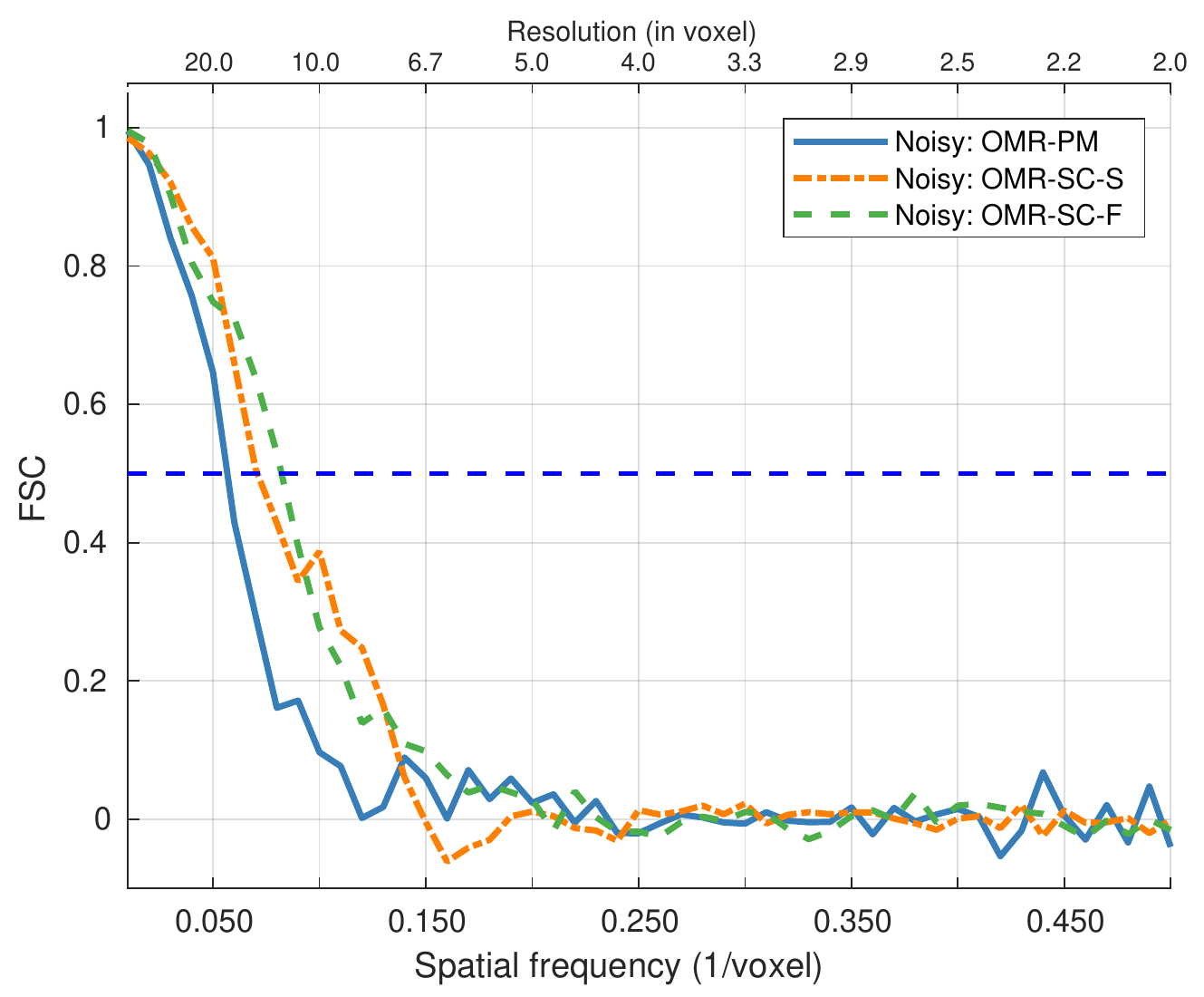}}
\caption{Density 9 (D9): FSC curves of recovered random density maps using the OMR-PM and OMR-SC approaches in the noiseless case and the noisy case (SNR=$0.1$). The cutoff-threshold of ``$0.5$'' is used to determine the resolution (in voxel). }
\label{fig:compare_fsc_d9}
\end{figure*}

\begin{figure*}[p]
\centering
\includegraphics[width=\textwidth]{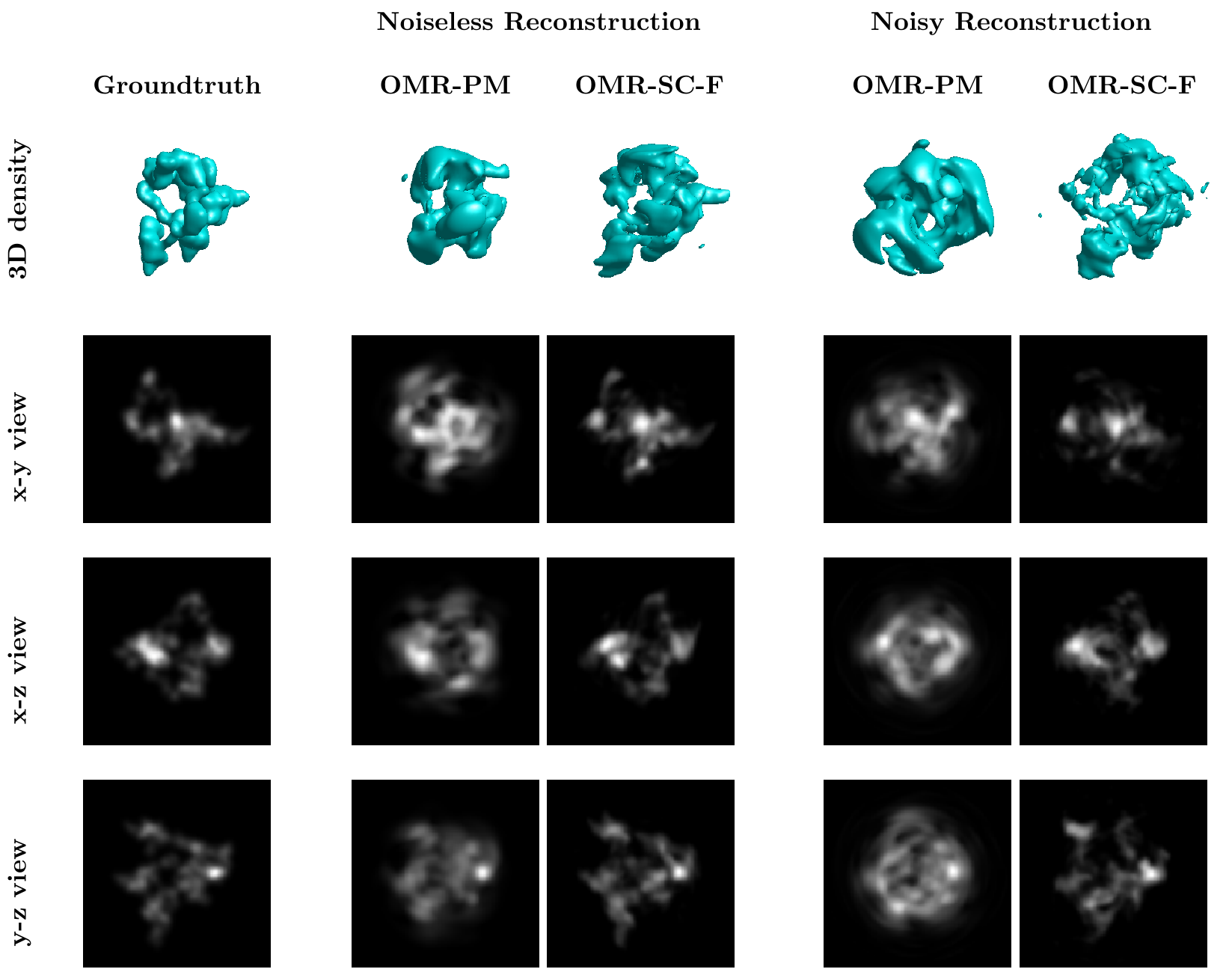}
\caption{Density 2 (D2): reconstructions using the OMR-PM and OMR-SC-F approaches in the noiseless case and the noisy case (SNR=$0.1$).}
\label{fig:compare_reconstruction_2th}
\end{figure*}

\begin{figure*}[p]
\centering
\includegraphics[width=\textwidth]{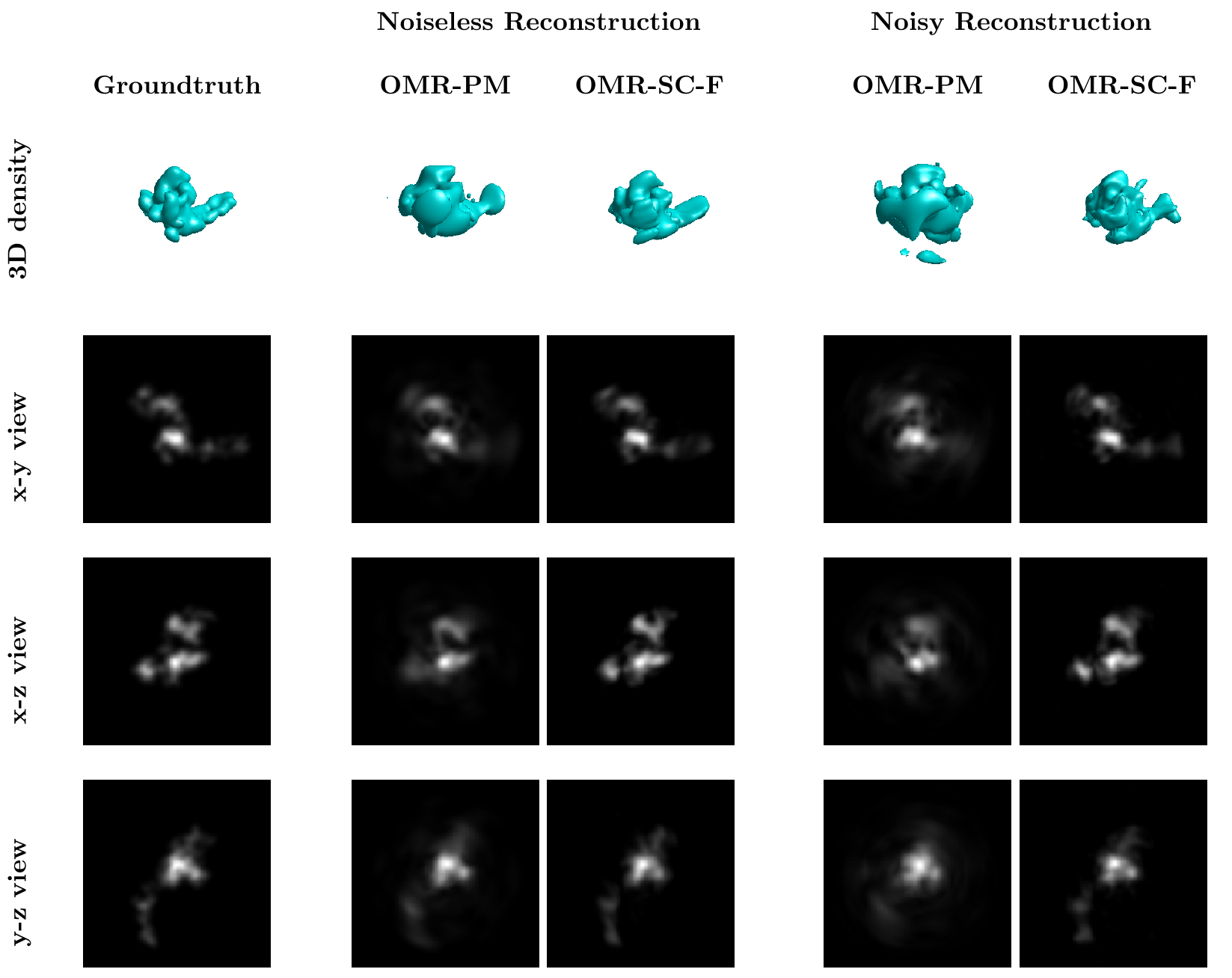}
\caption{Density 3 (D3): reconstructions using the OMR-PM and OMR-SC-F approaches in the noiseless case and the noisy case (SNR=$0.1$).}
\label{fig:compare_reconstruction_3rd}
\end{figure*}

\begin{figure*}[p]
\centering
\includegraphics[width=\textwidth]{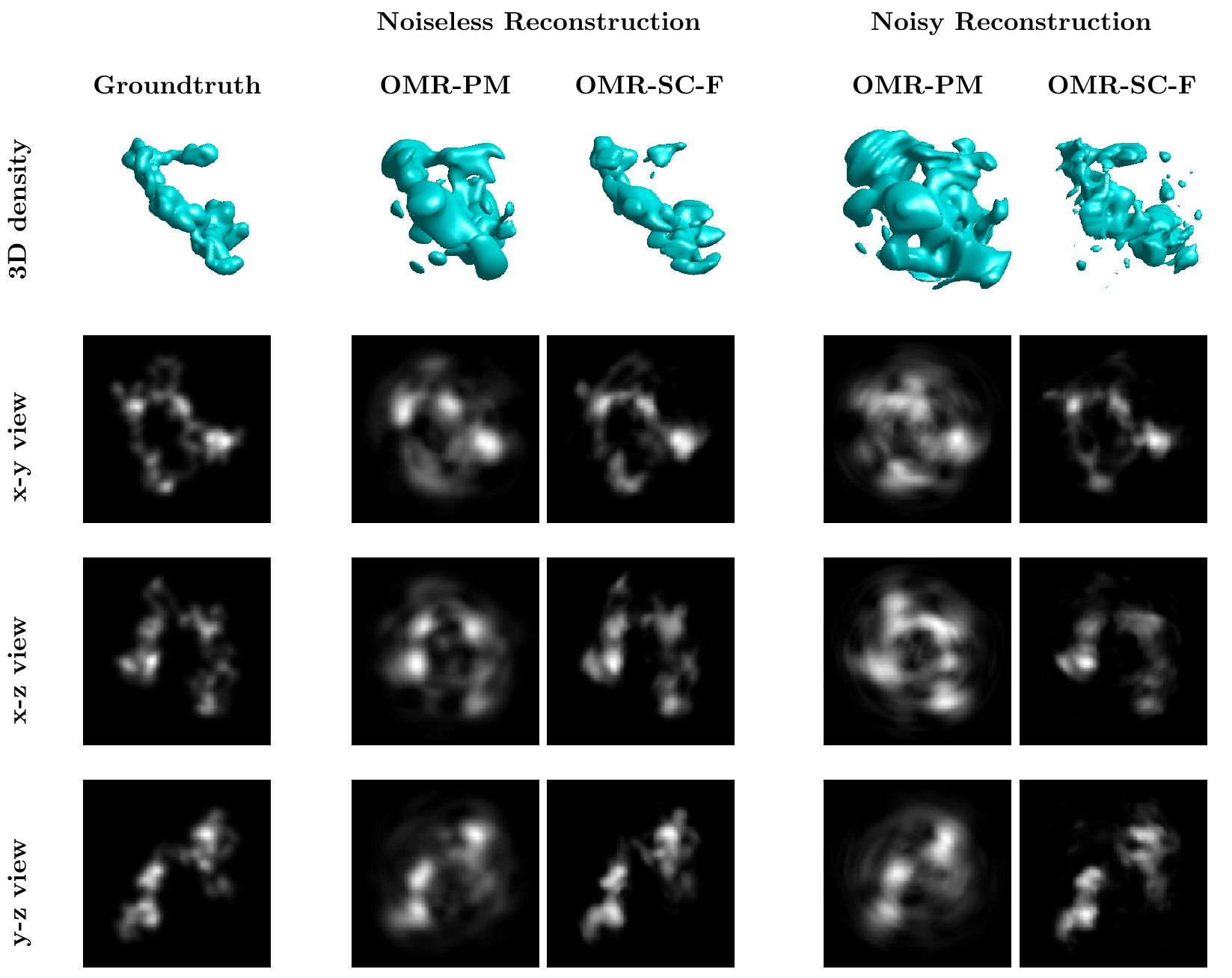}
\caption{Density 4 (D4): reconstructions using the OMR-PM and OMR-SC-F approaches in the noiseless case and the noisy case (SNR=$0.1$).}
\label{fig:compare_reconstruction_4th}
\end{figure*}

\begin{figure*}[p]
\centering
\includegraphics[width=\textwidth]{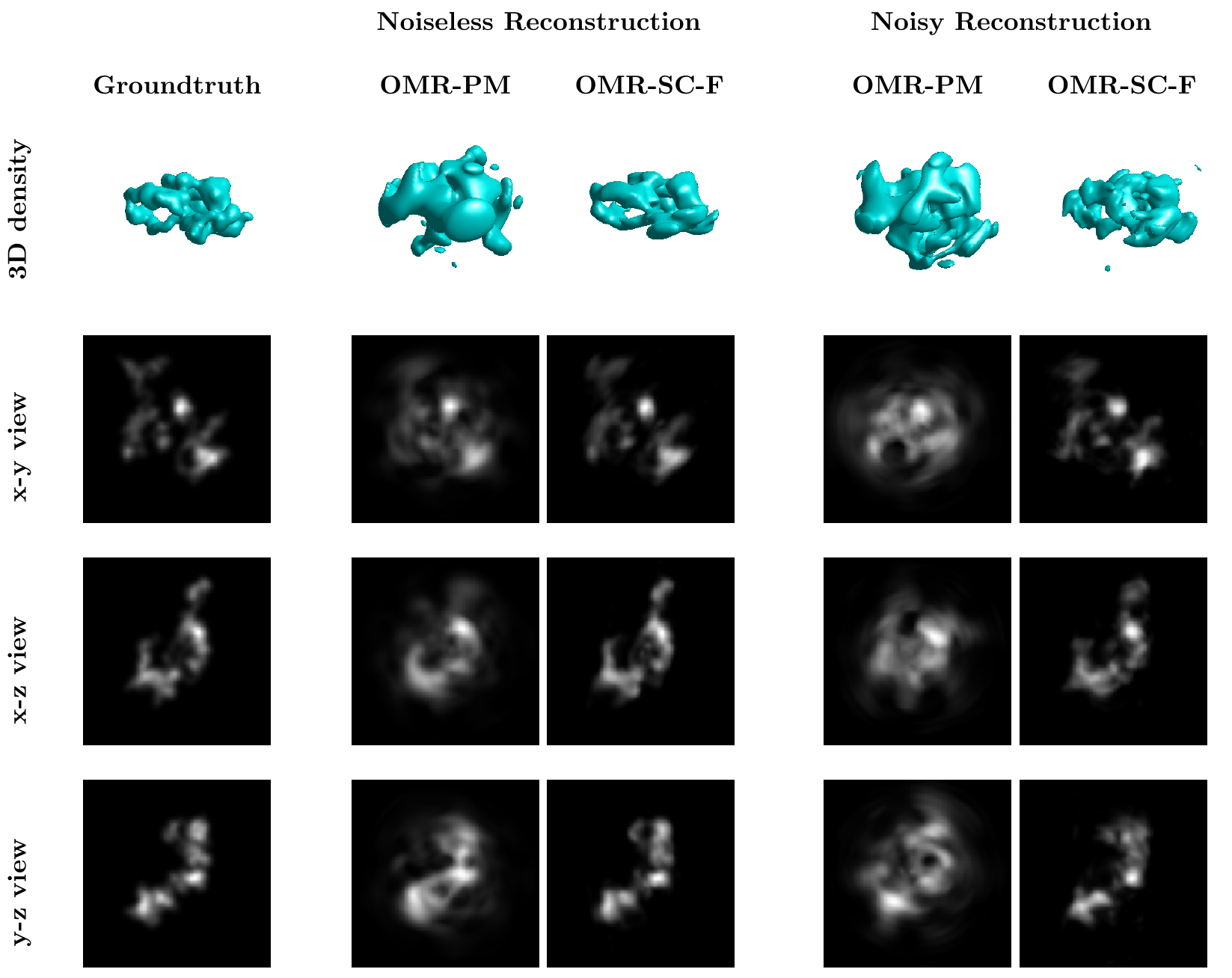}
\caption{Density 5 (D5): reconstructions using the OMR-PM and OMR-SC-F approaches in the noiseless case and the noisy case (SNR=$0.1$).}
\label{fig:compare_reconstruction_5th}
\end{figure*}

\begin{figure*}[p]
\centering
\includegraphics[width=\textwidth]{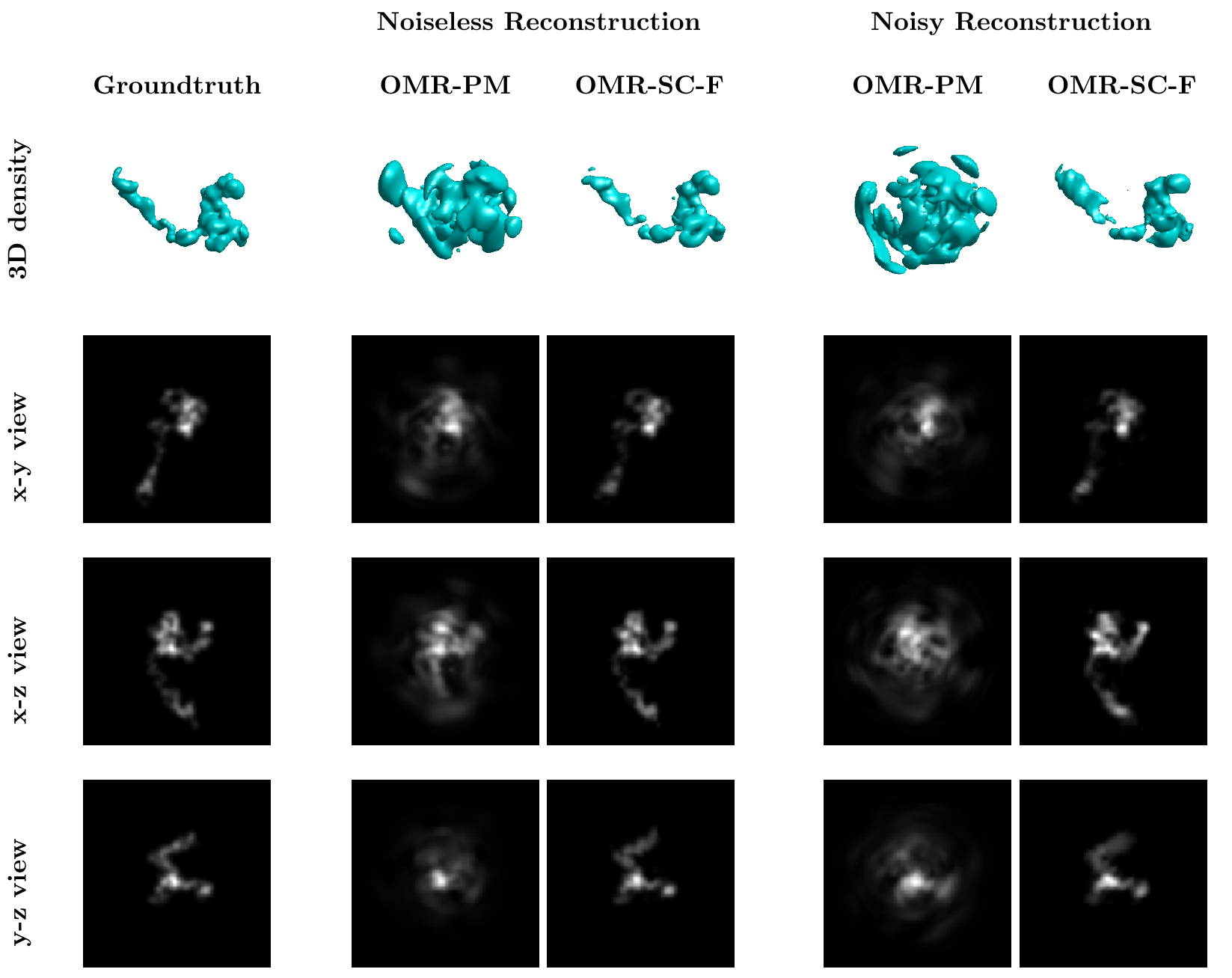}
\caption{Density 6 (D6): reconstructions using the OMR-PM and OMR-SC-F approaches in the noiseless case and the noisy case (SNR=$0.1$).}
\label{fig:compare_reconstruction_6th}
\end{figure*}

\begin{figure*}[p]
\centering
\includegraphics[width=\textwidth]{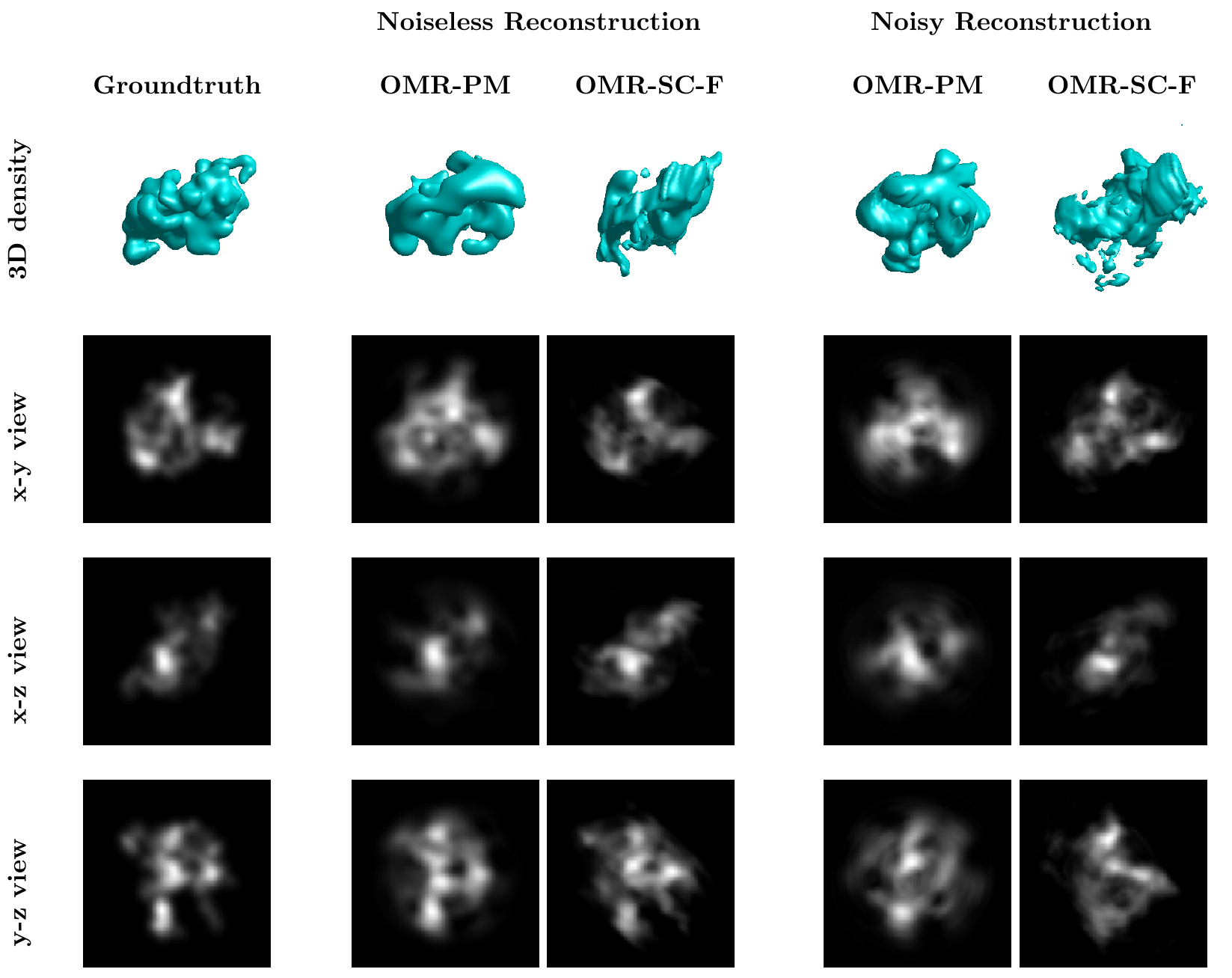}
\caption{Density 7 (D7): reconstructions using the OMR-PM and OMR-SC-F approaches in the noiseless case and the noisy case (SNR=$0.1$).}
\label{fig:compare_reconstruction_7th}
\end{figure*}

\begin{figure*}[p]
\centering
\includegraphics[width=\textwidth]{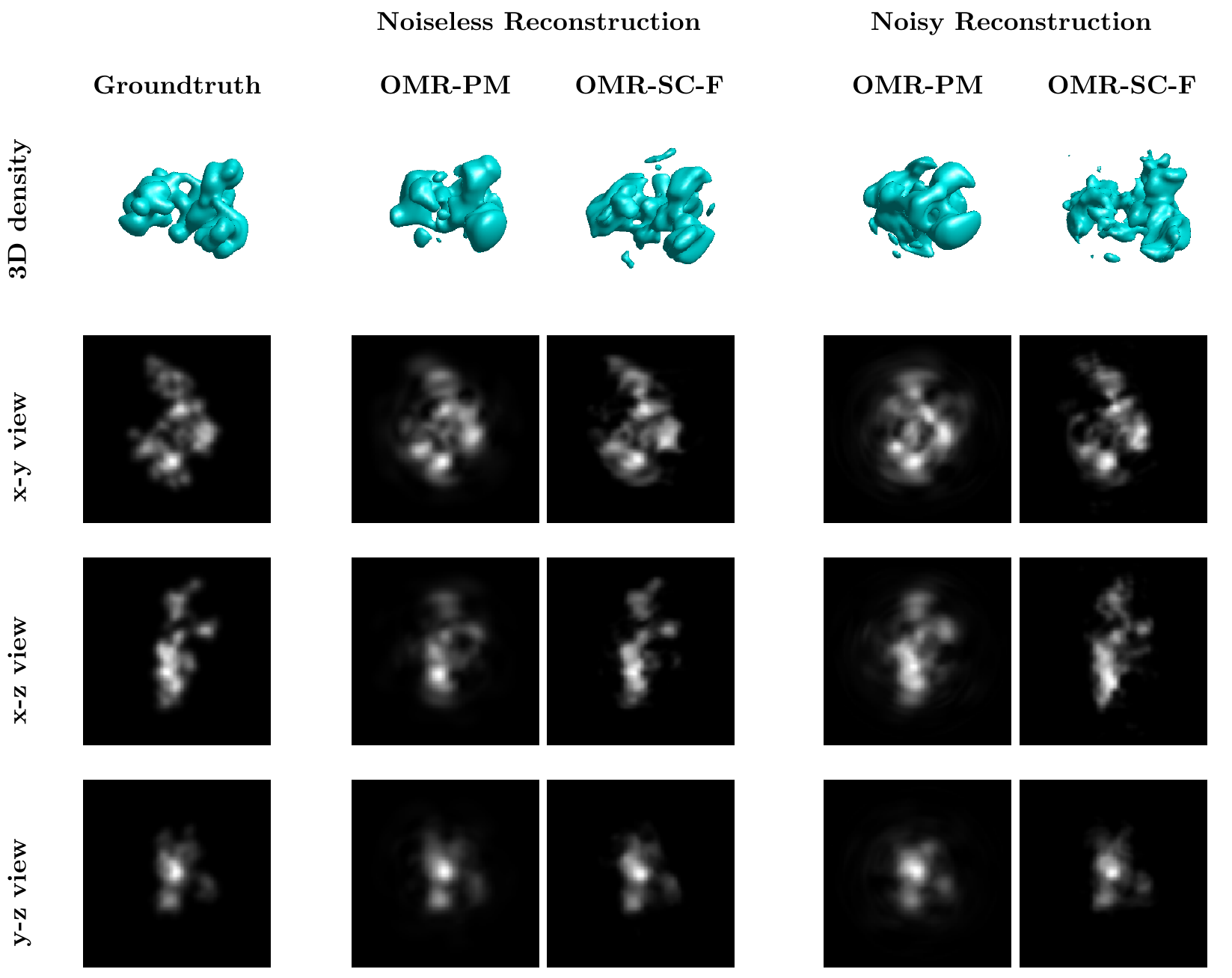}
\caption{Density 8 (D8): reconstructions using the OMR-PM and OMR-SC-F approaches in the noiseless case and the noisy case (SNR=$0.1$).}
\label{fig:compare_reconstruction_8th}
\end{figure*}

\begin{figure*}[p]
\centering
\includegraphics[width=\textwidth]{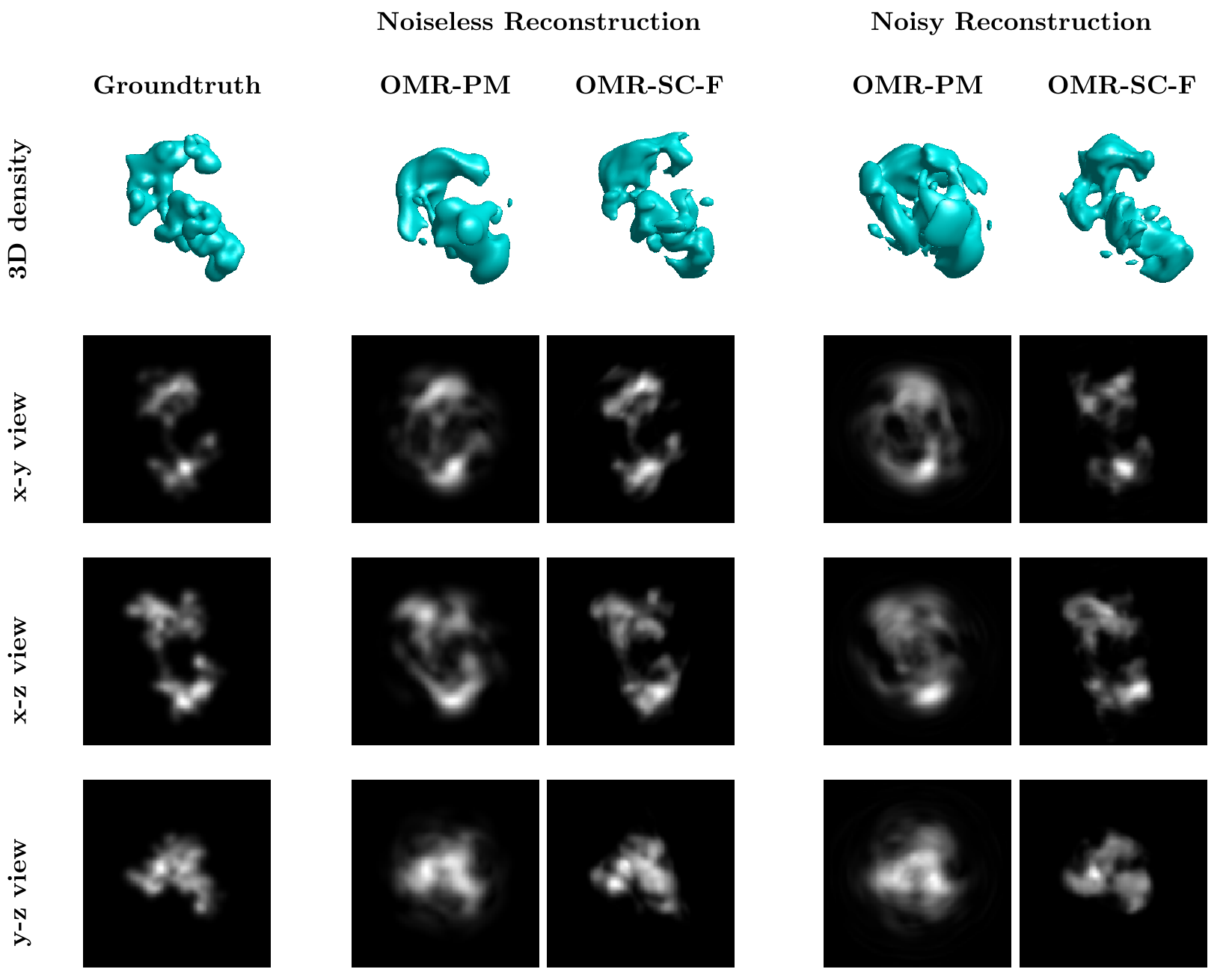}
\caption{Density 9 (D9): reconstructions using the OMR-PM and OMR-SC-F approaches in the noiseless case and the noisy case (SNR=$0.1$).}
\label{fig:compare_reconstruction_9th}
\end{figure*}

\end{document}

%% file: ex_shared.tex
\usepackage{lipsum}
\usepackage{amsfonts}
\usepackage{graphicx}
\usepackage{epstopdf}
\usepackage{algorithmic}
\ifpdf
  \DeclareGraphicsExtensions{.eps,.pdf,.png,.jpg}
\else
  \DeclareGraphicsExtensions{.eps}
\fi

\newsiamremark{remark}{Remark}
\newsiamremark{hypothesis}{Hypothesis}
\crefname{hypothesis}{Hypothesis}{Hypotheses}
\newsiamthm{claim}{Claim}

\headers{Orthogonal Matrix Retrieval with Spatial Consensus}{Shuai Huang, Mona Zehni, Ivan Dokmani\'c, and Zhizhen Zhao}

\title{Orthogonal Matrix Retrieval with Spatial Consensus\\ for 3D Unknown-View Tomography\thanks{Submitted to the editors on March 8th, 2023.
\funding{Ivan Dokmani\'c was supported by the European Research Council (ERC) Starting Grant 852821$-$SWING.  Mona Zehni and Zhizhen Zhao were partly supported by NSF DMS-1854791, NSF OAC-1934757, and Alfred P. Sloan Foundation.}}}

\author{Shuai Huang\thanks{Department of Radiology and Imaging Sciences, Emory University, Atlanta, GA 30329 USA
  (\email{shuai.huang@emory.edu}).}
\and Mona Zehni\thanks{Coordinate Science Laboratory, University of Illinois at Urbana–Champaign, Urbana, IL 61801 USA 
  (\email{mzehni2@illinois.edu}, \email{zhizhenz@illinois.edu}).}
\and Ivan Dokmani\'c\thanks{Department of Mathematics and Computer Science, University of Basel, 4001 Basel, Switzerland
  (\email{ivan.dokmanic@unibas.ch}).}
\and Zhizhen Zhao\footnotemark[3]}

\usepackage{amsopn}
